\documentclass{amsart}
\input BoxedEPS																	
\SetepsfEPSFSpecial									
\HideDisplacementBoxes										

\newtheorem{theorem}{Theorem}[section]
\newtheorem{lemma}[theorem]{Lemma}

\theoremstyle{definition}

\theoremstyle{remark}

\numberwithin{equation}{section}

\newcommand{\HHH}{{\mathbb H}^3}    
\newcommand{\HH}{{\mathbb H}^2}     
\newcommand{\Rr}{{\mathbb R}}       
\newcommand{\NOE}{N^\circ_\epsilon} 

\begin{document}

\title{Injectivity Radius Bounds in Hyperbolic $I$-Bundle Convex Cores}

\author{Carol E.\ Fan}
\address{Department of Mathematics, Oklahoma State University, Stillwater,
Oklahoma  74078-1058}
\curraddr{Department of Mathematics, Loyola Marymount University, Los
Angeles, California  90045-8130}
\email{cfan@lmumail.lmu.edu}
\thanks{I would like to thank Richard Canary for advising me during this
project.  In addition, Francis Bonahon, Timothy Comar, Robert Myers, Mingqing Ouyang,
Peter Scott, and Edward Taylor were very helpful with comments and
discussions.  The results in this paper represent a portion of my Ph.D.\
thesis completed at The University of Michigan.}


\date{\today}



\begin{abstract}
A version of a conjecture of McMullen is as follows: \ Given a
hyperbolizable 3-manifold M with incompressible boundary, there exists a
uniform constant K such that if N is a hyperbolic 3-manifold homeomorphic to the
interior of M, then the injectivity radius based at points in the convex core
of N is bounded above by K.  This conjecture suggests that convex cores are
uniformly congested.  We will give a proof in the case when M is an I-bundle
over a closed surface, taking into account the possibility of cusps.   
\end{abstract}

\maketitle

\section{Introduction}

A version of a conjecture of McMullen is as follows: \ Given a
hyperbolizable 3-manifold $M$ with incompressible boundary, there exists a
uniform constant $K$ such that if $N$ is a hyperbolic 3-manifold homeomorphic to
the interior of $M$, then the injectivity radius based at points in the convex
core of $N$ is bounded above by $K$.  This conjecture suggests that convex cores
are uniformly congested.  We will give a proof in the case when $M$ is an
$I$-bundle over a closed surface, taking into account the possibility of
cusps. 

In this paper, we investigate the geometry of convex cores of
hyperbolic 3-manifolds which are homotopy equivalent to a closed surface. 
If $M$ is an $I$-bundle over a closed surface, then
there exists a uniform upper bound on injectivity radius for points in the
convex core of any hyperbolic 3-manifold homeomorphic to the interior of
$M$.

The main result is an extension of a theorem of Kerckhoff-Thurston
\cite{kerck/thur}.  Their theorem established the existence of an upper bound
on injectivity radius for points in the convex cores of hyperbolic
3-manifolds without cusps where the manifolds are homotopy equivalent to a
closed surface. Our main result includes the possibility of cusps:

%
%

\begin{theorem}
\label{main result}
Let $M$ be an $I$-bundle over a closed hyperbolic surface $S$.  Then there
exists a constant $K_S$ such that if $N$ is a hyperbolic 3-manifold homeomorphic
to the interior of $M$ and $x \in C(N)$, then $inj_N(x) \leq K_S$.
\end{theorem} 

The main theorem partially answers a conjecture of McMullen
\cite{bielefeld} which hypothesized the existence of an upper bound on the
radius of balls that can be embedded in the convex core of hyperbolic
3-manifolds with finitely generated fundamental group, where the bound
depends on the number of generators of the fundamental group. 

In a later paper \cite{fan1}, we will use the main result of this paper to
show that the injectivity radius in the convex core is uniformly bounded
for  hyperbolic 3-manifolds homeomorphic to the interior of a book of
$I$-bundles or an acylindrical, hyperbolizable 3-manifold.  The result will
prove McMullen's original conjecture for the case of a book of $I$-bundles.   It
should be noted that for McMullen, the compressible boundary case was the
motivating case, and it still remains unsolved.  The bound on injectivity
radius, together with a theorem of McMullen's, will show that the limit set
varies continuously over the space of hyperbolic 3-manifolds homeomorphic to the
interior of a fixed hyperbolizable 3-manifold
$M$ under the geometric topology. \cite{fan2}

The first two sections of this paper will introduce some  background
material in hyperbolic geometry and review simplicial hyperbolic surfaces
which were first introduced by Bonahon.  In the third section we will review
the theory of continuous families of simplicial hyperbolic surfaces developed
by Canary.  In the last section, we will prove the main theorem.

\section{Hyperbolic Geometry and Kleinian Groups}

In this section, we introduce some of the necessary background material in
hyperbolic geometry that will be used throughout this paper.  More details
about hyperbolic geometry and Kleinian groups can be found in Beardon
\cite{beardon}, Benedetti and Petronio \cite{benedetti}, Canary-Epstein-Green
\cite{can/ep/green}, and Maskit
\cite{maskit}.  We will also introduce some results about injectivity radius
and coverings of manifolds which will be referred to throughout the paper.


We say that $\Gamma$ is a Kleinian group if it is a discrete subgroup of ${\rm
Isom}^+(\HHH)$.  If $\Gamma$ is torsion-free, then 
$N =
\HHH /
\Gamma$ is a complete, orientable hyperbolic 3-manifold. The {\it
limit set\/} of $\Gamma$ is $\Lambda_\Gamma$, the set of accumulation points of the
orbit of any
$x \in \HHH$.  The {\it domain of discontinuity\/} $\Omega_\Gamma$ is the
complement of the limit set in $\partial \HHH$, the sphere at infinity.  The {\it
convex core\/} of a hyperbolic 3-manifold
$N =
\HHH /
\Gamma$ is the smallest convex submanifold such that the inclusion map is a homotopy
equivalence.  For a nonelementary Kleinian group, it is equivalently defined as
$CH(\Lambda_\Gamma)/\Gamma$, where $CH(\Lambda_\Gamma)$ denotes the convex hull of
the limit set $\Lambda_\Gamma$, the set of accumulation points of the orbit of any
$x \in \HHH$. (Lem 6.3, Morgan
\cite{morgan}) We denote the convex core by
$C(N)$.

A hyperbolic 3-manifold $N$ is {\it geometrically finite\/} if $C_1(N)$ has
finite volume, where $C_\delta(N) = CH_\delta(\Lambda_\Gamma) / \Gamma$ and
$CH_\delta(\Lambda_\Gamma)$ is the open $\delta$-neighborhood of the convex hull of
the limit set.
Equivalently,
$N$ is geometrically finite if
$C(N)$ has finite volume and $\pi_1(N)$ is finitely generated.  In the absence of
parabolic elements, a manifold is geometrically finite if and only if its convex
core is compact. 

Let $\partial C_1(N)$ be the boundary of $C_1(N)$.  Then $\partial C_1(N)$ is
strictly convex (see Prop 2.2.1, Epstein-Marden \cite{epstein/marden}), so
that for $x \in N - C(N)$, there exists a unique geodesic ray,
$g_x$, such that $g_x$ is perpendicular
to $\partial C_1(N)$ and passes through
$x$. Then there exists a {\it nearest point retraction map\/} $r: N
\rightarrow C(N)$ where $r(x) = x$ if $x \in C(N)$, and
$r(x) = g_x \cap \partial C(N)$ otherwise.  This retraction map
gives $N - C(N)$ a product structure $\partial C_1(N) \times [0,\infty)$ where
the first coordinate is $g_x \cap \partial C_1(N)$, and the real coordinate
measures the distance from $\partial C(N)$.  Let us call these coordinates the {\it
nearest point coordinates\/} on $N - C(N)$.  The metric on
$N - C(N)$ is quasi-isometric to the metric
$cosh^2t \ dx^2_{\partial C_1(N)} + dt^2$ where $dx^2_{\partial
C_1(N)}$ is a metric on $\partial C_1(N)$.  (see Perry \cite{perry} or Sec 2.2,
 Epstein-Marden
\cite{epstein/marden})
For $\epsilon > 0$,
a similar map
$r_\epsilon: N \rightarrow C_\epsilon(N)$ retracts $N$ onto the
$\epsilon$-neighborhood of $C(N)$ and induces a homeomorphism between $N$ and
$C_\epsilon(N)$. 
 
Consider another type of nearest point retraction.  Let $N^* = (\HHH \cup
\Omega_\Gamma) / \Gamma$ be the {\it  conformal extension\/} of $N$.  For $\delta >
0$ and $x
\in C_\delta(N)$, we define $\hat{r}(x) = x$.  For $x \in N - C_\delta(N)$, the
point $\hat{r}(x)$ is the unique point of contact between $C_\delta(N)$ and the
largest hyperbolic ball centered at $x$ whose interior is disjoint from
$C_\delta(N)$.  For
$x \in
\Omega_\Gamma / \Gamma$, the point $\hat{r}(x)$ is the unique point of contact
between
$C_\delta(N)$ and the largest horoball centered at $x$ whose interior is disjoint
from $C_\delta(N)$. This map is the nearest point retraction of
$N^*$ onto
$C_\delta(N)$.

The {\it injectivity radius at a point $x$ in $N$}, $inj_N(x)$, is
defined to be half the length of the shortest homotopically non-trivial loop
passing through $x$.  It is equivalently defined as the radius of the largest
open hyperbolic ball based at $x$ that can be embedded in $N$.  The following lemma
shows that  as
$x$ travels radially out the product structure of a geometrically finite end, $inj_N(x)$ strictly
increases.

\begin{lemma}
\label{expinj} Let $N = \HHH / \Gamma$ be a hyperbolic 3-manifold. Let $U$ be a
component of $N - C(N)$, and let
$S$ be the component of $\partial C_1(N)$ associated to $U$. Then
$U$ has a product structure $S \times [0,\infty)$ with nearest point
coordinates.  For $(x,t) \in S \times [0,\infty)$, the function
$inj_N(x,t)$ is strictly increasing in $t$. 
\end{lemma}

\begin{proof} For $(x,t_0) \in S \times [0,\infty)$,
there exists a unique geodesic ray $g_x$ such that $g_x$ is perpendicular
to $S$ and passes through $(x,t_0)$.  Parametrize $g_x$ by $t$, where $t$
measures the distance from $\partial C(N)$.  Recall that for the nearest point
retraction map
$r_{t_0}: N \rightarrow C_{t_0}(N)$ which retracts $(x,t_1) \in S
\times [0,\infty)$ along $g_x$ for $t_1 > t_0$, we have $r_{t_0}(x,t_1) =
(x,t_0)$.  Because $r_{t_0}$ is strictly distance decreasing on $N - C(N)$ (see
Lem 1.4.6 in Epstein-Marden \cite{epstein/marden}), given a homotopically
non-trivial loop
$\alpha$ based at $(x,t_1)$, $r_{t_0}(\alpha)$ will be a homotopically non-trivial
loop based at $(x,t_0)$ of strictly shorter length.  So $inj_N(x,t)$ is strictly
increasing in
$t$.
\end{proof}


We will be interested in finding an upper bound on injectivity radius. Given a
bound on injectivity radius in a covering space, one can also deduce a bound on
injectivity radius in the base manifold. We will often use this fact
to find an upper bound on injectivity radius in the base manifold.

\begin{lemma} 
\label{coverinj} Let $M$ be a cover of $N$ with covering map $p$.  Then for $x
\in M$, $inj_N(p(x)) \leq inj_M(x)$.
\end{lemma}

\begin{proof} For $x \in M$, there exists a
homotopically non-trivial loop $\gamma$ based at $x$ of length
$2 \ inj_M(x)$.  Because $p$ is a local isometry, the projection $p(\gamma)$ of
the loop to $N$ has length $2 \ inj_M(x)$.  If $p(\gamma)$ is
homotopically trivial in $N$, then $p(\gamma)$ bounds a disk $D$.  By lifting
this disk to $M$, we see that $\gamma$ bounds a disk $\tilde{D}$ in $M$, which
contradicts the assumption that $\gamma$ was homotopically non-trivial in $M$. 
Thus, the loop $p(\gamma)$ must be homotopically non-trivial in $N$ so that $inj_N(p(x))
\leq inj_M(x)$.  
\end{proof}

The {\it thick-thin decomposition\/} divides a hyperbolic manifold into two
submanifolds based upon injectivity radius.  For a manifold $N$, we define
$N_{thin(\epsilon)}$ to be the set of all points in $N$ with injectivity radius
$\leq
\epsilon$.  Similarly, we define
$N_{thick(\epsilon)}$ to be the set of all points in $N$ with injectivity radius $\geq
\epsilon$.  By the Margulis Lemma (see Chap D of \cite{benedetti}), there exists
a uniform constant
$\epsilon_3$, called the {\it Margulis
constant\/}, such that for
$\epsilon <
\epsilon_3$, every component of
$N_{thin(\epsilon)}$ is one of the following:

\begin{enumerate}
\item a torus cusp $\cong T^2 \times [0,\infty)$, i.e., a horoball in $\HHH$
modulo a ${\mathbb Z} \oplus {\mathbb Z}$-parabolic action,

\item a rank one cusp $\cong {\mathbb R} \times S^1 \times [0,\infty)$, i.e., a
horoball in $\HHH$ modulo a ${\mathbb Z}$-parabolic action, or 

\item an open solid torus neighborhood of a geodesic $\cong \bar{D}^2 \times S^1$.
\end{enumerate}   Let $\NOE$ denote $N -$\{Type 1 and Type 2
components of $N_{thin(\epsilon)}$\}.  Similarly, we let $C^\circ_\epsilon(N)$
denote the $C(N) -$\{Type 1 and Type 2 components of $N_{thin(\epsilon)}$\}.

The {\it parabolic locus $P$ associated to $\NOE$\/} is the collection of
infinite annuli and tori which are the components of $\partial \NOE$.  There is a
one-to-one correspondence between the components of $P$ and the conjugacy classes
of maximal parabolic subgroups in $\pi_1(N)$.  In particular, each
representative of a parabolic element of $\pi_1(N)$ is homotopic into a component of
$P$.

An {\it end\/} $E$ of a topological space $X$ is a maximal family of sets
\{$U_i$\} called {\it neighborhoods\/} of $E$ such that: 
\begin{enumerate}
\item each $U_i$ has compact boundary but non-compact closure, and
\item for each $i$ and $j$, there exists $k$ such that $U_k \subset U_i \cap
U_j$.
\end{enumerate} 

Let $E$ be an end of $\NOE$, and let $U$ be a neighborhood of $E$.  We will give
the definition and parametrization of the {\it parabolic extension\/} $U_P$ of
$U$ given in Sec 3, Canary \cite{canary1}.   
Let $P_i$ be a rank one cusp component of
$N_{thin(\epsilon)}$.  We can normalize $P_i$ so that the fixed point is at infinity
in the upper half space model of $\HHH$ with $(z,t)$ coordinates, and $P_i = \{(z,t):t
\geq c\} / (z \mapsto z + 2\pi)$, where the constant $c$ depends only on
$\epsilon$. 
Then we can parametrize
$P_i$ by $S^1 \times (-\infty, \infty) \times [0, \infty)$ using the
coordinates $(\theta_i, r_i, s_i)$ where $\theta_i = Re \ z \ (mod \ 2\pi)$, $r_i
= H(Im \ z)$, and $s_i=log(t)+J$ for some constants $H$ and $J$ depending only on
$\epsilon$.  Then the Riemannian metric on $P_i$ can be
expressed as $e^{-2s_i}(G^2d\theta_i^2+dr_i^2)+ds_i^2$ where $G$ is a constant
depending only on $\epsilon$.   

Express $C_i = U \cap P_i$ as a subset of $S^1 \times (-\infty,\infty) \times
\{0\}$.
Let $B(C_i) =\{(\theta_i,r_i,s_i):(\theta_i,r_i,0)\in C_i \ and \ s_i \geq 0\}$
denote the portion of $P_i$ lying ``above'' $C_i$. Then the parabolic extension
$U_P$ is $U \cup [\bigcup_{P_i} B(U \cap \partial P_i)]$ where the union is taken
over all rank one cusps adjacent to $U$.  If
$A$ is a subset of
$U$, then
$A_P = A \cup [\bigcup_{P_i} B(A \cap \partial P_i)]$.

In the case that $U$ is homeomorphic to $\bar{S} \times [0,\infty)$,
there are two possibilities.  If $\bar{S}$ is closed, then the parabolic
extension $U_P$ of $U$ is simply $U$.  If $\bar{S}$ is not closed, then
$U_P = U \cup [\bigcup_{P_i} B(\partial U \cap P_i)]$, where the union is taken
over all rank one cusps adjacent to $U$.  Here $U_P$ is homeomorphic to $S
\times [0,\infty)$ where $S$ is homeomorphic to the interior of $\bar{S}$.

We say an end $E$ of $\NOE$ is {\it geometrically finite\/} if
there exists a neighborhood $U$ of $E$ such that $U \cap C(N) = \emptyset$. 
Otherwise we say that the end $E$ is {\it geometrically infinite\/}.

A {\it pleated surface\/} is a map $f: S \rightarrow N$ from a surface
$S$ with hyperbolic metric $\sigma$ to a hyperbolic 3-manifold $N$ such that: 
\begin{enumerate}
\item $f$ is a pathwise isometry with respect to the metric $\sigma$ on $S$, and
\item for each point $x \in S$, there exists a geodesic arc $\gamma$ such that $x
\in int \ \gamma$ and $f(\gamma)$ is a geodesic arc.
\end{enumerate}  
A {\it geodesic lamination\/} on a surface $S$ is a closed subset of $S$ which is a
disjoint union of simple geodesics which are called {\it leaves\/}.  For a pleated
surface
$f:S
\rightarrow N$, its {\it pleating locus\/} is the smallest geodesic lamination
$\lambda$ such that $f$ is totally geodesic in the complement of $\lambda$, and $f$
maps the leaves of
$\lambda$ geodesically.   Epstein-Marden (Chap 1, \cite{epstein/marden}) have
shown that the boundary components of convex cores are images of pleated surfaces. 

A useful tool will be the following lemma which states that for a cover 
associated to a boundary component of the convex core of a manifold, a component
of the complement of the convex core of the base space lifts homeomorphically to a
component of the complement of the convex core of the cover.

\begin{lemma}
\label{GFcover}  Let $\delta > 0$. Let $N$ be a hyperbolic 3-manifold.   Let $U$
be a component of
$N - int \ C(N)$, and let $S$ be
the component of
$\partial C(N)$ associated to $U$. Let $M = \HHH / \pi_1(S)$ be a cover of $N$
with projection map $p$.  Then there exists a lift $\tilde{U}$ of $U$ such
that $p|_{\tilde{U}}: \tilde{U} \rightarrow U$ is a homeomorphism, 
$(p|_{{\tilde{U}}})^{-1}(S)$ is a component of $\partial C(M)$, and
$\tilde{U}$ is a component of $M - int \ C(M)$. 
Furthermore, if $T \subset U$ is a component
of $\partial C_\delta(N)$, then $(p|_{{\tilde{U}}})^{-1}(T)$ is a
component of $\partial C_\delta(M)$.
\end{lemma}

\begin{proof}  Because $U$ is the component of $N -
int \ C(N)$ associated to $S$, it has a product structure $S \times
[0,\infty)$ with nearest point coordinates. By the Lifting Theorem, the inclusion
map $i: S \times [0,\infty) \rightarrow U$ lifts to a map $\tilde{i}: S
\times [0,\infty) \rightarrow M$. Let $\tilde{U} = \tilde{i}(S \times [0,\infty))$. 
Then  $p|_{\tilde{U}}: \tilde{U} \rightarrow U$
is a homeomorphism.

Recall that $S$ is the image of some pleated surface $f: S^\prime \rightarrow
N$. Then there exists a lift $\tilde{f}: S^\prime
\rightarrow M$ of $f$ such that $\tilde{f}(S^\prime) =
(p|_{{\tilde{U}}})^{-1}(S)$.  Let $\lambda$ be the pleating locus of
$\tilde{f}$.  Recall that a pleated surface $\tilde{f}: S^\prime \rightarrow M$ is
totally geodesic in the complement of $\lambda$, and $\tilde{f}$ maps the leaves of
$\lambda$ geodesically.  Epstein-Marden (Cor 1.6.3 \cite{epstein/marden})
have shown that the lifts to the universal cover of the
leaves of $\lambda$ have endpoints that lie in
$\Lambda_{\pi_1(S)}$. Therefore,
$(p|_{{\tilde{U}}})^{-1}(S) \subset C(M)$. 

If $(p|_{{\tilde{U}}})^{-1}(S)$ is not a boundary component of $C(M)$,
then there exists a geodesic arc $\gamma$ in $int \ \tilde{U}$ with endpoints in
$(p|_{{\tilde{U}}})^{-1}(S)$, i.e.,
$p(\gamma)$ is a geodesic arc with endpoints in $S$ whose interior lies
outside of
$C(N)$.  This contradicts the fact that
$S$ is a boundary component of a convex set $C(N)$.  Therefore,
$(p|_{{\tilde{U}}})^{-1}(S)$ is a boundary component of $C(M)$.  
So $\tilde{U}$ is a component of
$M - int \ C(M)$.  

Let $T \subset U$
be a component of $\partial C_\delta(N)$.  Then for each point $x$ on $T$, there
exists a geodesic arc $g_x$ perpendicular to $T$ and passing through $x$.
Lift $x$ to a point $\tilde{x} \in \tilde{U}$.  Because
$p$ is an isometry on $U$, this geodesic arc lifts to a geodesic arc
$g_{\tilde{x}}$ that is perpendicular to $(p|_{{\tilde{U}}})^{-1}(S)$ and passes
through $\tilde{x}$.  Then the set of points in $U$ a distance $\delta$ from
$S$ lifts precisely to the set of points in $\tilde{U}$ a distance $\delta$ from
$(p|_{{\tilde{U}}})^{-1}(S)$.  Hence $(p|_{{\tilde{U}}})^{-1}(T)$ is a
component of $\partial C_\delta(M)$.
\end{proof}

We can also show that if a neighborhood of a geometrically finite end of a covering
space $M$ embeds as a neighborhood of a geometrically finite end in the base
manifold, then the neighborhood associated to $M - int \ C(M)$ embeds in the
base manifold as well.

\begin{lemma}
\label{GFproj}  Let $0 < \epsilon < \epsilon_3$, and let $\delta > 0$.  Let $M
= \HHH / \hat{\Gamma}$ be a hyperbolic 3-manifold which covers another hyperbolic
3-manifold $N = \HHH / \Gamma$ with projection map $p$.  Let
$E$ be a geometrically finite end of $M^\circ_\epsilon$, and let
$\tilde{V}$ be a component of $M - int \ C_\delta(M)$ such
that
$\tilde{V} \cap M^\circ_\epsilon$ is a neighborhood of $E$.  Let
$\tilde{S} = \partial C_\delta(M) \cap \tilde{V}$.   Suppose
there exists a neighborhood $\tilde{U}$ of $E$ such that $\tilde{U}
\subset \tilde{V}$ and $\tilde{U}$ embeds in $N$.  Then
$p(\tilde{V})$ is a component of $\NOE - C_\delta(N)$, and 
$S = p(\tilde{S})$ is a boundary component of $\partial
C_\delta(N)$. 
\end{lemma}

\begin{proof} Recall the deformation retract $\hat{r}:
N^* \rightarrow C_\delta(N)$.  This map lifts to the deformation retract
$\tilde{r}: \HHH \cup \Omega_\Gamma \rightarrow CH_\delta(\Lambda_\Gamma)$.  Then
given a component $\tilde{S} = \partial C_\delta(M) \cap \tilde{V}$ of $\partial
C_\delta(M)$, we can consider a component $\bar{S}$ of the preimage of $\tilde{S}$
in $\HHH$. Here,
$\bar{S}$ will be a boundary component of $CH_\delta(\Lambda_{\hat{\Gamma}})$. 
This component $\bar{S}$ corresponds to a component $\Omega_{S}$ of
$\Omega_{\hat{\Gamma}}$, where
$\Omega_{S} = \tilde{r}^{-1}(\bar{S}) \cap \Omega_\Gamma$.  Let
$\bar{V}$ be the portion of $\HHH$ that is bounded by $\bar{S}$
and $\Omega_{S}$.  Let $\bar{U}$ be the portion of the preimage of $\tilde{U}$
in
$\HHH$ that is contained in $\bar{V}$. Then
$\bar{V}$ covers $\tilde{V}$, and $\bar{U}$ covers
$\tilde{U}$. Furthermore because $\tilde{U}$ is a neighborhood of the geometrically
finite end $E$, $\bar{U}$ must contain a neighborhood of $\Omega_{S}$ in
$\bar{\mathbb H}^3$. 

Because $\Lambda_{\hat{\Gamma}} \subset \Lambda_\Gamma$, each component of
$\Omega_\Gamma$ is contained in a component of $\Omega_{\hat{\Gamma}}$.  Suppose
$\Omega_{S}$ is not a component of $\Omega_\Gamma$.  Then there exists a point $z$
in $\Lambda_{\Gamma}$ such that $z$ is also contained in $\Omega_{S}$.  We know
that pairs of hyperbolic fixed points are dense in $\Lambda_{\Gamma} \times
\Lambda_{\Gamma}$ (Prop V.E.5, Maskit \cite{maskit}), i.e., for a pair
of points $(a,b) \in \Lambda_{\Gamma} \times
\Lambda_{\Gamma}$, there exists a sequence of hyperbolic isometries whose
fixed points approach $(a,b)$.  So there exists a hyperbolic isometry
$\gamma \in \Gamma - \hat{\Gamma}$
whose fixed points lie in a neighborhood $W$ of $z$ such that $W \subset
\Omega_{S}$. Furthermore $\bar{U}$ must contain a neighborhood of $W$. 
Thus, 
$\gamma(\bar{U}) \cap \bar{U} \neq \emptyset$.  Then because $\bar{U}$ covers
$\tilde{U}$, we can conclude that
$\tilde{U}$ does not embed in
$N$, but this is a contradiction.  So $\Omega_{S}$ is a component of
$\Omega_\Gamma$.  

Then $\bar{S}$ is also a boundary component of $CH_\delta(\Lambda_\Gamma)$, and
for $\gamma \in \Gamma - \hat{\Gamma}$, $\gamma(\bar{V}) \cap \bar{V} = \emptyset$.
So we can conclude that $\tilde{V}$ embeds in $N$, $p(\tilde{V})$ is a component of
$N - C_\delta(N)$, and
$S = \bar{S}/ \Gamma = p(\tilde{S})$ is a boundary component of
$\partial
C_\delta(N)$. \end{proof}

It is conjectured that all geometrically infinite ends of a hyperbolic 3-manifold
with finitely generated fundamental group are simply degenerate.  To
explicitly define a simply degenerate end, we need to define the term ``simplicial
hyperbolic surface,'' but we will delay the explicit presentation this definition
until the next section where we will explore these objects more thoroughly.
Roughly, a simplicial hyperbolic surface is a map of a surface whose image has
intrinsic curvature $\leq -1$.   An end of
$\NOE$ is {\it simply degenerate\/} if:
\begin{enumerate}
\item there exists a neighborhood $U$ of the end that is homeomorphic to $\bar{S}
\times [0, \infty)$, where $\bar{S}$ is a compact surface with or without boundary,

\item there exists a sequence of simplicial hyperbolic surfaces $\{f_n:
S \rightarrow U_P\}$, where $S$ is the interior of $\bar{S}$, such
that
$f_n(S)$ is properly homotopic within $U_P$ to $S \times \{0\}$,
and 
\item the sequence \{$f_n(S)$\} exits the end of $N$.
\end{enumerate}

For a hyperbolic 3-manifold $N$, a {\it compact core\/} $R$ of $N$ is a
compact submanifold of $N$ such that the inclusion of $R$ into
$N$ is a homotopy equivalence.  A theorem of Scott \cite{scott} guarantees that
every hyperbolic 3-manifold with finitely generated fundamental group has a
compact core. 

For $0 < \epsilon < \epsilon_3$, a {\it relative compact core of $\NOE$\/} is a
compact core $R$ of $\NOE$ that intersects every rank one component of $\partial
\NOE$ in a $\pi_1$-injective annulus and intersects every rank two component of
$\partial
\NOE$ in a torus.  The existence of relative compact cores is
guaranteed by McCullough \cite{mccullough} and Kulkarni-Shalen
\cite{kul/shalen}.   The collection of annuli and tori in $\partial R \cap \partial
\NOE$ is denoted $P$ and is called the {\it parabolic locus associated to $R$\/}.
There is a one-to-one correspondence between the components of $P$, the parabolic
locus associated to $\NOE$, and the components of $P$, the parabolic locus
associated to $R$.    It will be clear from context to which object a given
parabolic locus is associated.  

We can also define a more general version of a relative compact core. For $0 <
\epsilon <
\epsilon_3$, let
${\mathcal P}$ be a collection of Type I and Type II components of
$N_{thin(\epsilon)}$.  A {\it relative compact core of $N - {\mathcal P}$\/}
is a compact core $R$ of $N - {\mathcal P}$ that intersects every rank one component of
$\partial {\mathcal P}$ in a $\pi_1$-injective annulus and intersects every rank two
component of
$\partial {\mathcal P}$ in a torus. In this case, the collection of annuli and tori in
$\partial R \cap \partial {\mathcal P}$ is also denoted $P$ and is also called the {\it
parabolic locus associated to
$R$\/}.  Again, it will be clear from the context to which compact core a given
parabolic locus is associated.

A hyperbolic 3-manifold $N$ with finitely generated fundamental group is said
to be {\it geometrically tame\/} if each end of $\NOE$ is either
geometrically finite or simply degenerate. A sufficient condition for a
manifold to be geometrically tame is Bonahon's Condition (B):

\begin{theorem}
\label{geotame} (Thm A, Bonahon \cite{bonahon})  Let $N$ be a complete
hyperbolic 3-manifold such that $\Gamma = \pi_1(N)$ is of finite type and
satisfies Condition (B), i.e., for every non-trivial free decomposition
$\Gamma = A \ast B$, there exists a parabolic element in $\Gamma$ which is
conjugate into neither $A$ nor $B$.  Then $N$ is geometrically tame.
\end{theorem}

Condition (B) is equivalent to requiring that the fundamental group of every
component of $\partial R - \partial \NOE$ inject into $\pi_1(N)$, where $R$ is
a relative compact core of $\NOE$.  

By contrast, a hyperbolic 3-manifold $N$ is {\it topologically tame\/} if it
is homeomorphic to the interior of a compact manifold.  By definition, a
manifold that is geometrically tame is topologically tame.   Canary
\cite{canary2} has proven the converse as well.  In this paper, all of our
manifolds will be homeomorphic to the interior of a compact manifold and hence
will be geometrically and topologically tame.

Given a collection ${\mathcal P}$ of Type I and Type II components of
$N_{thin(\epsilon)}$, and a relative compact core $R$ of $N - {\mathcal P}$ with
associated parabolic locus $P$,  Bonahon (Prop 1.6, \cite{bonahon}) has shown
that if $N$ satisfies Condition (B), then the ends of $N - {\mathcal P}$ are in
one-to-one correspondence with the components of $(N - {\mathcal P}) - R$, or
equivalently, with the components of $\partial R - P$.

A nice property of topologically tame manifolds is that given a compact
set in a topologically tame manifold, it is possible to find a compact core that
contains the compact set.

\begin{lemma}
\label{compactincore}  Let $N$ be a topologically tame hyperbolic 3-manifold.  Let
${\mathcal P}$ be a collection of Type I and Type II components of
$N_{thin(\epsilon)}$.  Let $K$ be a compact set in $N - {\mathcal P}$.  Then there
exists a relative compact core $R$ of
$N - {\mathcal P}$ such that $K \subset R$ and the components of $(N - \mathcal{P}) - R$ are
topologically a product.
\end{lemma}

\begin{proof} Because $N$ is topologically
tame, there exists a relative compact core
$R^*$ of $N - {\mathcal P}$ with associated parabolic locus $P$ such that each component
of $\partial R^* - P$ is incompressible and each component $U_j$ of $(N - {\mathcal
P}) - int \ R^*$ is homeomorphic to a component of $(\partial R^* - P) \times
[0,\infty)$. (see Chap 1, Bonahon \cite{bonahon})  Choose a level surface
$T_j$ in each $U_j$ so that $K$ lies to one side. 
Then each $T_j$ is incompressible and properly isotopic to a component of
$\partial R^* - P$.  Thus, the $T_j$ bound a
relative compact core $R$ of $N - {\mathcal P}$ such that $K \subset R$. 
Furthermore,
the components of $(N - \mathcal{P}) - R$ are
topologically a product.  \end{proof}

Another property of topologically tame manifolds is that for any relative compact
core $R$ of a topologically tame manifold $N$ where the components of $\partial R -
P$ are incompressible, the components of
$(N - {\mathcal P}) - int \ R$ possess a product structure.

\begin{lemma}
\label{prodstructure}  Let $N$ be a topologically tame hyperbolic 3-manifold.  Let
${\mathcal P}$ be a collection of Type I and Type II components of
$N_{thin(\epsilon)}$.  Let $R$ be a relative compact core of $N - {\mathcal P}$ with
associated parabolic locus $P$.  Let
$\{S_j\}$ be the components of
$\partial R - P$, and let $U_j$ be the component of $(N - {\mathcal P}) - int \ R$
associated to
$S_j$.  If the
\{$S_j$\} are incompressible, then the \{$U_j$\} are homeomorphic to $S_j \times [0,
\infty)$. 
\end{lemma}

\begin{proof}  By the argument given in Lemma
\ref{compactincore}, there exists a relative compact core
$R^*$ such that each component of $(N - {\mathcal P}) - int \ R^*$ is
 topologically a
product and $R \subset R^*$. Select a surface  $T_j$ in a component of
$(N - {\mathcal P}) - int \ R^*$.  Let
$A_j$ be the subset of $U_j$ bounded by $T_j$ and $S_j$. If $A_j$ is not
homeomorphic to $S_j \times [0,1]$, then 
$\pi_1(S_j)$ does not surject onto
$\pi_1(A_j)$. (Thm 10.5, Hempel \cite{hempel}) Then there exists a
homotopically non-trivial loop $\gamma$ in $A_j$ that is not homotopic in $A_j$ to
a curve in $S_j$.  Because $S_j \subset \partial R - P$ is incompressible and
separating in
$N -
\mathcal{P}$,
$\gamma$ is not homotopic in
$N - {\mathcal P}$ to a curve in $R$. But $R$ is a compact core of $N$, so this is a
contradiction. Then for each $j$, $U_j$ is homeomorphic to $S_j \times [0,\infty)$. 
\end{proof}

Now we present a general fact about ends of topologically tame manifolds and their covers
which will be used throughout the paper.

\begin{lemma}
\label{not in C(N)} Let $N$ be a topologically tame hyperbolic 3-manifold.  Let $0 <
\epsilon < \epsilon_3$, and let $\delta \geq 0$.  Let ${\mathcal P}$ be a collection of
Type I components of $N_{thin(\epsilon)}$.  Let $S$ be an
incompressible separating surface of $N - {\mathcal P}$ such that $\partial S \subset
\partial {\mathcal P}$.  Let $M = \HHH / \pi_1(S)$ be a cover of $N$ with
covering map
$p$.  Let $U$ be the closure of a component of $(N - {\mathcal P}) - S$ such that
there exists
$\tilde{U} \subset M$ such that
$p|_{\tilde{U}}:\tilde{U} \rightarrow U$ is a homeomorphism and
$(p|_{{\tilde{U}}})^{-1}(S) \subset \overline{C_\delta(M)}$.  Then
$$(p|_{{\tilde{U}}})^{-1}[(N
- C_\delta(N)) \cap U] = (M - C_\delta(M)) \cap
\tilde{U}$$ 
\end{lemma}

\begin{proof}  $\subseteq$:  We know that
$\pi_1(S)$ is a subgroup of $\pi_1(N)$ so that $\Lambda_{\pi_1(S)} \subset
\Lambda_{\pi_1(N)}$.  Thus,
$CH(\Lambda_{\pi_1(S)}) \subset CH(\Lambda_{\pi_1(N)})$, and hence $p(C(M))
\subset C(N)$.  Therefore,
$p^{-1}(N - C_\delta(N)) \subset M - C_\delta(M)$.

$\supseteq$:  We will show that for a point $x \in [M - C_\delta(M)]
\cap \tilde{U}$, we have $p(x) \in [N - C_\delta(N)] \cap U$. 
If $\tilde{U} \subset C_\delta(M)$, then this is automatically true.  

If $\tilde{U} \not \subset C_\delta(M)$, then let $B$ be a component of $\partial
C_\delta(M)$ that separates $x$ from $(p|_{{\tilde{U}}})^{-1}(S)$.  Let $A$ be the
component of $M - B$ that does not contain
$C_\delta(M)$.  So in particular, $A$
contains $x$.  Then $M - A$ is a locally convex submanifold of $M$ with
boundary $B$, and hence $M - A$ is
convex. (Cor 1.3.7, Canary-Epstein-Green \cite{can/ep/green}) 

Using the $(\theta_i,r_i,s_i)$ coordinates on each component of ${\mathcal
P}$ adjacent to $\partial S$, we can extend the homeomorphism
$p|_{\tilde{U}}$ to a homeomorphism
$\bar{p}:\tilde{U}_P \rightarrow U_P$ where $U_P$ is the parabolic
extension of $U$.  Without loss of generality, we can assume that $\epsilon$ is
small enough so that $B$ intersects each component of $U_P - U$ totally
geodesically. Thus, we know that
$B
\subset
\tilde{U}_P$.  Let
$X_\nu$ be a component of
$M - int \ C_\nu(M)$ which intersects
$\tilde{U}$.  Then because $\bar{p}$ is a homeomorphism, for large enough $\nu$,
$X_\nu$ embeds in $N$.  Then by Lemma \ref{GFproj}, $p(\bar{A})$ is a component of
$N - int \ C_\delta(N)$ with boundary $p(B)$.  Therefore, $N -
p(A)$ is a convex submanifold of $N$.   We can
construct a deformation retract
$R: N \rightarrow N - p({A})$ by letting
$R|_{N - p({A})}= id|_{N - p({A})}$ and $R|_{p(\bar{A})} =
r_\delta|_{p(\bar{A})}:p(\bar{A}) \rightarrow p(B)$ where $r_\delta: N \rightarrow
C_\delta(N)$ is the nearest point retraction map.  Then
$N - p({A})$ is a convex submanifold of $N$ which carries the homotopy of
$N$.  Therefore, $C(N) \subset N - p({A})$. 

By Lemma \ref{GFproj}, we know that for the component $T$ of
$\partial C(M)$ corresponding to $B$, $p(T)$ is a component of $\partial C(N)$. Then
because
$p|_{U}$ is an isometry,
$[C_\delta(M) - C(M)] \cap \tilde{U}$ projects homeomorphically 
to $[C_\delta(N) -
C(N)] \cap U$. 
Thus we can conclude that $C_\delta(N) \subset N - p({A})$.

So if $x \in (M - C_\delta(M)) \cap \tilde{U}$, then $x \in A$.  Then $p(x) \in p(A)
\subset U$ and hence $p(x) \not \in C_\delta(N)$.  
\end{proof}

\section{Simplicial Hyperbolic Surfaces}

\subsection{Triangulations} 

In order to define a simplicial hyperbolic surface, we must first define a
triangulation on a punctured surface.  Let $\hat{S}$ be a (not necessarily 
orientable) surface with a finite number of distinguished points $\{p_1, \ldots,
p_n\}$.  Let
$S$ denote
$\hat{S} - \{p_1, \ldots, p_n\}$.  Then $S$ is a {\it punctured surface\/}, and the
$p_i$ are {\it punctures\/}.  Let $V$ denote the punctures
$\{p_1,\ldots, p_n\}$ and a finite collection of vertices on $S$.  The punctures are
{\it ideal vertices\/}, and the vertices that are not punctures are {\it internal
vertices\/}.  We require that $V$ include at least one internal vertex. 

A {\it curve system\/} $\{\alpha_i, \ldots, \alpha_m\}$ is a collection of arcs in
$\hat{S}$ which have endpoints in $V$, which are disjoint except at their endpoints,
no two of which are ambient isotopic relative to $V$, and none of which is homotopic
to a point relative to
$V$. Then a {\it triangulation\/} of a surface $S$ with associated vertices $V$ is
the restriction of a maximal curve system on $(\hat{S}, V)$ to $(S, V)$.  Two
triangulations are equivalent if they are ambient isotopic relative to
$V$. 

\begin{figure}[ht]
\begin{center}
\BoxedEPSF{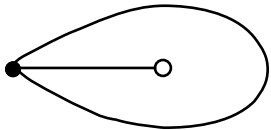}
\caption{A Punctured Monogon} \label{fig:puncmonogon}
\end{center}
\end{figure}

Note that the ``triangles'' in these triangulations may not have three distinct
edges nor three distinct vertices.  For example, one edge may count as two sides of
a triangle, or there may only be one vertex in a triangle.  In this paper, we will
often use a certain type of triangle called a punctured monogon.  A {\it punctured
monogon\/} consists of an ideal vertex, an internal vertex, an edge joining the
internal vertex and the ideal vertex, and an edge with both endpoints at the
internal vertex that encloses precisely the ideal vertex and the edge joining the
internal and ideal vertices.  See Figure ~\ref{fig:puncmonogon}. If each ideal
vertex is contained in a punctured monogon, then we say the triangulation is {\it
suitable.\/}  

\begin{figure}[ht]
\begin{center}
\BoxedEPSF{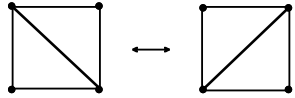 scaled 700}
\caption{A Type I Elementary Move} \label{fig:type1move}
\end{center}
\end{figure}

\begin{figure}[ht]
\begin{center}
\BoxedEPSF{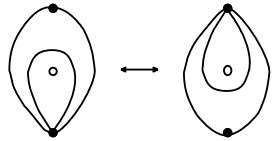 scaled 700}
\caption{A Type II Elementary Move} \label{fig:type2move}
\end{center}
\end{figure}

Harer and Hatcher \cite{harer1,harer2,hatcher} use a more general definition of a
triangulation of a surface with or without punctures: \ A {\it generalized
triangulation\/} of a surface
$S$ is a maximal curve system on $(S,V)$, where $V$ may or may not include all of the
ideal vertices. 
 They have shown that any two generalized triangulations of a surface can be related
by a finite series of {\it elementary moves\/}, and that there are exactly two types
of elementary moves, Type I and Type II.  Given a quadrilateral with a diagonal, a
Type I move replaces the given diagonal with the other diagonal.  See Figure
~\ref{fig:type1move}. Given a punctured bigon containing an edge with both
endpoints at one vertex, a Type II move replaces the given edge with the edge with
both endpoints at the other vertex.  See Figure ~\ref{fig:type2move}.

\begin{figure}[ht]
\begin{center}
\BoxedEPSF{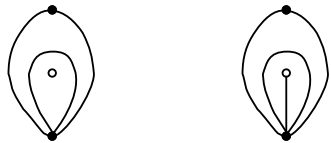 scaled 700}
\caption{A Punctured Bigon in a Generalized Triangulation (left) and in a Suitable
Triangulation (right)}
\label{fig:2puncbigons}
\end{center}
\end{figure}

In this paper, our definition of triangulation is more restrictive, because we
require that all punctures be vertices of the triangulation.  Thus, under our
definition of triangulation, a punctured bigon containing an edge with both
endpoints at one vertex must also contain an edge connecting the ideal vertex to the
internal vertex.  Hence, in our situation, there are no punctured bigons.  See Figure
~\ref{fig:2puncbigons}. In other words, Type II elementary moves are unnecessary. 
Therefore, in this paper, given a triangulation, we can obtain any other
triangulation using only Type I elementary moves.

\begin{figure}[ht]
\begin{center}
\BoxedEPSF{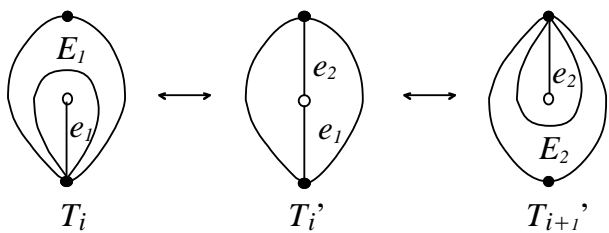 scaled 700}
\caption{A Type II Elementary Move Can Be Replaced By Two Type I Elementary Moves}
\label{fig:twotype1}
\end{center}
\end{figure}

Notice that in a suitable triangulation, there are no {\it doubly ideal edges\/},
i.e., edges with both endpoints at ideal vertices.  In fact, given any suitable
triangulation, we can obtain any other suitable triangulation using a specific
sequence of Type I elementary moves in which each triangulation in the sequence has
no doubly ideal edges.  Suppose we are given two suitable triangulations
$T_A$ and $T_B$.  Then a
specific sequence of solely Type I elementary moves altering
$T_A$ to $T_B$ can be realized by replacing each Type II elementary move by two
Type I elementary moves as follows:  
\begin{enumerate}
\item convert the suitable triangulations $T_A$ and $T_B$ into generalized
triangulations ${\mathcal T}_A$ and ${\mathcal T}_B$ by removing the set of edges
with one endpoint at an ideal vertex;

\item the result of Harer and Hatcher guarantees the existence of a sequence of
Type I and Type II elementary moves that alters the triangulation ${\mathcal T}_A$
on $(S,V^\prime)$ to the triangulation
${\mathcal T}_B$ on $(S,V^\prime)$ where $V^\prime$ contains only internal vertices; 

\item at each step in the Harer--Hatcher sequence, there is either a Type I move or
a Type II move:
			\begin{enumerate}
			\item for a Type I move, perform a Type I move on ${\mathcal
									T}_i$ to obtain ${\mathcal	T}_{i+1}$.  Then let	$T_{i+1} = {\mathcal
									T}_{i+1} \cup (T_i - {\mathcal	T}_i)$.  The	resulting	suitable triangulation
									$T_{i+1}$ contains no doubly	ideal edges; and  
			\item for a Type II move, perform a sequence of two
									Type I moves on $T_i$.   Consider the punctured bigon in ${\mathcal T}_i$ on
									which	the Type II move is to be performed.  In the corresponding
									punctured bigon in
									$T_i$, let $E_1$ be the edge with both endpoints at the internal
									vertex $v_1$ which
									surrounds the
									puncture $p$, and let	$e_1$
									be	the edge joining $p$ and $v_1$.  Let	$e_2$
									denote	an edge joining $p$ and the other internal vertex $v_2$,
									and
									let	$E_2$ denote the edge with both
									endpoints at $v_2$ which surrounds $p$. 
									Let
									$T_i^\prime$ be the triangulation $T_i - E_1 + e_2$, and let $T_{i+1}$
									be
									the triangulation $T_i^\prime - e_1 + E_2$.  Then $T_{i+1}$ is a
									suitable triangulation, and $T_i^\prime$ and
									$T_{i+1}$ contain no doubly ideal edges. See
									Figure ~\ref{fig:twotype1}.
			\end{enumerate}
\end{enumerate}

After renumbering the sequence, each of the triangulations $T_j$ in the resulting
sequence   contains no doubly ideal edges.  Let us call this sequence of
triangulations a {\it modified HH sequence\/}. In this paper, we can restrict our
attention to triangulations that contain no doubly ideal edges.

\subsection{Simplicial Hyperbolic Surfaces}

A homomorphism $\alpha: \Theta \rightarrow \Gamma$ between two Kleinian groups is
{\it type-preserving\/} if $\alpha(\gamma)$ is parabolic if and only if
$\gamma$ is parabolic.  By contrast, a homomorphism {\it weakly preserves
parabolicity\/} if $\alpha$ has the property that $\alpha(\gamma)$ is parabolic if
$\gamma$ is parabolic.  We say that $\alpha(\gamma)$ is an  {\it accidental
parabolic\/} if $\alpha(\gamma)$ is parabolic, whereas $\gamma$ is hyperbolic. 
Recall that each accidental parabolic is associated to a
closed curve in the parabolic locus
$P$ of
$\NOE$.  Sullivan \cite{sullivan} has shown that for a finitely generated Kleinian
group, there can be at most a finite number of conjugacy classes of maximal abelian
subgroups associated to accidental parabolics.

For any $\pi_1$-injective map between hyperbolic manifolds, we can consider a map
$f:S \rightarrow N$ and its induced map $f_*$ on fundamental groups.  Then
$f$ {\it weakly preserves parabolicity\/} if $f_*$ weakly preserves parabolicity. 
Similarly, a homotopically non-trivial curve $\gamma$ in $S$ is an {\it accidental
parabolic\/} if $\gamma$ is not homotopic into a cusp of $S$ but $f(\gamma)$ is
homotopic into a cusp of $N$.  This definition agrees with our previous group
theoretic definition.

Now we can present the definition of a simplicial hyperbolic surface.  Let $f: S
\rightarrow N$ be a proper map from a (not necessarily orientable) surface $S$ into a
3-manifold
$N$ where
$S$ has triangulation
$T$. We will use the notation $f: (S,T) \rightarrow N$ to denote such a map.  Suppose
$f$ has the following properties: 
\begin{enumerate}
\item $f$ weakly preserves parabolicity, 
\item $f$ maps every edge $e$ in $T$ to a geodesic arc, and
\item $f$ maps each face of $T$ to a non-degenerate, totally geodesic triangle in
$N$.
\end{enumerate}    Then $f$ is a {\it simplicial pre-hyperbolic surface\/}.  Let
$ang \ f(v)$ be the total angle about a vertex $f(v)$, that is, the sum of the
angles based at $f(v)$ in each of the geodesic triangles which share $f(v)$ as a
vertex.  If $ang \ f(v) \geq 2\pi$ for every internal vertex $v$, then
$f$ is a {\it simplicial hyperbolic surface\/}.  This angle condition guarantees
that the intrinsic geometry of $f(S)$ is like that of a surface of curvature
$\leq -1$.

The map $f$ induces a piecewise Riemannian metric on $S$ called the {\it simplicial
hyperbolic structure\/}, denoted $\tau$.  The surface
$(S,\tau)$ has curvature $\leq -1$ at every point except at the vertices which have
``concentrated'' negative curvature, because the total angle about each vertex is
$\geq 2\pi$.  By the Gauss-Bonnet Theorem for hyperbolic triangles, the area of
$(S,\tau)$ is
$$area(S,\tau) = 2\pi|\chi(S)|-\sum_{v\in T}(ang \ f(v)-2\pi)$$  (Lem 1.13, Bonahon
\cite{bonahon}) In particular,
$area(S,\tau) \leq 2\pi|\chi(S)|$.


We will only be interested in $\pi_1$-injective simplicial hyperbolic surfaces. 
Using the area bound of
$(S,\tau)$, we can deduce a bound on injectivity radius for points in the image of a
$\pi_1$-injective simplicial hyperbolic surface.

\begin{lemma} 
\label{shsinj} Let $f: S \rightarrow N$ be a $\pi_1$-injective simplicial hyperbolic
surface.  Then there exists a constant $K_S$ that depends on the Euler
characteristic of $S$ such that for $x \in f(S)$, $inj_N(x) \leq K_S$.  
\end{lemma}

\begin{proof} For each point
$x \in S$, there exists an embedded disk $D$ of radius $inj_N(x)$ in $S$.  Because
the curvature of $(S, \tau)$ is $\leq -1$, in $(S,\tau)$, $area(D) \geq 4\pi
sinh^2[\frac{1}{2}(inj_N(x))]$, the latter being the area of a disk with the same
radius in hyperbolic space.  Since 
$area \ [f(S)] \leq 2\pi|\chi(S)|$, for any $x$ in the image of a simplicial
hyperbolic surface,
$inj_N(x) \leq K_S$ where the constant
$K_S$ depends on the Euler characteristic of $S$. \end{proof}

Now let us describe an explicit construction due to Bonahon \cite{bonahon} of a
simplicial pre-hyperbolic surface which will be used later to build a continuous
family of simplicial hyperbolic surfaces.  Let $T$ be a triangulation of
$S$ with no doubly ideal edges, and let $f: S \rightarrow N$ be a proper map that
weakly preserves parabolicity, maps every edge $e$ in $T$ with both endpoints at the
same internal vertex to a homotopically non-trivial loop in $N$, and maps no two
vertices of $T$ to the same point.  First, we homotop $f$, keeping $f(V)$ fixed, to
a map $f_1: S
\rightarrow N$ such that if $e$ is a finite edge in $T$, then $f_1(e)$ is the unique
geodesic arc in its homotopy class.  Then we properly homotop $f_1$, fixing $\bigcup
f_1(e)$ over all finite edges $e$ to a map $f_2: S \rightarrow N$ such that if
$e^\prime$ is a half infinite edge, then
$f_2(e^\prime)$ is the half infinite geodesic ray which has the same endpoint and is
properly homotopic to $f(e^\prime)$.  Note that in the universal cover, each
component of the  lift of
$f_2(e^\prime)$ connects a lift of $f(v)$ and a fixed point of a parabolic element
of $\pi_1(N)$.  Recall that we are assuming that there are no doubly ideal edges, so
there are no edges in the triangulation with endpoints at two distinct ideal
vertices.  Finally, we properly homotop
$f_2$, fixing $f_2(T)$, to a map $F(f, T): S \rightarrow N$ such that each face of
$T$ is taken to the totally geodesic triangle spanned by the images of its edges. 
This homotopy results in a well-defined map $F(f,T)$ which is a simplicial
pre-hyperbolic surface.  

Canary \cite{canary1} proved that if each vertex of a simplicial pre-hyperbolic
surface $f: S \rightarrow N$ is {\it NLSC\/}, then $f: S \rightarrow
N$ is a simplicial hyperbolic surface. The criterion $NLSC$---$NLSC$ stands for Not
Locally Strictly Convex---can be described as follows: \ For each vertex
$v$ and an edge $e$ with endpoint at $v$, we associate a point in
$T^{1}_{f(v)}(\HHH)$ which is the unit tangent vector in the direction of
$f(e)$ (pointing away from $f(v)$).  Let $S_{f(v)}$ denote the set of all points in 
$T^{1}_{f(v)}(\HHH)$ associated to edges with an endpoint at $v$. Note that if an
edge has both endpoints at $v$, then we will associate two points in
$T^{1}_{f(v)}(\HHH)$ to this edge.  We say that $v$ is {\it NLSC with respect to\/}
$f$ if there is no open hemisphere in $T^{1}_{f(v)}(\HHH)$ which contains
$S_{f(v)}$.  Alternatively, a vertex $v$ is {\it NLSC\/} if in the universal cover,
there exists no open half-space $H$ containing $\tilde{f}(\tilde{v})$ such that all
the edges originating from $\tilde{f}(\tilde{v})$ lie in $H$.

A useful way to check for {\it NLSC\/} is via the following lemma:

\begin{lemma}
\label{tetrahedron} (Lem 4.3, Canary \cite{canary1}) Let $f: (S,T) \rightarrow N$ be
a simplicial pre-hyperbolic surface.  Let $\tilde{f}:
\tilde{S}
\rightarrow \HHH$ be a lift of $f$.  Let $v$ be a vertex of $T$, and let
$\tilde{v}$ be one of its lifts.  Let $\tilde{e_1}, \tilde{e_2}, \tilde{e_3}$, and
$\tilde{e_4}$ be edges of $\tilde{T}$ which have $\tilde{v}$ as an endpoint.  If
$\tilde{f}(\tilde{v})$ lies in the tetrahedron spanned by the other endpoints of
$\tilde{f}(\tilde{e_1}), \tilde{f}(\tilde{e_2}),
\tilde{f}(\tilde{e_3})$, and
$\tilde{f}(\tilde{e_4})$, then $v$ satisfies {\it NLSC\/}, and the total angle about
$f(v)$ in $f(S)$ is at least
$2\pi$.
\end{lemma}

Although the ideal vertices are technically not in the domain of a simplicial
hyperbolic surface $f: S \rightarrow N$, we can define an extension map that can be
used to investigate the behavior of the ideal vertices.  Let $S = \HH / \Theta$ be a
surface with punctures, and let $\hat{S} = S \cup \{ideal \ vertices\}$.  Let
$\hat{\mathbb H}^2$ denote $\HH
\cup$ {\it \{parabolic fixed points\}\/} $\subset \bar{\mathbb H}^2$ so that
$\hat{S} =
\hat{\mathbb H}^2 / \Theta$.  Let $N = \HHH / \Gamma$, and let $f:S \rightarrow N$
be a simplicial hyperbolic surface.  Let $\tilde{f}: \HH
\rightarrow \HHH$ be a lift of $f$.  Define the {\it SH extension map of $f$\/} to be
the map $\bar{f}:
\hat{\mathbb H}^2 \rightarrow \bar{\mathbb H}^3$ such that if $v$ is the fixed point
of a parabolic element
$\alpha \in \Theta$, then $\bar{f}(v)$ is the fixed point of $f_*(\alpha) \in
\Gamma$. 

Although $\Theta$ gives $S$ a Riemannian metric, we instead consider the simplicial
hyperbolic structure $\tau$ induced by the map $f$.  Note that
$\Theta$ was needed solely to determine the ideal vertices, and that $\tau$ and the
metric given by $\Theta$ are unrelated.  Under the metric $\tau$,  $\bar{f}$ maps a
ray with endpoint $v \in \hat{\mathbb H}^2$ to a piecewise geodesic ray in
$\bar{\mathbb H}^3$ with endpoint $\bar{f}(v)$.  In fact, $\bar{f}$ maps each lift
of an edge of $\hat{S}$ geodesically and maps each lift of a face of
$\hat{S}$ totally geodesically.  Thus, given a face with (possibly ideal) endpoints
$a, b, c$ in the lift $\tilde{T}$ of the triangulation $T$,
$\bar{f}$ maps this face to a totally geodesic triangle in $\bar{\mathbb H}^3$
spanned by the points
$\bar{f}(a), \bar{f}(b), \bar{f}(c)$.

The following lemma guarantees that the $SH$ extension map of a $\pi_1$-injective
simplicial hyperbolic surface does not map two parabolic fixed points to the same
point.  This property will be useful when altering one type of simplicial hyperbolic
surface into another, as we will see in the next section.


\begin{lemma}
\label{extension} 
Let $f:S \rightarrow N$
be a simplicial hyperbolic surface, and  let $\bar{f}: \hat{\mathbb H}^2 \rightarrow
\bar{\mathbb H}^3$ be the SH extension map of $f$.  If $f$ is
$\pi_1$-injective, then $\bar{f}$ maps no two parabolic fixed points to the same
point. 
\end{lemma}

\begin{proof} Let $\alpha$ and $\beta$ be parabolic elements in $\Theta$
with distinct fixed points.  Suppose for contradiction that
$\bar{f}$ maps the fixed point of
$\alpha$ and the fixed point of $\beta$ to the same point in
$\bar{\mathbb H}^3$.  Then $f_*(\alpha)$ and $f_*(\beta)$ commute in
$\Gamma$.  But $\alpha$ and $\beta$ represent parabolic elements which do not
commute in $\Theta$.  Therefore, $f$ is not
$\pi_1$-injective, but this is a contradiction. So indeed, $\bar{f}$ maps no two
parabolic fixed points to the same point. \end{proof}


%

\subsection{Useful and Practical Simplicial Hyperbolic Surfaces}

Canary does extensive work with a specific type of simplicial hyperbolic surface.  A
simplicial hyperbolic surface $f: S \rightarrow N$ with triangulation $T$ is {\it
useful\/} if $T$ has exactly one internal vertex $v$, and $T$ contains a {\it
distinguished edge\/} which passes through $v$ and is mapped to a closed geodesic.

Given a surface $S$, a
$\pi_1$-injective map
$f:S
\rightarrow N$, and an edge which can be mapped to a closed geodesic, we can
construct a
$\pi_1$-injective useful simplicial hyperbolic surface which is properly homotopic
to $f$.

\begin{lemma} 
\label{makeuseful} Let $S$ be a surface, and let $f: S \rightarrow N$ be a
$\pi_1$-injective map.  Let
$T$ be a suitable triangulation on $S$ with exactly one internal vertex $v$ and an 
edge $\hat{e}$ such that $f(\hat{e} \cup v)$ has a closed geodesic representative.
Then there exists a
$\pi_1$-injective useful simplicial hyperbolic surface $h: (S,T)
\rightarrow N$ with distinguished edge $\hat{e}$ which is properly homotopic to $f$.
\end{lemma}

\begin{proof}  Let $\tilde{S}$ be the universal cover of $S$, and let
$\tilde{v} \in \tilde{S}$ be a lift of $v$.  Let
$\tilde{e}$ be the component of the pre-image of
$(\hat{e} \cup v)$ to $\tilde{S}$ which contains $\tilde{v}$. Let $\tilde{f}:
\tilde{S}
\rightarrow
\HHH$ be a lift of $f$.

By hypothesis, there exists a geodesic axis
$\tilde{f}(\tilde{e})^*$ which is properly homotopic in
$\bar{\mathbb H}^3$ to $\tilde{f}(\tilde{e})$. We can properly homotop $f$ to a map
$g: (S,T) \rightarrow N$ which maps $v$ to a point on the
geodesic covered by $\tilde{f}(\tilde{e})^*$.  Using Bonahon's construction,
$h = F(g,T): (S,T) \rightarrow N$ is a $\pi_1$-injective simplicial pre-hyperbolic
surface.  Because the lift $\tilde{h}: \tilde{S} \rightarrow \HHH$ of
$h$ maps $\tilde{v}$ to a point on $\tilde{f}(\tilde{e})^*$, the
internal vertex
$v$ satisfies
$NLSC$ with respect to $h$ so that
$h$ is a $\pi_1$-injective useful simplicial hyperbolic surface with distinguished
edge $\hat{e}$. By construction, $f$ is properly homotopic to $h$.
\end{proof}

In the case that $S$ is a surface with at least two punctures, we can define an
analogue of a useful simplicial hyperbolic surface.  A simplicial
hyperbolic surface
$f:S \rightarrow N$ with triangulation $T$ is {\it practical\/} if the following
conditions are satisfied:
\begin{enumerate}
\item $S$ is a surface with more than two punctures,
\item $T$ has exactly one internal vertex,
\item $T$ has no doubly ideal edges, and
\item there exists a union of two edges in $T$, called {\it marked edges\/}, each
joining the internal vertex to an ideal vertex, such that the $SH$ extension map
$\bar{f}:
\hat{\mathbb H}^2 \rightarrow \bar{\mathbb H}^3$ maps each component of the pre-image
of the marked edges to a geodesic axis which has endpoints at the parabolic fixed
points associated to the ideal vertices.
\end{enumerate}
    
Note that given a  surface $S$ with more than two punctures and a $\pi_1$-injective
map $f: S
\rightarrow N$, we can construct a $\pi_1$-injective practical simplicial
hyperbolic surface which is properly homotopic to $f$.

\begin{lemma}
\label{makepractical} Let $S$ be a surface with more than two punctures, and let
$T$ be a suitable triangulation on $S$ with exactly one internal vertex $v$. Let
$e_1$ and $e_2$ be two edges in $T$ joining the internal vertex
$v$ to two distinct punctures of $S$.  If $f: S \rightarrow N$ is $\pi_1$-injective,
then there exists a
$\pi_1$-injective practical simplicial hyperbolic surface $h: (S,T)
\rightarrow N$ with marked edges $e_1$ and $e_2$ which is properly homotopic to $f$.

\end{lemma}

\begin{proof}  Let $\tilde{S}$ be the universal cover of $S$,
and let
$\tilde{v} \in \tilde{S}$ be a lift of $v$. Let $\tilde{e}_1$ and $\tilde{e}_2$ be
the components  of the pre-images of $(e_1 \cup v)$ and $(e_2 \cup v)$ to $\tilde{S}$
which intersect at
$\tilde{v}$. 
 
Using the notation of the $SH$ extension map, because $f$ is $\pi_1$-injective, we
can construct a well-defined extension map
$\bar{f}: \hat{\mathbb H}^2 \rightarrow \bar{\mathbb H}^3$ such that if $w$ is the
fixed point of a parabolic element
$\alpha \in \Theta$, then $\bar{f}(w)$ is the fixed point of $f_*(\alpha) \in
\Gamma$.  By an argument similar to that in Lemma \ref{extension}, we can show that
whereas $\bar{f}(\tilde{e}_1)$ and $\bar{f}(\tilde{e}_2)$ share
$\bar{f}(\tilde{v})$ as an endpoint, the other endpoints of
$\bar{f}(\tilde{e}_1)$ and $\bar{f}(\tilde{e}_2)$ do not coincide.  Hence, there
exists a geodesic axis 
$\tilde{f}(\tilde{e}_1 \cup \tilde{e}_2)^*$ which is properly homotopic in
$\bar{\mathbb H}^3$ to $\tilde{f}(\tilde{e}_1
\cup \tilde{e}_2)$. 

We can properly homotop $f$  to a map $g: (S,T) \rightarrow N$ which maps $v$ to a
point on the geodesic covered by $\tilde{f}(\tilde{e}_1 \cup
\tilde{e}_2)^*$.  Using Bonahon's construction,
$h = F(g,T): (S,T) \rightarrow N$ is a $\pi_1$-injective simplicial pre-hyperbolic
surface.  Because the lift $\tilde{h}: \tilde{S} \rightarrow \HHH$ of
$h$ maps $\tilde{v}$ to a point on $\tilde{f}(\tilde{e}_1 \cup \tilde{e}_2)^*$, the
internal vertex $v$ satisfies
$NLSC$ with respect to $h$ so that
$h$ is a $\pi_1$-injective practical simplicial hyperbolic surface with marked edges
$e_1$ and $e_2$.  By construction,
$f$ is properly homotopic to $h$.
\end{proof}
 
Both useful and practical simplicial hyperbolic surfaces have the property that
their images lies in the convex core as shown in the following lemma.

\begin{lemma}
\label{shsincore} Let $f:(S,T) \rightarrow N$ be a simplicial hyperbolic surface. 
If $f$ maps every internal vertex of $T$ into the convex core, then $f(S) \subset
C(N)$.  In particular, if $f$ is a useful or practical simplicial hyperbolic
surface, then $f(S) \subset C(N)$.
\end{lemma}

\begin{proof}  
Consider the $SH$ extension map of $f$, $\bar{f}: \hat{\mathbb H}^2
\rightarrow \bar{\mathbb H}^3$.
By definition, $\bar{f}$ maps the parabolic fixed points of
$\pi_1(S)$ into $\Lambda_{\pi_1(N)}$. By hypothesis $f$ maps all the internal
vertices of $T$ into $C(N)$.  Thus, $\bar{f}$ maps all lifts of internal vertices
into $CH(\Lambda_{\pi_1(N)}) \subset \HHH$.  Recall that $\bar{f}$ maps a face with 
endpoints $a, b, c$ to a geodesic triangle in $\bar{\mathbb H}^3$ spanned by the
points
$\bar{f}(a), \bar{f}(b), \bar{f}(c)$.  Because $CH(\Lambda_{\pi_1(N)})$ is convex,
we can conclude that $\bar{f}(\hat{\mathbb H}^2) \subset (CH(\Lambda_{\pi_1(N)})
\cup
\Lambda_{\pi_1(N)})$, and hence $f(S) \subset C(N)$.

Suppose $f:S \rightarrow N$ is a useful simplicial hyperbolic surface with
distinguished edge $\hat{e}$.  Then $f$ maps $\hat{e}$ to a closed geodesic. 
Because all closed geodesics lie in the convex core, and the sole internal vertex
$v$ lay on $\hat{e}$, we know that $f(v) \in C(N)$.  Then by the above argument,
$f(S)
\subset C(N)$.

Suppose $f:S \rightarrow N$ is a practical simplicial hyperbolic surface with marked
edges $\hat{e}_1$ and $\hat{e}_2$.  By definition,
$\bar{f}$ maps each component of the pre-image of ($\hat{e}_1 \cup v \cup \hat{e}_2$)
to
$\tilde{S}$ to a geodesic that has endpoints at parabolic fixed
points associated to the ideal vertices contained in the marked edges.  Thus,
$\bar{f}$ maps each component of the pre-image of ($\hat{e}_1
\cup v \cup \hat{e}_2$) to $\tilde{S}$ to a geodesic axis in 
$(CH(\Lambda_{\pi_1(N)})
\cup
\Lambda_{\pi_1(N)})$.  Since the lift $\tilde{v}$ of the sole internal vertex
lay on this axis,
$\bar{f}(\tilde{v}) \in CH(\Lambda_{\pi_1(N)})$, and hence $f(v) \in C(N)$.  Then by
the above argument,
$f(S)
\subset C(N)$. \end{proof}

\section{Continuous Families of Simplicial Hyperbolic Surfaces}

One of the basic tools in proving the main theorem is the theory of continuous
families of simplicial hyperbolic surfaces.  A similar theory for pleated
surfaces was first developed by Thurston \cite{thurston}, and this theory was
developed for simplicial hyperbolic surfaces by Canary \cite{canary1}.  This
tool provides a way of interpolating between pairs of simplicial hyperbolic
surfaces that differ in an elementary way. 

A basic continuous family of simplicial pre-hyperbolic surfaces can be
constructed by dragging a map along a path.  This construction will be used
throughout this section.  We start with a simplicial pre-hyperbolic map $f: (S,
\tau) \rightarrow N$, and a path $\gamma: [0, 1] \rightarrow N$ such that for some
vertex $v \in T$, we have $f(v) =
\gamma(0)$.  We will construct a homotopy $H: S \times [0, 1] \rightarrow N$
such that $H_t = H|_{S \times \{t\}}$ is a simplicial pre-hyperbolic surface for
all $t \in [0, 1]$, $H_t(v) = \gamma(t)$, and $H_0 = f$.  We will say that $H$ is a
continuous family of simplicial pre-hyperbolic surfaces obtained by {\it
dragging $f$ along $\gamma$}.

To construct $H$, we first choose $\delta < inj_{(S, \tau)}(v)$ such that no
other vertex of $T$ lies within the ball $B_\delta(v)$ in $(S, \tau)$ of radius
$\delta$ based at $v$.  Construct $R: S \times [0, 1] \rightarrow N$ with the
following properties for all $t \in [0,1]$: \  (1) $R(x, t)$ agrees with $f$ in the
complement of $B_\delta(v)$, (2) $R(v, t) = \gamma(t)$, and (3)
$R(B_\delta(v), t) \subset f(B_\delta(v)) \cup
\gamma([0, t])$.  Then
$H(x, t) = F(R(x, t), T)$ is a continuous family of simplicial pre-hyperbolic
surfaces, where
$F$ is the map from Bonahon's construction described in the previous section which
takes each face of
$T$ of
$R(x,t)$ to a totally geodesic triangle spanned by the images of its edges. 
Note that the resulting continuous family is independent of the choice of
$R$, because under the same restrictions, the homotopy classes of the images of
the edges remain the same.  Also note that if one can show that the internal 
vertex $v$ satisfies $NLSC$ with respect to $H_t$ for $t \in [0,1]$, then $H$ will be
a continuous family of simplicial hyperbolic surfaces.

Given a useful simplicial hyperbolic surface $f$, Canary has
shown how to construct a continuous family of simplicial hyperbolic surfaces in a
variety of situations: \ changing the distinguished edge, performing an elementary
move, and shrinking an accidental parabolic.   These are Lemmas 5.2 to 5.4 in
\cite{canary1}. We state these lemmas here for reference, and in addition, we show
that the image of each of these continuous families lies in $C(N)$.

For the remainder of this section, for a homotopy $H: S \times [0,1] \rightarrow N$,
the notation
$H_t = H|_{S \times \{t\}}$ will be used.

\begin{lemma}
\label{edge} (Lem 5.2, Canary \cite{canary1})  Let
$h: (S,T) \rightarrow N$ be a useful simplicial hyperbolic surface with 
distinguished edge
$\hat{e}$.  Let
$\bar{e}$ be another edge of $T$ such that $h(\bar{e} \cup v)$ has a closed
geodesic representative
$h(\bar{e} \cup v)^*$.  Then there exists a continuous family $J: S \times [0,2]
\rightarrow N$ of simplicial hyperbolic surfaces with only one internal vertex (which
satisfies {\it NLSC\/}) joining $h$ to a useful simplicial hyperbolic surface
with associated triangulation $T$ and distinguished edge $\bar{e}$. 
Furthermore,
$J(S,[0,2]) \subset C(N)$.
\end{lemma}

\begin{proof} We begin by outlining the construction of
the continuous family given by Canary, then we will show that its image lies in
the convex core.  Let
$W: S^1
\times [0,1]
\rightarrow N$ be a homotopy between $h(\bar{e} \cup v)$ and its geodesic
representative
$h(\bar{e} \cup v)^*$.  Let $\tilde{S}$ be the universal cover of $S$, and let
$\tilde{v}
\in \tilde{S}$ be a lift of $v$.  Let
$\alpha_1$ and
$\alpha_2$ be the components of the pre-images of $(\hat{e} \cup v)$ and $(\bar{e}
\cup v)$ to
$\tilde{S}$ which intersect at $\tilde{v}$.  Lift both
$h$ and
$W$ to maps $\tilde{h}: \tilde{S} \rightarrow \HHH$ and $\tilde{W}:
\Rr
\times [0,1] \rightarrow \HHH$ such that $\tilde{W}(\Rr,0) =
\tilde{h}(\alpha_2)$.  Let $\gamma_1 = \tilde{h}(\alpha_1)$, and let
$\gamma_2 = \tilde{W}(\Rr,1)$.  Let $\tilde{R}$ be the unique
common perpendicular joining $\gamma_1$ to $\gamma_2$.  Let $x = \tilde{R} \cap
\gamma_1$, and let $\tilde{P}$ be a geodesic arc in $\gamma_1$ joining
$\tilde{h}(\tilde{v})$ to $x$.  Let $R$ and $P$ be the projections to $N$ of
$\tilde{R}$ and $\tilde{P}$ respectively.  By dragging $h$ along $P$, we obtain a
continuous family of simplicial hyperbolic surfaces $G: S
\times [0,1]
\rightarrow N$.  In fact, this is a continuous family  of useful simplicial
hyperbolic surfaces because $\tilde{P}$ is a geodesic arc contained in a lift of
$h(\hat{e} \cup v)$ and hence $G$ maps $(\hat{e} \cup v) \times \{t\}$  to a
closed geodesic for $t \in [0,1]$. By dragging
$G_1$ along
$R$, we obtain a continuous family of simplicial hyperbolic surfaces $H:S \times [1,2]
\rightarrow N$ such that
$H_2(v)$ lies on the closed geodesic in the homotopy class of $h(\bar{e} \cup
v)$.  Then $H_2$ is a useful simplicial hyperbolic surface with associated
triangulation $T$ and distinguished edge $\bar{e}$.  Finally, $J$ is obtained by
concatenating $G$ and $H$. This completes the sketch of the argument given by Canary.

Because $J_t$  are useful simplicial hyperbolic surfaces for $t \in [0,1]$,  their
image automatically lies in $C(N)$.  Since the sole internal vertex of each
simplicial hyperbolic surface in the family $J_t$ for $t \in [1,2]$ is dragged along
the arc 
$R$, by Lemma
\ref{shsincore}, it suffices to show that this arc lies in the convex core.
 
%
Let $\bar{e}$ be the distinguished edge
of $J_2: S \rightarrow N$.  Here, $\tilde{R}$ is defined to be the unique 
common perpendicular between two components of lifts of
$h(\hat{e} \cup v)$ and $h(\bar{e} \cup v)$, each of which is a collection of
axes of hyperbolic  elements.  Then $\tilde{R}$ is a geodesic that joins two
axes of  hyperbolic elements.  Because $CH(\Lambda_\Gamma)$ is convex, we know that
$\tilde{R}$ lies in
$CH(\Lambda_\Gamma)$, and hence $R \subset C(N)$.  Thus, the image of $J$ lies
in the convex core of $N$. \end{proof}

\begin{lemma}
\label{elem move}  (Lem 5.3, Canary \cite{canary1})  Let 
$h: (S,T) \rightarrow N$ be a useful simplicial hyperbolic surface with 
distinguished edge $\hat{e}$.  Let $e_1, e_2, e_3$ and $e_4$ bound a quadrilateral
in $T$ with diagonal $e_5 \neq \hat{e}$.  Let $e_6$ be the other diagonal of this
quadrilateral, and let $T^\prime$ be the triangulation obtained by replacing $e_5$
with $e_6$.  Assume that $h(e_6 \cup v)$ is homotopically non-trivial.  Then there
exists a continuous family of simplicial hyperbolic surfaces $K: S \times [0,1]
\rightarrow N$ with at most two internal vertices (each having property {\it
NLSC\/}) joining
$h$ to a useful simplicial hyperbolic surface with associated
triangulation
$T^\prime$ and distinguished edge $\hat{e}$.  Furthermore, $K(S,[0,1]) \subset
C(N)$.
\end{lemma}

\begin{proof} We begin by outlining the construction of
the continuous family given by Canary, then we will show that its image lies in
the convex core. Let $\tilde{S}$ be the universal cover of $S$, and let
$\tilde{e}_1$, $\tilde{e}_2$, $\tilde{e}_3$, $\tilde{e}_4$, $\tilde{e}_5$, and
$\tilde{e}_6$ be lifts of $e_1$, $e_2$, $e_3$,
$e_4$,
$e_5$, and
$e_6$ to
$\tilde{S}$ such that $\tilde{e}_1$, $\tilde{e}_2$, $\tilde{e}_3$, and
$\tilde{e}_4$ form a quadrilateral with diagonals $\tilde{e}_5$ and
$\tilde{e}_6$.  Let $\tilde{h}:
\tilde{S} \rightarrow \HHH$ be a lift of $h$.  Let
$\tilde{h}(\tilde{e}_6)^*$ denote the geodesic arc joining the endpoints of
$\tilde{h}(\tilde{e}_6)$.  Notice that $\tilde{h}(\tilde{e}_1)$,
$\tilde{h}(\tilde{e}_2)$, $\tilde{h}(\tilde{e}_3)$, $\tilde{h}(\tilde{e}_4)$,
$\tilde{h}(\tilde{e}_5)$, and $\tilde{h}(\tilde{e}_6)$ form a tetrahedron $X$ in
$\HHH$.  Let $v^\prime$ denote the intersection of $e_5$ and $e_6$, and
let $\tilde{v}^\prime = \tilde{e}_5 \cap \tilde{e}_6$.  Let $\tilde{R}$ be the
geodesic in $X$ joining $\tilde{h}(\tilde{v}^\prime)$ and
$\tilde{h}(\tilde{e}_6)^*$.  We obtain a new triangulation $\bar{T}$ of $S$ by
adding the edge $e_6$ to $T$ so that $\bar{T}$ now has an additional vertex at
$v^\prime = e_5
\cap e_6$.

Let $R$ be the projection to $N$ of $\tilde{R}$.  By dragging $h$ along $R$ we
obtain a continuous family of simplicial hyperbolic surfaces
$K: S \times [0,1] \rightarrow N$ with associated triangulation $\bar{T}$. 
After removing the vertex $v^\prime$ and the edge $e_5$,
$K_1$ is a useful simplicial hyperbolic surface with associated triangulation
$T^\prime$.  This completes the sketch of the argument given by Canary.

Note that $K_0$ and
$K_1$ are useful simplicial hyperbolic surfaces, and hence their image lies in
the convex core.  The triangulation $\bar{T}$ has two internal vertices, $v$ and
$v^\prime$.  To show that $K(S,[0,1]) \subset
C(N)$, by Lemma \ref{shsincore}, it suffices to show that
for $t \in [0,1]$, the two vertices $v$ and $v^\prime$ are mapped  into the convex
core.

Note that for all $t \in [0,1]$, $K_t(\hat{e})$ is a fixed closed geodesic in $N$. 
Because $v \in \hat{e}$, the vertex $v$ is mapped into $C(N)$ for all $t$.  In the
construction, the vertex
$v^\prime$ is dragged along a geodesic arc $R$, so it suffices to show that $R
\subset C(N)$.  Since the endpoints of $\tilde{R}$ are in
$CH(\Lambda_{\pi_1(N)})$, it is clear that $R \subset C(N)$.  \end{proof}

The following lemma extends Canary's version of the lemma by not only showing that
the image of the continuous family lies in $C(N)$, but also by showing that the
image lifts to a $\delta$-neighborhood of the convex core of a cover.

\begin{lemma}
\label{length to 0} (Extension of Lem 5.4, Canary \cite{canary1}) Let 
$h: (S,T)
\rightarrow N$ be a useful simplicial hyperbolic surface with distinguished edge
$\hat{e}$.  Let
$e^\prime \neq \hat{e}$ be an edge of $T$, and let $v$ be
the sole internal vertex of
$T$ such that
$\gamma = (e^\prime \cup v)$ is an accidental parabolic.  Then:
\begin{enumerate}
\item there exists a continuous family of simplicial hyperbolic surfaces $H: S
\times [0,1) \rightarrow N$ with only one vertex (which satisfies {\it NLSC\/})
such that as $t \rightarrow 1$, the length of $H_t(e^\prime \cup v)$ decreases
and converges to $0$;
\item $H(S,[0,1)) \subset C(N)$; 
\item $H(S,[0,1))$ has compact closure in  $N - \bigcup P_i$, where \{$P_i$\}
is the set of rank one cusps of $N - \NOE$; and
\item for a component $S_i$ of $S - \gamma$, $H_t(S_i)$ lifts to
$C_{\delta(t)}(M_i)$, where $M_i = \HHH / h_*(\pi_1(S_i))$ and where 
$\delta: [0,1) \rightarrow \Rr$ is a function such that as $t \rightarrow 1$,
$\delta(t)
\rightarrow 0$.  Furthermore, for $t \in [0,1)$, the map $H_t: S_i \rightarrow
N$ has fewer
accidental parabolics than $H_0: S \rightarrow N$.
\end{enumerate}
\end{lemma}
%

\begin{proof} First let us give Canary's proof of the
existence of a continuous family of simplicial hyperbolic surfaces.  Let
$W : S^1 \times [0, \infty) \rightarrow N$ be a proper homotopy of $h(e^\prime
\cup v)$ into a cusp of $N$.   Let $\tilde{S}$ be the universal cover of $S$.  Let
$\tilde{v}$ be a lift of $v$, and let $\alpha$ be the component of the pre-image of
$(e^\prime \cup v)$  to $\tilde{S}$ that contains
$\tilde{v}$.  Lift both $h$ and $W$ to
maps
$\tilde{h}: \tilde{S}
\rightarrow
\HHH$ and $\tilde{W}: \Rr \times [0, \infty) \rightarrow \HHH$ such that
$\tilde{h}(\alpha) = \tilde{W}(\Rr,0)$.

Let $t \in \Rr$ such that $\tilde{h}(\tilde{v}) = \tilde{W}(t, 0)$. 
Let $\tilde{R}$ be a geodesic ray properly homotopic to $\tilde{W}(t,[0,
\infty))$, relative to $\tilde{W}(t, 0)$.  Note that $\tilde{R}$ has an endpoint
in the sphere at infinity that corresponds to the fixed point of a parabolic
element. Drag $h$ along the projection $R$ of
$\tilde{R}$ to obtain a continuous family of simplicial pre-hyperbolic surfaces
$H: S \times [0, 1) \rightarrow N$. 

For the remainder of the proof, we will consider the upper half plane model of $\HHH$,
normalized so that $\tilde{R}$ is a geodesic ray with endpoint at infinity. 

To check that we have a continuous family of simplicial hyperbolic surfaces, we
will utilize Lemma \ref{tetrahedron} to verify that the internal vertex $v$
satisfies {\it NLSC\/} with respect to $H_t$ for $t \in [0,1)$.  Recall that the
distinguished edge $(\hat{e} \cup v)$ is mapped to a closed geodesic.   Let $\beta$
be the component of the pre-image of $(\hat{e} \cup v)$ that intersects $\tilde{v}$. 
Let
$\gamma_1 =
\tilde{h}(\beta)$.  Let $g_1$ be the hyperbolic isometry with axis $\gamma_1$,
and let $g_2$ be the parabolic isometry whose fixed point is the endpoint of
$\tilde{R}$.  There are
four edges of $\tilde{H}_t(\tilde{T})$ which have one endpoint at
$\tilde{H}_t(\tilde{v})$ have other endpoints at $g_1(\tilde{H}_t(\tilde{v})),
g_{1}^{-1}(\tilde{H}_t(\tilde{v})), g_2(\tilde{H}_t(\tilde{v}))$, and
$g_{2}^{-1}(\tilde{H}_t(\tilde{v}))$.

By construction $\tilde{H}_t(\tilde{v})$ lies on $\tilde{R}$.  Note that 
$\tilde{H}_t(\tilde{v})$ lies on the same horocycle as
$g_2(\tilde{H}_t(\tilde{v}))$ and $g_{2}^{-1}(\tilde{H}_t(\tilde{v}))$. 
Therefore, the geodesic spanned by $g_2(\tilde{H}_t(\tilde{v}))$ and
$g_{2}^{-1}(\tilde{H}_t(\tilde{v}))$ intersects $\tilde{R}$ above
$\tilde{H}_t(\tilde{v})$.  Now because $\tilde{v}$,
$g_1(\tilde{H}_t(\tilde{v}))$, and $g_{1}^{-1}(\tilde{H}_t(\tilde{v}))$ are all
exactly a radial distance $r$ from the axis of $g_1$, we know that the geodesic
joining
$g_1(\tilde{H}_t(\tilde{v}))$ and $g_{1}^{-1}(\tilde{H}_t(\tilde{v}))$ lies
inside the $r$-neighborhood of the axis of $g_1$.  Then this geodesic will
intersect $\tilde{R}$ below
$\tilde{H}_t(\tilde{v})$.  Hence, $\tilde{H}_t(\tilde{v})$ lies in the
tetrahedron spanned by the four other endpoints.  By Lemma \ref{tetrahedron},
the vertex $v$ satisfies {\it NLSC\/} with respect to $H_t$ for $t \in [0,1)$.  Thus,
$H$ is a continuous family of simplicial hyperbolic surfaces. 

Now let us show that the length of $H_t(\gamma)$ converges to 0.  As $t
\rightarrow 1$, $h$ is dragged along $R$, and for each $t$, $H_t$ maps $\gamma$ to a
 geodesic arc joined at its endpoints.  In the universal cover $\HHH$, for each
$t$, the component of the pre-image of $H_t(\gamma)$ which intersects $\tilde{R}$ is a
piecewise-geodesic line.  Within the fundamental domain in $\HHH$ containing
$\tilde{R}$, the pre-image of
$H_t(\gamma)$ is a geodesic arc of bounded Euclidean
length whose height approaches infinity as
$t$ approaches 1.  So the hyperbolic length of
$H_t(\gamma)$ converges to 0.  This completes the sketch of the argument given by
Canary.

Now we will show that the image of the continuous family $H$ lies
in the convex core.  Since $h: S \rightarrow N$ is a useful simplicial hyperbolic
surface, by Lemma \ref{shsincore}, its image lies in $C(N)$.   Because
the vertex is dragged along $R$, by Lemma \ref{shsincore}, it suffices to show that
$R \subset C(N)$.

One endpoint of $\tilde{R}$ is a lift of $h(v)$.  Since the distinguished
edge $\hat{e}$ passes through the vertex $v$, 
$h(v)$ lies in the convex core of $N$.  The other endpoint of $\tilde{R}$ is a
parabolic fixed point.  Because $R$ is a geodesic ray, we know that
$\tilde{R} \subset CH(\Lambda_\Gamma)$, and hence $R \subset C(N)$.  Thus,
$H(S,[0,1)) \subset C(N)$.

Now consider the image $H(S,[0,1))$.  Each triangle in the finite triangulation $T$ on
$S$ is mapped to a geodesic triangle in $N$ with sole internal vertex $h(v) \in R$. 
Let $H(A,0)$ be the image in $N$ of a triangle $A$ in $T$.  Let $\tilde{H}: S \times
[0,1) \rightarrow \HHH$ be a lift of $H$.  Then
$\tilde{H}(A,0) \subset \HHH$ consists of triangles 
whose vertices lie on three geodesic rays which are conjugates of $\tilde{R}$.   As $t
\rightarrow 1$, we drag the vertices of $\tilde{H}(A,0)$ along the conjugates of
$\tilde{R}$ towards  conjugates of
$\infty$.  In doing so, each component of $\tilde{H}(A,0)$ converges to an ideal
triangle with vertices at three conjugates of $\infty$.  Thus, each component of
$\tilde{H}(A,[0,1))$ lies in a  totally geodesic triangular prism with three endpoints
at lifts of
$h(v)$ and the three remaining (ideal) endpoints at conjugates of $\infty$.  Hence
$\tilde{H}(S,[0,1))$ lies in a finite collection of such geodesic prisms.  
Here, the ideal endpoints of the prisms are conjugates of
$\infty$ which are the fixed points of the parabolic elements associated to rank one
cusps
$P_i$ of
$N - \NOE$.  So we can conclude that $H_2(S,[0,1))$ has compact closure in $N -
\bigcup P_i$.

Now let $S_i$ be a component of $S - \gamma$.  Here, we consider $S_i$ to be a
surface with an (additional) ideal vertex associated to $\gamma$.  Note that the
restriction to $S_i$ of the triangulation of $S$ forms a system of curves on
$S_i$ that does not constitute a triangulation, because there are homotopic
edges.  However, for $t \in [0,1)$, $H_t$ still maps each edge of the system of
curves to a geodesic arc and maps each face geodesically.  

%

Let us investigate the behavior of the map $H_t: S \rightarrow N$ restricted to
the subsurface $S_i$.  Let $\bar{H}_t: \hat{\mathbb H}^2 \rightarrow
\bar{\mathbb H}^3$ be the $SH$ extension
map of $H_t$.  We will consider the restriction of $\bar{H}_t$ to the
component $\tilde{S}_i$  of the preimage of $S_i$ in $\tilde{S} \subset \hat{\mathbb
H}^2$.  
We can lift the triangulation $T$ on $S$ to a triangulation $\tilde{T}$ on
$\tilde{S}$ and consider the restriction $\tilde{T}_i$ of the triangulation
$\tilde{T}$ to
$\tilde{S}_i$.

Let $M_i = \HHH / h_*(\pi_1(S_i))$.     Because $H_t$ maps
each edge to a geodesic arc and maps each face geodesically, by an argument similar
to that in Lemma \ref{shsincore}, to show that $H_t(S_i)$ lifts to
$C_{\delta(t)}(M_i)$,  it suffices to show that $\bar{H}_t$ maps all internal
vertices of $\tilde{T}_i$ into $CH_{\delta(t)}(\Lambda_{h_*(\pi_1(S_i))})$. 
Furthermore, since
$\bar{H}_t$ respects the group action of $\pi_1(S_i)$,
it suffices to show that
$\bar{H}_t(\tilde{v})$ maps $\tilde{v}$ into
$CH_{\delta(t)}(\Lambda_{h_*(\pi_1(S_i))})$.

Because
$C(M_i)$ is nonempty, the convex hull of $\Lambda_{h_*(\pi_1(S_i))}$ must
contain a geodesic ray $\tilde{Y}$ with endpoint at infinity.  Let
$\bar{H}_t(\tilde{v})$ be a point of height $r_t$ on
$\tilde{R}$.  Recall that $\tilde{R}$ is also a geodesic ray with endpoint at
infinity.  Then for a fixed $t$, there exists $\delta(t) > 0$ such that
$d(\bar{H}_t(\tilde{v}), \tilde{Y}) < \delta(t)$ so that $\bar{H}_t(\tilde{v})
\in CH_{\delta(t)}(\Lambda_{h_*(\pi_1(S_i))})$.  Therefore,
$H_t(S_i)$
lifts to $C_{\delta(t)}(M_i)$.

Because
$\tilde{R}$ is a fixed Euclidean distance from $\tilde{Y}$, as
$t \rightarrow 1$, $r_t \rightarrow \infty$, and hence
$d(\tilde{v}_t,\tilde{Y}) \rightarrow 0$.  Therefore, as $t \rightarrow 1$,
$\delta(t)
\rightarrow 0$. 
Furthermore, because we have eliminated the accidental parabolic associated to
$\gamma$, for $t \in [0,1)$, the map $H_t: S_i \rightarrow N$ has fewer accidental
parabolics than
$H_0: S \rightarrow N$.
\end{proof}

Now we will prove the analogues of Lemmas \ref{edge}, \ref{elem move}, and \ref{length
to 0} for practical simplicial hyperbolic surfaces.   These lemmas alter the
marked edges, perform an elementary move, and shrink an accidental parabolic.  
Furthermore, the images of these continuous families lie in $C(N)$.

\begin{lemma}
\label{3 edge}  Let $h: (S,T)
\rightarrow N$ be a $\pi_1$-injective practical simplicial hyperbolic surface with
 marked edges $\hat{e}_1$ and $\hat{e}_2$.  Let
$\hat{e}_3$ and $\hat{e}_4$ be two edges of $T$ joining the sole internal vertex to
two ideal vertices where at least one of $\hat{e}_3$ and $\hat{e}_4$ is different
from $\hat{e}_1$ and $\hat{e}_2$.  Then there exists a continuous family $J: S
\times [0,2] \rightarrow N$ of simplicial hyperbolic surfaces with only one
internal vertex (which satisfies {\it NLSC\/}) joining $h$ to a practical simplicial
hyperbolic surface with associated triangulation $T$ and marked edges
$\hat{e}_3$ and $\hat{e}_4$.  Furthermore, $J(S,[0,2]) \subset C(N)$.
\end{lemma}

\begin{proof} The proof of this lemma mimics that
of Lemma \ref{edge}; we give a brief outline here.  Let $\tilde{S}$ be the universal
cover of
$S$, and let $\tilde{v} \in \tilde{S}$ be a lift of $v$.  Let
$\alpha_1$ and
$\alpha_2$ be the components of the pre-images of
$(\hat{e}_1 \cup \hat{e}_2 \cup v)$ and $(\hat{e}_3 \cup \hat{e}_4 \cup v)$ to
$\tilde{S}$ which intersect at a point $\tilde{v}$. Lift
$h$ to a map
$\tilde{h}: \tilde{S} \rightarrow \HHH$.  Let $\gamma_1 =
\tilde{h}(\alpha_1)$.

Because $h$ is a $\pi_1$-injective simplicial hyperbolic surface, Lemma
\ref{extension} guarantees that the $SH$ extension $\bar{h}$
maps no two parabolic fixed points to the same point.  Thus there exists a
geodesic axis $\gamma_2$ in the proper homotopy class of $\bar{h}(\alpha_2)$.

Take $\tilde{R}$ to be the unique common perpendicular
joining $\gamma_1$ to $\gamma_2$.  Let $x = \tilde{R} \cap \gamma_1$, and let
$\tilde{P}$ be a geodesic arc in $\gamma_1$ joining $\tilde{h}(\tilde{v})$ to
$x$.  Let $R$ and
$P$ be the projections of $\tilde{R}$ and $\tilde{P}$ respectively.  We construct
$J$ by dragging $h$ along $P$, and then $R$.  The remainder of the proof is
similar to that of Lemma \ref{edge}.  \end{proof}

\begin{lemma}
\label{3 elem move}  Let $h:(S,T) \rightarrow N$ be a practical simplicial
hyperbolic surface with marked edges
$\hat{e_1}$ and $\hat{e_2}$.  Let $e_1, e_2, e_3$, and $e_4$ bound a quadrilateral
$T$ with diagonal $e_5$ where $e_5$ is not a marked edge.  Let $e_6$ be the other
diagonal of this quadrilateral, and let $T^\prime$ be the triangulation obtained
by replacing $e_5$ with $e_6$.  Assume that $h(e_6 \cup v)$ is homotopically
non-trivial.  Then there exists a continuous family of simplicial hyperbolic
surfaces $K:S \times [0,1] \rightarrow N$ with at most two internal vertices (each
having property {\it NLSC\/}) joining $h$ to a practical simplicial hyperbolic
surface with associated triangulation $T^\prime$ and marked edges $\hat{e_1}$ and
$\hat{e_2}$.  Furthermore, $K(S,[0,1]) \subset C(N)$.
\end{lemma}

\begin{proof} 
The proof of this lemma mimics that of Lemma \ref{elem move}; we give a brief
outline here.  Let
$\tilde{e}_1$,
$\tilde{e}_2$,
$\tilde{e}_3$,
$\tilde{e}_4$,
$\tilde{e}_5$, and
$\tilde{e}_6$ be lifts of $e_1$, $e_2$, $e_3$,
$e_4$,
$e_5$, and
$e_6$ to
$\tilde{S}$ such that $\tilde{e}_1$, $\tilde{e}_2$, $\tilde{e}_3$, and
$\tilde{e}_4$ form a quadrilateral with diagonals $\tilde{e}_5$ and
$\tilde{e}_6$.  Let $\tilde{h}: \tilde{S} \rightarrow \HHH$ be a lift of $h$. 
Let $\tilde{h}(\tilde{e}_6)^*$ denote the geodesic arc joining the endpoints of
$\tilde{h}(\tilde{e}_6)$.  Notice that $\tilde{h}(\tilde{e}_1)$,
$\tilde{h}(\tilde{e}_2)$, $\tilde{h}(\tilde{e}_3)$, $\tilde{h}(\tilde{e}_4)$,
$\tilde{h}(\tilde{e}_5)$, and $\tilde{h}(\tilde{e}_6)$ form a tetrahedron $X$ in
$\HHH$.  Let $v^\prime$ denote the intersection of $e_5$ and $e_6$, and
let $\tilde{v}^\prime = \tilde{e}_5 \cap \tilde{e}_6$.  Let $\tilde{R}$ be the
geodesic in $X$ perpendicular to $\tilde{h}(\tilde{e}_6)^*$ and passing
through $\tilde{h}(\tilde{v}^\prime)$.  We obtain a new triangulation
$\bar{T}$ of
$S$ by adding $e_6$ to $T$ which now has an additional vertex at $v^\prime = e_5
\cap e_6$.

Let $R$ be the projection to $N$ of $\tilde{R}$.  By dragging $h$ along $R$ we
obtain a continuous family of simplicial pre-hyperbolic surfaces
$K: S \times [0,1] \rightarrow N$ with associated triangulation $\bar{T}$. 
The remainder of the proof is similar to that of Lemma \ref{elem move}.
\end{proof}

\begin{lemma}
\label{3 length to 0} Let 
$h: (S,T) \rightarrow N$ be a practical simplicial hyperbolic surface with
marked edges
$\hat{e_1}$ and $\hat{e_2}$.  Let $e^\prime$ be an edge of $T$ which is different
from $\hat{e_1}$ and $\hat{e_2}$, and let
$v$ be the sole internal vertex of
$T$ such that $\gamma = (e^\prime \cup v)$ is an accidental parabolic.  Then:

\begin{enumerate}
\item there exists a continuous family of simplicial hyperbolic surfaces $H: S
\times [0,1) \rightarrow N$ with only one internal vertex (which satisfies {\it
NLSC\/}) such that as $t \rightarrow 1$, the length of $H_t(e^\prime \cup v)$
decreases and converges to $0$;

\item $H(S,[0,1)) \subset C(N)$;
\item $H(S,[0,1))$ has compact closure in  $N - \bigcup P_i$, where \{$P_i$\}
is the set of rank one cusps of $N - \NOE$; and
\item for a component $S_i$ of $S - \gamma$, $H_t(S_i)$ lifts to
$C_{\delta(t)}(M_i)$, where $M_i = \HHH / h_*(\pi_1(S_i))$ and where 
$\delta: [0,1) \rightarrow \Rr$ is a function such that as $t \rightarrow 1$,
$\delta(t)
\rightarrow 0$.  Furthermore, for $t \in [0,1)$, the map $H_t: S_i \rightarrow
N$ has fewer
accidental parabolics than $H_0: S \rightarrow N$.

\end{enumerate}
\end{lemma}

\begin{proof} The proof of this lemma is virtually identical to that of Lemma
\ref{length to 0}. \end{proof}

Sometimes it will be necessary to construct a continuous family  that interpolates
between a useful  and a practical simplicial hyperbolic
surface.

\begin{lemma}
\label{usefulpractical}  Let $S$ be a surface with at least two punctures, and let
$h:(S,T) \rightarrow N$ be a $\pi_1$-injective useful simplicial hyperbolic surface
with distinguished edge $\hat{e}$.  Let $\hat{e}_1,
\hat{e}_2$ be two edges in $T$ that join the internal vertex to two distinct
ideal vertices. Then there exists a continuous family of simplicial hyperbolic
surfaces
$K:S \times [0,1]
\rightarrow N$ with only one internal vertex (which satisfies {\it NLSC\/})
joining $h$ to a practical simplicial hyperbolic surface with
associated triangulation $T$ and marked edges $\hat{e_1}, \hat{e_2}$. 
Furthermore,
$K(S,[0,1]) \subset C(N)$.
\end{lemma}

\begin{proof}  The proof of this lemma
also mimics that of Lemma \ref{edge}.  We outline the argument here.  Let
$\tilde{S}$ be the universal cover of $S$, and let
$\alpha_1$ and $\alpha_2$ be components of pre-images of $(\hat{e} \cup v)$
and
$(\hat{e}_1 \cup \hat{e}_2 \cup v)$ to $\tilde{S}$ which intersect at a point
$\tilde{v}$ which is a lift of $v$.  Lift $h$ to a map $\tilde{h}: \tilde{S}
\rightarrow \HHH$. Let $\gamma_1 = \tilde{h}(\alpha_1)$. 

Because $h$ is a $\pi_1$-injective simplicial hyperbolic surface, Lemma
\ref{extension} guarantees that the $SH$ extension $\bar{h}$
maps no two parabolic fixed points to the same point.  Then there exists a
geodesic axis $\gamma_2$ in the proper homotopy class of $\bar{h}(\alpha_2)$.

Let $\tilde{R}$ be the unique common perpendicular joining $\gamma_1$ to
$\gamma_2$.  Let $x = \tilde{R} \cap
\gamma_1$, and let $\tilde{P}$ be a geodesic arc in $\gamma_1$ joining
$\tilde{h}(\tilde{v})$ to $x$.  Let $R$ and $P$ be the projections to $N$ of
$\tilde{R}$ and $\tilde{P}$ respectively.  We construct $J$ by dragging $h$ along
$P$, then $R$.  The remainder of the proof is similar to that of Lemma
\ref{edge}.  \end{proof}

Finally, let $\Theta$ be a cofinite area torsion-free Fuchsian group, and let $S =
\HH / \Theta$.  Let $\Gamma$ be a Kleinian group such that there exists a weakly
type-preserving isomorphism between $\Theta$ and $\Gamma$.  Let $N = \HHH / \Gamma$.
A compact core $\hat{R}$ of $N$ will be homeomorphic to $\bar{S} \times [0,1]$ (Thm
10.5, Hempel \cite{hempel}) via a homeomorphism $\bar{h}$ where $\bar{S}$ is a compact
core of $S$.   Suppose $N$ has accidental parabolics
associated to both $\bar{h}(\bar{S},0)$ and $\bar{h}(\bar{S},1)$. 

If a triangulation on $S$ has edges which represent accidental parabolics associated
to both
$\bar{h}(\bar{S},0)$ and $\bar{h}(\bar{S},1)$, the following lemma enables us to
construct a continuous family which shrinks both the accidental parabolics.    

\begin{lemma}
\label{para to para} Let $h: S \rightarrow N$ be a $\pi_1$-injective, proper, and
weakly type-preserving map. Let
$v$ be the sole internal vertex of $T$ such that $\gamma^\prime = (e^\prime \cup v)$
is an accidental parabolic associated to one end, and let $\gamma^{\prime \prime} =
(e^{\prime
\prime} \cup v)$ be an accidental parabolic associated to the other end.   Then:
\begin{enumerate}
\item there exists a continuous family of
simplicial hyperbolic surfaces $H: S \times (-1,1) \rightarrow N$ with only one vertex
(which satisfies NLSC) such that:
\begin{enumerate}
\item $H_0$ is properly homotopic to $h$, 
\item as $t \rightarrow 1$, the length of $H_t(e^\prime
\cup v)$ decreases and converges to 0, and
\item as $t \rightarrow -1$, the length of
$H_t(e^{\prime \prime} \cup v)$ decreases and converges to 0;
\end{enumerate}  
\item $H(S, (-1,1)) \subset C(N)$;
\item $H(S,(-1,1))$ has compact closure in  $N - \bigcup P_i$, where \{$P_i$\}
is the set of rank one cusps of $N - \NOE$; 
\item for a component $S_i^\prime$ of $S - \gamma^\prime$, $H_t(S_i^\prime)$ lifts to
$C_{\delta^\prime(t)}(M_i)$, where $M_i = \HHH / h_*(\pi_1(S_i^\prime))$ and where 
$\delta^\prime: (-1,1) \rightarrow \Rr$ is a function such that as $t \rightarrow 1$,
$\delta^\prime(t)
\rightarrow 0$.  Furthermore, for $t \in (-1,1)$, the map $H_t: S_i^\prime \rightarrow
N$ has fewer
accidental parabolics than $H_0: S \rightarrow N$; and
\item for a component $S_i^{\prime \prime}$ of $S - \gamma^{\prime \prime}$,
$H_t(S_i^{\prime \prime})$ lifts to
$C_{\delta^{\prime \prime}(t)}(M_i)$, where $M_i = \HHH / h_*(\pi_1(S_i^{\prime
\prime}))$ and where 
$\delta^{\prime \prime}: (-1,1) \rightarrow \Rr$ is a function such that as $t
\rightarrow -1$,
$\delta^{\prime \prime}(t)
\rightarrow 0$.  Furthermore, for $t \in (-1,1)$, the map $H_t: S_i^{\prime \prime}
\rightarrow N$ has fewer
accidental parabolics than $H_0: S \rightarrow N$.
\end{enumerate}
\end{lemma}

\begin{proof}
The proof of this lemma mimics that of Lemma \ref{length to 0}.  We give a brief
outline here.  Let $W^\prime: S^1 \times [0,\infty) \rightarrow N$ be a proper
homotopy of $h(e^\prime \cup v)$ into a cusp of $N$, and let  $W^{\prime \prime}: S^1
\times (-\infty,0] \rightarrow N$ be a proper homotopy of $h(e^{\prime \prime} \cup
v)$ into a cusp of $N$.  Let $\alpha^\prime$ and $\alpha^{\prime \prime}$ be
components of the pre-images of $(e^\prime \cup v)$ and $(e^{\prime \prime} \cup
v)$ to $\tilde{S}$, respectively, that intersect at $\tilde{v}$,  a lift of $v$.
Lift the maps $h$,  $W^\prime$, and $W^{\prime \prime}$ to maps $\tilde{h}: \tilde{S}
\times \HHH$, $\tilde{W}^\prime : \Rr \times [0,\infty) \rightarrow \HHH$,
and $\tilde{W}^{\prime \prime}: \Rr \times (-\infty,0] \rightarrow \HHH$ such
that $\tilde{h}(\alpha^\prime) = \tilde{W}^\prime(\Rr,0)$ and
$\tilde{h}(\alpha^{\prime \prime}) = \tilde{W}^{\prime \prime}(\Rr,0)$, and such
that $\tilde{h}(v) = \tilde{W}^\prime(t,0) = \tilde{W}^{\prime \prime}(t,0)$ for
some $t \in \Rr$.  

Let $\tilde{R}$ be the geodesic line properly homotopic to
$[\tilde{W}^\prime(\{t\} \times [0,\infty))]
\cup [\tilde{W}^{\prime \prime}(\{t\} \times (\infty,0])]$.  Then the endpoints of
$\tilde{R}$ at the sphere at infinity are distinct parabolic fixed points.  Note that
$\tilde{v}$ may not lie on $\tilde{R}$.

Let $Y: \Rr \times [0,1] \rightarrow \HHH$ be a proper homotopy between
 $Y(\Rr \times  \{0\}) = [\tilde{W}^\prime(\{t\} \times [0,\infty))]
\cup [\tilde{W}^{\prime \prime}(\{t\} \times (\infty,0])]$ and $Y(\Rr \times \{1\}) =
\tilde{R}$. 
 Let $s \in \Rr$ be such that $Y(s,0) = \tilde{h}(\tilde{v})$.  We can properly
homotop $h$ to a map $f:(S,T) \rightarrow N$ which maps $v$ to $Y(s,1)$.  Using
Bonahon's construction, $f=F(h,T): (S,T) \rightarrow N$ is a $\pi_1$-injective
simplicial pre-hyperbolic surface.  We then drag $f$ along the projection $R$ of
$\tilde{R}$ to obtain a continuous family of simplicial pre-hyperbolic surfaces $H: S
\times (-1,1) \rightarrow N$.

Now we will show that $H$ is in fact a continuous family of simplicial hyperbolic
surfaces.  We will use Lemma \ref{tetrahedron} to verify that
the internal vertex $v$ satisfies $NLSC$ with respect to $H_t$ for $t \in (-1,1)$. 
Without loss of generality, we can renormalize so that $\tilde{R}$ has endpoints at
$\infty$ and $0$, and that the parabolic element with fixed point at
$\infty$ is $g_1$ and the other parabolic element with fixed point at $0$ is $g_2$.  

There are
four edges of $\tilde{H}_t(\tilde{T})$ which have one endpoint at
$\tilde{H}_t(\tilde{v})$ have other endpoints at $g_1(\tilde{H}_t(\tilde{v})),
g_{1}^{-1}(\tilde{H}_t(\tilde{v})), g_2(\tilde{H}_t(\tilde{v}))$, and
$g_{2}^{-1}(\tilde{H}_t(\tilde{v}))$.
By construction $\tilde{H}_t(\tilde{v})$ lies on $\tilde{R}$.  Note that 
$\tilde{H}_t(\tilde{v})$ lies on the same horocycle based at $\infty$ as
$g_1(\tilde{H}_t(\tilde{v}))$ and $g_1^{-1}(\tilde{H}_t(\tilde{v}))$. 
Therefore, the geodesic spanned by $g_1(\tilde{H}_t(\tilde{v}))$ and
$g_1^{-1}(\tilde{H}_t(\tilde{v}))$ intersects $\tilde{R}$ above
$\tilde{H}_t(\tilde{v})$.  Note that 
$\tilde{H}_t(\tilde{v})$ lies on the same horocycle based at $0$ as
$g_2(\tilde{H}_t(\tilde{v}))$ and $g_{2}^{-1}(\tilde{H}_t(\tilde{v}))$.
Therefore the geodesic spanned by $g_2(\tilde{H}_t(\tilde{v}))$ and
$g_2^{-1}(\tilde{H}_t(\tilde{v}))$ intersects $\tilde{R}$ below
$\tilde{H}_t(\tilde{v})$.  Hence, $\tilde{H}_t(\tilde{v})$ lies in the
tetrahedron spanned by the four other endpoints.  By Lemma \ref{tetrahedron},
the vertex $v$ satisfies {\it NLSC\/} with respect to $H_t$ for $t \in (-1,1)$.  Thus,
$H$ is a continuous family of simplicial hyperbolic surfaces. 

The remainder of the proof is similar to that of Lemma \ref{length to 0}.

\end{proof}

\section{Bounds for Surface Groups}

The following theorem of Kerckhoff-Thurston \cite{kerck/thur} states that given a
Kleinian group that is type-preserving isomorphic to a Fuchsian group, then there
exists an upper bound on the injectivity radius for points in the convex core of the
associated hyperbolic 3-manifold.

\begin{theorem}
\label{KT thm} (Kerckhoff-Thurston \cite{kerck/thur}) Let $\Theta$ be a cofinite area
torsion-free Fuchsian group, and let $S = \HH / \Theta$.  Then there exists a constant
$K_S$ such that for any Kleinian group $\Gamma$ such that there exists a
type-preserving isomorphism between $\Theta$ and $\Gamma$, and for $x \in C(N)$ where
$N = \HHH / \Gamma$,
$inj_N(x) \leq K_S$.
\end{theorem}

A proof of their theorem also appears in Canary \cite{canary1}.  A special case of the
above theorem is the version given in the introduction: \ Given an
$I$-bundle
$M$ over a closed surface $S$, then there
exists an upper bound on injectivity radius for points in the convex core of any
hyperbolic 3-manifold
$N$ without cusps such that $N$ is homeomorphic to the interior of $M$.  The
special case excludes all parabolic elements.

Our goal in this section is to prove an extension of this theorem that includes the
possibility of accidental parabolics. 

\begin{theorem}
\label{surface bound} Let $\Theta$ be a cofinite area torsion-free Fuchsian group, and
let $S = \HH / \Theta$.  There exists a constant $L_S$ such that for any  Kleinian
group $\Gamma$ such that there exists a weakly type-preserving isomorphism between
$\Theta$ and $\Gamma$, and for $x \in C(N)$ where
$N =
\HHH / \Gamma$, $inj_N(x) \leq L_S$.
\end{theorem}

\begin{proof}  We first give an outline of the proof.  We will use induction
on the number of accidental parabolics in the weakly type-preserving isomorphism. 
This will divide the proof into three cases, each of which will be proved by dividing
the convex core into a compact core and its complement.  Using the induction hypothesis
and continuous families of simplicial hyperbolic surfaces, points in both the compact
core of the convex core and the complement will have bounded injectivity radius.
This completes the outline of the proof.

Recall that Sullivan \cite{sullivan} has shown that for a finitely generated Kleinian
group, there can be at most a finite number of conjugacy classes of maximal abelian
subgroups associated to accidental parabolics.  We will
use induction on the number of accidental parabolics. 


The base step in the induction is the case where there are no accidental parabolics. 
In this case, $\Theta$ is isomorphic to
$\Gamma$ via a type-preserving isomorphism.  Then by Theorem \ref{KT thm}, we know
that there exists a constant
$K_S$ such that the injectivity radius for points in the convex core is bounded above
by
$K_S$.  Thus, we are done.

Now let us deal with the induction step.  Suppose that the induction hypothesis is
true for a weakly type-preserving isomorphism with fewer than $p$ accidental
parabolics.  Let $\Gamma_j$ be a Kleinian group that is weakly type-preserving
isomorphic to a Fuchsian group $\Theta_j$ with fewer than $p$ accidental parabolics.
Let $M_j = \HHH / \Gamma_j$ and $S_j = \HHH / \Theta_j$.  Then by the induction
hypothesis, for $x \in C(M_j)$, $inj_{M_j}(x) \leq L_{S_j}$.  Formally, we are
assuming the existence of an upper bound $L_{S_j,p^\prime}$ on injectivity radius for
points in $C(M_j)$ where the bound depends on the number
$p^\prime < p$ of accidental parabolics in the isomorphism between $\Theta_j$ and
$\Gamma_j$, and we will use this bound to establish an upper bound $L_{S,p}$ on
injectivity radius for points in $C(N)$ where $\Theta$ and $\Gamma$ are weakly
type-preserving isomorphic with $p$ accidental parabolics.  Here,
$L_{S_j,0} = K_{S_j}$, the bound obtained from Theorem \ref{KT thm}.  For
ease of exposition, we will let $L_{S_j} = max_{p^\prime < p}
\{L_{S_j,p^\prime}\}$ denote the upper bound on injectivity radius for points in
$C(M_j)$ for all weakly type-preserving isomorphisms between $\Theta_j$ and $\Gamma_j$
with less than
$p$ accidental parabolics, and we will use this bound to establish an upper bound
$L_S$ on injectivity radius for points in $C(N)$, where $\Theta$ and $\Gamma$ are 
weakly type-preserving isomorphic with $p$ accidental
parabolics.

Let $0 < \epsilon < \epsilon_3$.  Let $\hat{R}$ be a relative compact core of $\NOE$
with associated parabolic locus $\hat{P}$, and let $\bar{S}$ be a compact core of
$S$.  Then because $\Gamma$ is isomorphic to
$\Theta$, $\hat{R}$ is homeomorphic to $\bar{S} \times [0,1]$. (Thm 10.5,
Hempel
\cite{hempel})  Let $\bar{h}: \bar{S}
\times [0,1] \rightarrow \hat{R}$ be a homeomorphism where $\bar{h}$ is in the
homotopy class induced by the weakly type-preserving isomorphism and
$\bar{h}(\partial \bar{S} \times [0,1])$ is a collection of components of $\hat{P}$. 
We can extend the homeomorphism $\bar{h}|_{(\bar{S},0)}$ to a homeomorphism $h:S
\rightarrow N$ on
$S$ as follows: \  Let $h$ be a proper embedding such that
$h(x) = \bar{h}(x,0)$ for $x \in \bar{S}$, and let $h(S) =
(\bar{h}(\bar{S},0))_P$ where
$(\bar{h}(\bar{S},0))_P$ is the parabolic extension of $\bar{h}(\bar{S},0)$.

Recall that every parabolic element of $\Gamma$ can be associated to a closed curve in
the parabolic locus
$\hat{P}$ of $\hat{R}$.  In particular, every accidental parabolic can be associated to
a closed curve in a component of $[\hat{P} \cap
\bar{h}(\bar{S},0)] 
\cup [\hat{P} \cap \bar{h}(\bar{S},1)]$.  Let the accidental parabolics that are
associated to a curve on
$[\hat{P} \cap \bar{h}(\bar{S},0)]$ be called {\it accidental parabolics  associated 
to 
$A$\/}; let the others be called {\it accidental parabolics  associated to
$B$\/}.  Let us first assume that there are accidental parabolics associated to both
$A$ and $B$; we will discuss the proof in the case where accidental parabolics are
associated to only one of $A$ and $B$ later.  \\

Case I.  Accidental Parabolics Associated to $A$ and $B$.  \\
 
Let us give an outline of the proof in this case.  We will divide the convex core into
a compact core and its complement.  Then we will cover the convex core with continuous
families of simplicial hyperbolic surfaces and projections of convex cores of covers
with fewer accidental parabolics (which, by the induction hypothesis, have bounded
injectivity radius), so that points in the compact core will have bounded injectivity
radius.  The remaining portions of the convex core will either lie in an
$\epsilon$-thin part, or will lie in the projection of the convex core of a cover with
fewer accidental parabolics. This completes the outline of the proof in this case.

The first step of the proof in Case I will be to find a continuous family of simplicial
hyperbolic surfaces which shrinks an accidental parabolic associated to $A$ and an
accidental parabolic associated to $B$.  We will want points in the image of this
continuous family to lie in the convex core and to have bounded injectivity radius.

Consider a triangulation $T_1$ containing an edge which represents an accidental
parabolic associated to $A$.  Either there exists an edge in $T_1$ which represents
an accidental parabolic associated to $B$, or not.  In the latter case, we can perform
the modified $HH$ sequence of elementary moves altering $T_1$ to a triangulation $T_2$
containing an edge
 representing an accidental parabolic associated to $B$.  Let $T_A$ be the ``last''
triangulation in the sequence in which there is an edge $\gamma_A$ representing an
accidental parabolic associated to $A$.  Let
$T_B$ be the ``first'' triangulation in the sequence in which there is an edge
$\gamma_B$ representing an accidental parabolic associated to $B$. There are two
possibilities: (1) $T_A$ and $T_B$ are adjacent, or (2) they are separated by
triangulations whose edges are not accidental parabolics associated to either $A$ or
$B$.  

In either of these cases, we say that the sequence of triangulations has
$Property \ T$, that is, there exist triangulations $T_A$ and $T_B$ in our sequence such
that: 
\begin{enumerate}
\item the edges of $T_A$ may only represent accidental parabolics associated to $A$,
\item the edges of $T_B$ may only represent accidental parabolics associated to $B$, and
\item the edges of each
intermediary triangulation do not represent accidental parabolics associated to either
$A$ or $B$.
\end{enumerate}
Note that when $T_A$ and $T_B$ are adjacent, there are no intermediary triangulations.

Thus, we can see that either a sequence of triangulations has $Property \ T$, or
there exists a triangulation with edges which represent accidental parabolics associated
to both $A$ and $B$.
\\ 

Case IA.  There Exists a Sequence of Triangulations with $Property  \ T$. \\

Suppose our sequence of triangulations has $Property \ T$.   
Now there are two cases as well:  either $T_A$ and $T_B$ are separated by at least one
intermediary triangulation, or not. \\ 

Case IA(1).   $T_A$ and $T_B$ Separated by at Least one Intermediary Triangulation.
\\

In this case, there exist triangulations
$T_A$ and $T_B$ such that exactly one edge $\gamma_A$ of $T_A$ is an accidental
parabolic associated to $A$, and exactly one edge $\gamma_B$ of $T_B$ is an accidental
parabolic associated to $B$, such that the edges of each intermediary triangulation are
not accidental parabolics associated to either
$A$ or $B$. The fact that there is exactly one edge in each
triangulation associated to an accidental parabolic follows directly from the fact that
an elementary move alters exactly one edge in the triangulation.

If $S$ is a 4-punctured sphere, then there is exactly one nonperipheral edge in its
suitable, one-vertex triangulation.  Thus, in $T_A$, this nonperipheral edge is
$\gamma_A$; similarly, in $T_B$, this nonperipheral edge is $\gamma_B$.  Let
$\hat{e}_1$ and $\hat{e}_2$ be any two edges in $T_A$ joining distinct ideal vertices
to the sole internal vertex.  We can use Lemma
\ref{makepractical} to construct a $\pi_1$-injective practical simplicial hyperbolic
surface
$f:(S,T_A)
\rightarrow N$ with marked edges $\hat{e}_1$ and $\hat{e}_2$ such that $f$ is properly
homotopic to $h:S \rightarrow N$. Let
$\hat{e}_3$ and $\hat{e}_4$ be any two edges in $T_B$ joining distinct ideal vertices
to the sole internal vertex.    Similarly, construct a $\pi_1$-injective practical
simplicial hyperbolic surface $g:(S,T_B)
\rightarrow N$ with marked  edges $\hat{e}_3$ and $\hat{e}_4$ such that $g$ is properly
homotopic to $h:S \rightarrow N$.  

By repeatedly using Lemmas \ref{3 elem move} and \ref{3 edge}, we can
 construct the $HH$ interpolation 
$H_1: S
\times [0,1] \rightarrow N$ such that
$H_1|_{(S,0)}:S \rightarrow N$ is the map $f:(S,T_A) \rightarrow N$ and 
$H_1|_{(S,1)}:S \rightarrow N$ is the map $g:(S,T_B) \rightarrow N$. (Note that we have
henceforth abandoned the notation $H_t = H|_{(S,t)}$ from Section 4.)  Here, we are
using Lemma \ref{3 elem move} to perform the elementary moves in the modified $HH$
sequence of triangulations which do not involve the existing marked edges, and we are
using Lemma \ref{3 edge} to change the marked edges when necessary.

Now we can use Lemma \ref{3 length to 0}, setting $H_1|_{(S,1)} = H_2|_{(S,0)}$, to
obtain a continuous family of simplicial hyperbolic surfaces $H_2: S \times [1,2)
\rightarrow N$ that shrinks the accidental parabolic $\gamma_A$.  Similarly, we can
apply Lemma
\ref{3 length to 0}, setting $H_1|_{(S,0)} = H_3|_{(S,0)}$, to obtain
a continuous family of simplicial hyperbolic surfaces $H_3: S \times (-1,0] \rightarrow
N$ that shrinks the accidental parabolic $\gamma_B$.

If $S$ is not a 4-punctured sphere, then  there exists a
nonperipheral edge $\hat{e}_1$ in $T_A$ which does not represent an accidental
parabolic, and hence can be mapped to a closed geodesic.  Similarly,  there
exists a nonperipheral edge $\hat{e}_2$ in $T_B$ which does not represent an accidental
parabolic, and hence can be mapped to a closed geodesic.
Using Lemma \ref{makeuseful}, there exists a $\pi_1$-injective useful simplicial
hyperbolic surface $f:(S,T_A) \rightarrow N$ with distinguished edge $\hat{e}_1$ such
that $f$ is properly homotopic to $h:S \rightarrow N$.    Similarly, there exists a
$\pi_1$-injective useful simplicial hyperbolic surface $g:(S,T_B) \rightarrow N$ with
distinguished edge $\hat{e}_2$  such that $g$ is
properly homotopic to $h:S \rightarrow N$.  

By repeatedly using Lemma \ref{elem move}, we can
 construct a continuous family of simplicial hyperbolic surfaces that performs the
elementary moves indicated in the modified $HH$ sequence of
triangulations that interpolates between $T_A$ and $T_B$.  We can also use Lemma
\ref{edge} to construct a continuous family that alters the distinguished edge from
$\hat{e}_1$ to $\hat{e}_2$ as needed in the modified $HH$ sequence.  By concatenating
all of these continuous families, we construct the
$HH$ interpolation 
$H_1: S
\times [0,1] \rightarrow N$ such that
$H_1|_{(S,0)}:S \rightarrow N$ is the map $f:(S,T_A) \rightarrow N$ and 
$H_1|_{(S,1)}:S \rightarrow N$ is the map $g:(S,T_B) \rightarrow N$. 

Now we can use Lemma \ref{length to 0}, setting $H_1|_{(S,1)} = H_2|_{(S,0)}$, to obtain
a continuous family of simplicial hyperbolic surfaces $H_2: S \times [1,2) \rightarrow
N$ that shrinks the accidental parabolic $\gamma_A$.  Similarly, we can apply Lemma
\ref{length to 0}, setting $H_1|_{(S,0)} = H_3|_{(S,0)}$, to obtain
a continuous family of simplicial hyperbolic surfaces $H_3: S \times (-1,0] \rightarrow
N$ that shrinks the accidental parabolic $\gamma_B$. \\

Case IA(2).  $T_A$ and $T_B$ Not Separated by at Least one Intermediary
Triangulation. 
\\

In this case, we can conclude that $T_A$ and $T_B$ are adjacent triangulations, that is,
the elementary move alters the edge $\gamma_A$ to the edge $\gamma_B$. 

If  $S$ is a 4-punctured sphere, then because its suitable, one-vertex triangulation
contains exactly one nonperipheral edge,  we know that in $T_A$,
$\gamma_A$ is that nonperipheral edge, and that in $T_B$, $\gamma_B$ is that
nonperipheral edge.   Let
$\hat{e}_1$ and $\hat{e}_2$ be any two edges in $T_A$ joining distinct ideal vertices
to the sole internal vertex.  We can use Lemma
\ref{makepractical} to construct a $\pi_1$-injective practical simplicial hyperbolic
surface
$f:(S,T_A)
\rightarrow N$ with marked edges $\hat{e}_1$ and $\hat{e}_2$ such that $f$ is properly
homotopic to $h:S \rightarrow N$.   Let
$\hat{e}_3$ and $\hat{e}_4$ be any two edges in $T_A$ joining distinct ideal vertices
to the sole internal vertex.   Similarly, construct a $\pi_1$-injective practical
simplicial hyperbolic surface $g:(S,T_B)
\rightarrow N$ with marked  edges $\hat{e}_3$ and $\hat{e}_4$ such that $g$ is properly
homotopic to $h:S \rightarrow N$.  

We can use Lemma \ref{3 elem move} to perform the elementary move that separates the
triangulations $T_A$ and $T_B$.  We can also use Lemma \ref{3 edge} to construct a
continuous family that alters the marked edges from $\hat{e}_1$ and $\hat{e}_2$ to
$\hat{e}_3$ and $\hat{e}_4$.  After concatenating all of these continuous families, we
construct the
$HH$ interpolation 
$H_1: S
\times [0,1] \rightarrow N$ such that
$H_1|_{(S,0)}:S \rightarrow N$ is the map $f:(S,T_A) \rightarrow N$ and 
$H_1|_{(S,1)}:S \rightarrow N$ is the map $g:(S,T_B) \rightarrow N$.

Now  we can
apply Lemma
\ref{3 length to 0}, setting $H_1|_{(S,0)} = H_3|_{(S,0)}$, to obtain
a continuous family of simplicial hyperbolic surfaces $H_3: S \times (-1,0] \rightarrow
N$ that shrinks the accidental parabolic $\gamma_B$.

If $S$ is not a 4-punctured sphere, then there exist edges
$\hat{e}_1 \neq \gamma_A$ of
$T_A$  and $\hat{e}_2 \neq \gamma_B$ of
$T_B$, each of which can be mapped to a closed geodesic.  Therefore, we can use Lemma
\ref{makeuseful} to construct a $\pi_1$-injective useful simplicial hyperbolic surface
$f:(S,T_A)
\rightarrow N$ with distinguished edge $\hat{e}_1$ 
such that $f$ is properly homotopic to $h:S \rightarrow N$.    Similarly, construct
a $\pi_1$-injective useful simplicial hyperbolic surface $g:(S,T_B) \rightarrow N$ with
distinguished edge $\hat{e}_2$ such that $g$ is
properly homotopic to $h:S \rightarrow N$.  

We can use the Lemma \ref{elem move} to perform the elementary move that separates the
triangulations $T_A$ and $T_B$.  We can also use Lemma \ref{edge} to construct a
continuous family that alters the distinguished edge from $\hat{e}_1$ to
$\hat{e}_2$.  After concatentating all of these continuous families, we 
construct the
$HH$ interpolation 
$H_1: S
\times [0,1] \rightarrow N$ such that
$H_1|_{(S,0)}:S \rightarrow N$ is the map $f:(S,T_A) \rightarrow N$ and 
$H_1|_{(S,1)}:S \rightarrow N$ is the map $g:(S,T_B) \rightarrow N$.

Now we can use Lemma \ref{length to 0}, setting $H_1|_{(S,1)} = H_2|_{(S,0)}$, to obtain
a continuous family of simplicial hyperbolic surfaces $H_2: S \times [1,2) \rightarrow
N$ that shrinks the accidental parabolic $\gamma_A$.  Similarly, we can apply Lemma
\ref{length to 0}, setting $H_1|_{(S,0)} = H_3|_{(S,0)}$, to obtain
a continuous family of simplicial hyperbolic surfaces $H_3: S \times (-1,0] \rightarrow
N$ that shrinks the accidental parabolic $\gamma_B$.  \\

Case IB.  There Does Not Exist a Sequence of Triangulations with $Property  \ T$.
\\

If we cannot construct a sequence of triangulations with $Property \ T$, then as we
argued before, there exists a triangulation that contains edges that are accidental
parabolics associated to both
$A$ and $B$.  Using Lemma
\ref{para to para} we can construct a continuous family of simplicial hyperbolic
surfaces $H_1: S
\times (-1,2) \rightarrow N$ that shrinks the accidental parabolic $\gamma_A$ as $t
\rightarrow 2$, and that  shrinks the accidental parabolic $\gamma_B$ as $t \rightarrow
-1$. \\

Let $W = H_1(S,[0,1]) \cup H_2(S,[1,2)) \cup H_3(S,(-1,0])$ (or
$W = H_1[S,(-1,2)]$ in Case IB) be the union of the images of the
continuous families of simplicial hyperbolic surfaces.   
Lemmas \ref{edge}, \ref{elem move}, \ref{length to 0}, \ref{3 edge}, \ref{3 elem move},
\ref{3 length to 0}, and \ref{para to para}
 guarantee that $W \subset C(N)$. 
 Furthermore, we
also know that for any $x \in W$ lies in the image of a
simplicial hyperbolic surface, and hence, by Lemma
\ref{shsinj}, $inj_N(x) \leq K_S$ where the constant $K_S$ depends on the Euler
characteristic of $S$. Now we have bounded the injectivity radius in a portion of the
convex core. 


Because $\Gamma$ is isomorphic to a surface group which contains no ${\mathbb Z}
\oplus {\mathbb Z}$ subgroups,
$N - \NOE$ possesses only rank one cusps.  There exists a retraction map, $\phi: N
\rightarrow \NOE$, such that $\phi$ is the identity on $\NOE$, and $\phi$ retracts each
component of $N - \NOE$ onto $\partial \NOE$.  More explicitly, $\phi$ can be
defined as follows: \ Let $P_i$ denote a component of $N -
\NOE$. Each of the $P_i$ can be parametrized as $S^1 \times (-\infty, \infty) \times
[0,\infty)$ using the coordinates $(\theta_i, r_i, s_i)$.  Then $\phi|_{P_i}: P_i
\rightarrow P_i$ is the map $\phi(\theta_i, r_i, s_i) = (\theta_i, r_i, 0)$. Note that
$\phi(W) \subset C(N)$, because $C(N)$ is convex. 

Now we will show that $\overline{\phi(W)}$ is compact. 
Canary (Lem 4.4
\cite{canary1}) has shown that for a simplicial hyperbolic surface $f$ whose vertices
are $NLSC$, 
$$\max_{x \in f^{-1}(P_i)} r_i(x)  =  \max_{x \in f^{-1}(\partial P_i)} r_i(x)$$ where
$r_i(x)$ denotes the $r_i$-coordinate of $f(x)$.  Therefore, to show that
$\overline{\phi(W)}$ is compact, it suffices to show that $W$ has compact closure in 
$N - \bigcup P_i$.  Because $H_1(S \times \{t\})$ has compact closure in $N - \bigcup
P_i$ for
$t
\in [0,1]$ and $[0,1]$ is compact, we know that the image $H_1(S,[0,1])$ in $N - \bigcup
P_i$ is compact.  Thus, we need only consider the images
$H_2(S,[1,2))$, $H_3(S,(-1,0])$, and $H_1(S,(-1,2))$. By Lemmas \ref{length
to 0}, \ref{3 length to 0}, and \ref{para to para}, we know that $H_2(S,[1,2))$,
$H_3(S,(-1,0])$, and $H_1(S,(-1,2))$ have compact closure in $N -
\bigcup P_i$, and hence $\overline{\phi(W)}$ is compact.

The next step in the proof of Case I involves finding a compact
core $R$ of $N - P_S$, where $P_S$ represents the collection of rank one cusps of
$\NOE$ associated to $\partial S$, such that points in $R$ have bounded injectivity
radius. 

Let $\xi: S \rightarrow N$
denote the map $H_1|_{(S,1)}: S \rightarrow N$.  Let \{$S_{j}$\} be the components of
$S - \gamma_A$, and let \{$S_{k}$\} be the components of
$S - \gamma_B$.  Let \{$M_j=\HHH / \xi_*(\pi_1(S_{j}))$\} and \{$M_k=\HHH /
\xi_*(\pi_1(S_{k}))$\} be covers of $N$ with covering maps \{$p_j$\} and \{$p_k$\}
respectively.  Let
$P_A$ and $P_B$ be the rank one cusps of $\NOE$ associated to $\gamma_A$ and
$\gamma_B$ respectively.  Let
${\mathcal P} = P_A \cup P_B \cup P_S$. 

Given $0 < \epsilon < \epsilon_3$.  By  Lemma \ref{compactincore}, there exists a
relative compact core
$R^*$ of
$\overline{C_{\epsilon}(N)} - {\mathcal P}$ such that $\partial C_\epsilon(N) \cap
C^\circ_\epsilon(N)
\subset
\partial R^*$ and such that $\overline{\phi(W)} \subset R^*$.  
Let
\{$U_j^\prime$\} be the components of $[N - {\mathcal P}] - int \ R^*$
such that each $U_j^\prime$ is adjacent to
$P_A$, and let \{$U_k^\prime$\} be the components of $[N - {\mathcal P}] -
int \ R^*$ such that each $U_k^\prime$ is adjacent to $P_B$.
Let
\{$U_j$\} be the components of $[C_{2 \epsilon}(N) - {\mathcal P}] - int \ R^*$ such
that each $U_j$ is adjacent to
$P_A$, and let \{$U_k$\} be the components of $[C_{2 \epsilon}(N) - {\mathcal P}] -
int \ R^*$ such that each $U_k$ is adjacent to $P_B$.
Then by construction, $U_j \subset U_j^\prime$, and $U_k \subset U_k^\prime$.

In the following lemma,  we will choose a compact core
$R$ of $N - P_S$ such that $R$ contains $R^*$,
$R$ lies in $C_{2 \epsilon}(N)$, and points in $R$ have bounded injectivity radius. 
Figure ~\ref{fig:surfacecase} may help to keep track of all of the different surfaces
used in the proof.

\begin{figure}[ht]
\begin{center}
\BoxedEPSF{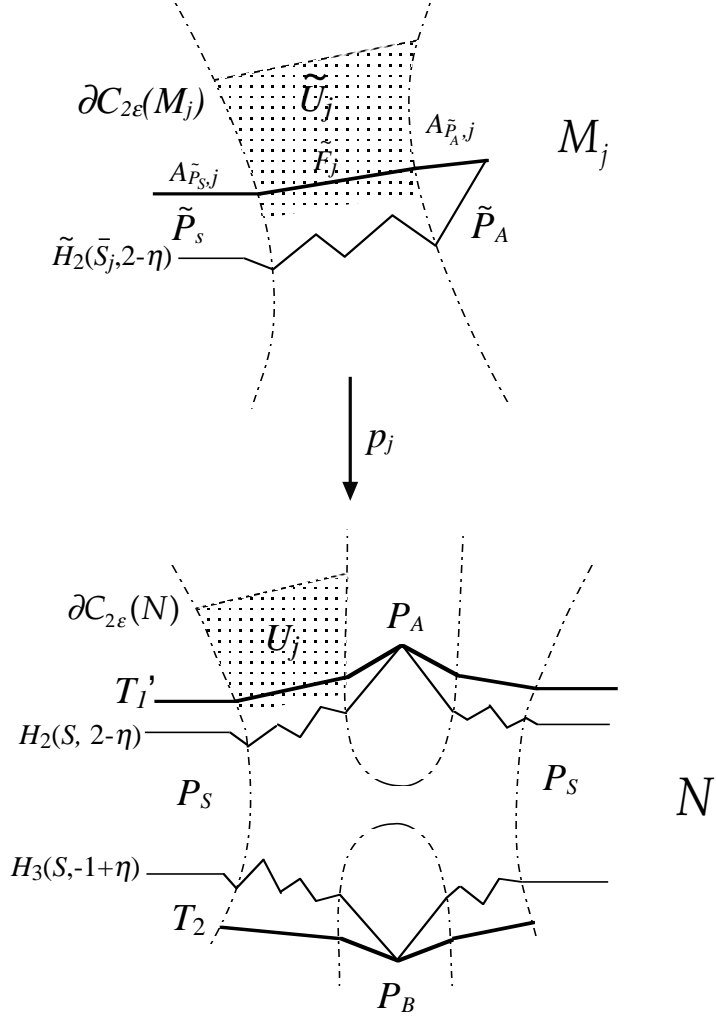 scaled 600}
\caption{Surfaces Used in Lemma ~\ref{surface core}} \label{fig:surfacecase}
\end{center}
\end{figure}

\begin{lemma}
\label{surface core} For all $0 < \epsilon < \epsilon_3$, there exists a compact core
$R$ of $N - P_S$ bounded by surfaces
$T_1$ and $T_2$ such that the following are true:
\begin{enumerate}
\item $R^* \subset R \subset C_{2 \epsilon}(N)$, 

\item if $x \in R$, then $inj_N(x) \leq max \{K_S, L_{S_j} + 2 \epsilon, L_{S_k} + 2
\epsilon\}$

\item for each $j$, $\partial (T_1 \cap U_j^\prime) \subset \partial {\mathcal P}$, and
for each
$k$,  $\partial (T_2 \cap U_k^\prime) \subset \partial {\mathcal P}$,

\item for each $j$, there exist $\tilde{U}_j^\prime \subset M_j$ such that
$p_j|_{\tilde{U}_j^\prime}: \tilde{U}_j^\prime \rightarrow U_j^\prime$ is a
homeomorphism, and
$(p_j|_{\tilde{U}_j^\prime})^{-1}(T_1 \cap U_j^\prime) \subset \overline{C_{2
\epsilon}(M_j)}$, and 

\item for each $k$, there exist $\tilde{U}_k^\prime \subset M_k$ such that
$p_k|_{\tilde{U}_k^\prime}: \tilde{U}_k^\prime \rightarrow U_k^\prime$ is a
homeomorphism, and
$(p_k|_{\tilde{U}_k^\prime})^{-1}(T_2 \cap U_k^\prime) \subset \overline{C_{2
\epsilon}(M_k)}$.  
\end{enumerate}
\end{lemma}

\begin{proof} Let us begin with an outline of the proof.  We will first choose ${\mathcal
T}_1$ and ${\mathcal T}_2$ to be incompressible surfaces in $C(N)$ such that subsurfaces
of
${\mathcal T}_1$ and ${\mathcal T}_2$ lift to
$C_\epsilon(M_j)$ and
$C_\epsilon(M_k)$, and ${\mathcal T}_1$ and ${\mathcal T}_2$ bound a submanifold
 of
$N$.  For each $i$, we will construct an embedded surface $T_i^\prime \subset C_{2
\epsilon}(N)$ that is properly homotopic to ${\mathcal T}_i$ via a homotopy
$\zeta_i$ whose image lies in $\bigcup p_j(C_{2 \epsilon}(M_j))$ such that the
\{$T_i^\prime$\} bound a submanifold
that contains
${\mathcal T}_1$ and ${\mathcal T}_2$.  For each $i$, let $T_i = T_i^\prime \cap (N -
P_S)$.  The
\{$T_i$\} bound a compact core $R$ of $N - P_S$  such that $R$ contains $R^*$.  We
will show that $R$ is contained in the union of
$W$ and the images of the homotopies \{$\zeta_i$\}, so that points in
$R$ will have bounded injectivity radius.

Because $\gamma_A$ is not an accidental parabolic
in $M_j$, $M_j$ has fewer accidental parabolics than
$N$. Therefore, by the induction hypothesis, for every $x \in C(M_j)$,
$inj_{M_j}(x) \leq L_{S_j}$ where $L_{S_j}$ depends on the surface
$S_j$.   Similarly, for every $x \in C(M_k)$, $inj_{M_k}(x) \leq L_{S_k}$, where
$L_{S_k}$ depends on $S_k$.

Let $\bar{S}_j$ be a subsurface of $S$ obtained by cutting $S$ along $\gamma_A$.  Then
$\bar{S}_j$ is a surface with boundary containing of one or two copies of $\gamma_A$.
Let $\bar{S}_k$ be the analogous subsurface of $S$, obtained by cutting
$S$ along $\gamma_B$, with boundary consisting
of one or two copies of $\gamma_B$.  By Lemma \ref{length to 0} (or Lemma \ref{3 length
to 0} if $S$ is a 4-punctured sphere), there exists
$\eta > 0$ such that for every $j$ and $k$,
$H_2(\bar{S}_j,2-\eta)$ lifts to $C_\epsilon(M_j)$ and $H_3(\bar{S}_k,-1+\eta)$ lifts
to
$C_\epsilon(M_k)$. We also require that $H_2(v,2-\eta) \subset P_A$ and $H_3(v,-1+\eta)
\subset P_B$. 

Let  ${\mathcal T}_1 = H_2(S,2-\eta) \subset C(N)$, and let ${\mathcal T}_2 = H_3(S,-1+\eta)
\subset C(N)$. 

Now we will construct the surfaces $T_1^\prime$ and $T_2^\prime$.

Let
$P^*$ be the parabolic locus associated to $R^*$.  Because the fundamental group of
$N$ is isomorphic to a surface group, we know that $R^*$ is homeomorphic
to $\bar{S} \times [0,1]$. (Thm 10.5, Hempel \cite{hempel})  Because $P^*$ is a disjoint
collection of incompressible annuli in $\partial R^*$ such that $\partial \bar{S}
\times I \subset P^*$, we can conclude that the components of
$\partial R^* - P^*$  are incompressible in $N$ and hence also in
$C_{2 \epsilon}(N)$. 

Note that $R^*$ is also a relative compact core for $N - {\mathcal
P}$ where the components of $\partial R^* - P^*$ are incompressible in $N$.  Then we can apply Lemma
\ref{prodstructure} to show that each
$U_j^\prime$ possesses a product structure.  By the Lifting Theorem, the inclusion map
of the product structure into $N$ lifts to an inclusion map  of the product
structure into
$M_j$.    Let
$\tilde{U}_j^\prime$ be the image of the lifted inclusion map in $M_j$.  Then
$\tilde{U}_j^\prime$ has a product structure, and
$p_j|_{\tilde{U}_j^\prime}:
\tilde{U}_j^\prime \rightarrow U_j^\prime$ is a homeomorphism. 

Note that $R^*$ is also a relative compact core for ${C_{2 \epsilon}(N)} - {\mathcal
P}$ where the components of $\partial R^* - P^*$ are incompressible in $C_{2
\epsilon}(N)$.   Then we can apply the argument of Lemma
\ref{prodstructure} to $C_{2 \epsilon}(N) - {\mathcal P}$ to show that each
$U_j$ possesses a product structure. By the Lifting Theorem, the inclusion map
of the product structure into $N$ lifts to an inclusion map  of the product
structure into
$M_j$.    Let
$\tilde{U}_j$ be the image of the lifted inclusion map in $M_j$.  Then
$\tilde{U}_j$ has a product structure, and
$p_j|_{\tilde{U}_j}:
\tilde{U}_j \rightarrow U_j$ is a homeomorphism. 

Let $\tilde{P}_A$ be the rank one cusp of $(M_j)^\circ_\epsilon$ associated to
$\gamma_A$ where $p_j(\tilde{P}_A) = P_A$.  Similarly, let $\tilde{P}_S$ be the rank
one cusp(s) of $(M_j)^\circ_\epsilon$ associated to the ideal vertices of $\bar{S}_j$
where $p_j(\tilde{P}_S) = P_S$.  We will let
$\tilde{\mathcal P}_j = \tilde{P}_A \cup \tilde{P}_S$.  

Now we will show that for each $j$, $\tilde{U}_j \subset C_{2 \epsilon}(M_j)$. 
Suppose not.  Then there exists $\tilde{x} \in \tilde{U}_j - C_{2 \epsilon}(M_j)$ and
a component
$\tilde{A}_j$ of $\partial C_{2 \epsilon}(M_j)$ such that $\tilde{A}_j$ separates
$\tilde{x}$ from $C_{2 \epsilon}(M_j)$.  Consider a geodesic ray
$\bar{g}_{\tilde{x}}$ in $M_j - C_{2 \epsilon}(M_j)$ that is perpendicular to
$\tilde{A}_j$ and passes through $\tilde{x}$.  Let
$g_{\tilde{x}}$ be the portion of $\bar{g}_{\tilde{x}}$ beginning at $\tilde{x}$. 
Note that by Lemma
\ref{expinj}, the injectivity radius strictly increases out a geometrically finite
end, so the ray
$g_{\tilde{x}}$ is entirely contained in $(M_j)^\circ_\epsilon$. Let $Q_j = U_j \cap
R^*$, and let $\tilde{Q}_j$ be the lift of $Q_j$ to
$\tilde{U}_j$.  Without loss of generality, let us assume that
$g_{\tilde{x}}$ intersects
$\tilde{Q}_j$, transversely. 
Then one of the three following cases occurs:
\
(1) $g_{\tilde{x}}$ intersects
$\tilde{Q}_j$ an odd number of times,
(2) $g_{\tilde{x}}$ intersects a different boundary component of the closure of
$\tilde{U}_j$, or (3) all but a compact portion of 
$g_{\tilde{x}}$ is completely contained in $\tilde{U}_j$.

Suppose $g_{\tilde{x}}$ intersects $\tilde{Q}_j$ an odd number of
times.  Because
$\tilde{Q}_j$ is a properly embedded surface separating
$M_j - \tilde{\mathcal P}_j$, $\tilde{Q}_j$ has two sides.  Let the component of
$[M_j - \tilde{\mathcal P}_j] - \tilde{Q}_j$ containing
$\tilde{U}_j$ be the positive side of
$\tilde{Q}_j$.  Then the last time $g_{\tilde{x}}$ intersects $\tilde{Q}_j$,
$g_{\tilde{x}}$ passes from the positive to the negative side of 
$\tilde{Q}_j$.  

Let $Q_j = p_j(\tilde{Q}_j) \subset U_j$.  Downstairs, we can also let the side of
$[C_{2 \epsilon}(N) - {\mathcal P}] - Q_j$ containing
$U_j$ be the positive side of $Q_j$.  Recall that because
$\phi[H_2(\bar{S}_j,2-\eta)] \subset \phi(W) \subset R^*$, and $U_j$ is a component of
$[C_{2
\epsilon}(N) - {\mathcal P}] - int \ R^*$, then we know
that $\phi[H_2(\bar{S}_j,2-\eta)]$ lies on the negative side of
$Q_j$. 

Because $p_j|_{\tilde{U}_j}$ is a homeomorphism and $Q_j$ lifts to
$\tilde{Q}_j$, we can conclude that $\phi[\tilde{H}_2(\bar{S}_j,2-\eta)]$ lies to the
negative side of $\tilde{Q}_j$.  Note that
$\phi[\tilde{H}_2(\bar{S}_j,2-\eta)]$ is an incompressible separating surface of
$M_j - \tilde{\mathcal P}_j$, and that
$\phi[\tilde{H}_2(\bar{S}_j,2-\eta)]$ is homologous to $\tilde{Q}_j$.  Also note that
$g_{\tilde{x}}$ leaves every compact set of
$(M_j)^\circ_\epsilon$ and all but a compact portion of
$g_{\tilde{x}}$ lies in $(M_j)^\circ_\epsilon - \tilde{U}_j$.  Therefore, if
$g_{\tilde{x}}$ intersects $\tilde{Q}_j$, then $g_{\tilde{x}}$ must also intersect
$\phi[\tilde{H}_2(\bar{S}_j,2-\eta)]$.  But $g_{\tilde{x}} \subset M_j - C_{2
\epsilon}(M_j)$ and
$\phi[\tilde{H}_2(\bar{S}_j,2-\eta)] \subset C_{2
\epsilon}(M_j) \cup N_{thin(\epsilon)}$.  So this is a contradiction.  

Suppose $g_{\tilde{x}}$ intersects a different boundary component of the closure of
$\tilde{U}_j$.  Then $p_j(g_{\tilde{x}})$ intersects some component $B_j$ of $\partial
C_{2 \epsilon}(N)$ such that $B_j$ is a boundary component of the closure of $U_j$ in
$N$.  By Lemma \ref{GFcover},
$B_j$ lifts to a component $\tilde{B}_j$ of $\partial C_{2 \epsilon}(M_j)$.  Because
$g_{\tilde{x}}$ intersects $\tilde{B}_j$, there is a portion of 
$\bar{g}_{\tilde{x}}$ that contains $\tilde{x}$ and joins two components of $\partial
C_{2 \epsilon}(M_j)$.  Because $C_{2
\epsilon}(M_j)$ is convex, $\tilde{x} \in g_{\tilde{x}}$ must lie in $C_{2
\epsilon}(M_j)$, which is a contradiction.

Then all but a compact portion of $g_{\tilde{x}}$ is contained in
$\tilde{U}_j$.  Using the coordinates $(\theta_i,r_i,s_i)$ on each component of
${\mathcal P}$ adjacent to $\tilde{U}_j^\prime$, we can extend the homeomorphism
$p_j|_{\tilde{U}_j^\prime}$ to a homeomorphism
$p_j|_{(\tilde{U}_j^\prime)_P}:(\tilde{U}_j^\prime)_P
\rightarrow (U_j^\prime)_P$ where $(\tilde{U}_j^\prime)_P$ and $(U_j^\prime)_P$ are the
parabolic extensions of
$\tilde{U}_j^\prime$ and $U_j^\prime$ respectively.  Let
$\tilde{V}_j^\prime$ be the closure of the component of
$M_j - \tilde{A}_j$ that contains $g_{\tilde{x}}$.  Let
$X_\delta = [M_j - C_\delta(M_j)] \cap \tilde{V}_j^\prime$.  Let $\bar{U}_j^\prime$ be
the closure of $\tilde{U}_j^\prime$ in $M_j$.  Then because $\partial \bar{U}_j^\prime -
\partial
\tilde{\mathcal P}_j$ is compact, and
$g_{\tilde{x}}
\cap X_\delta
\neq
\emptyset$ for every $\delta > 0$, we can guarantee that
$X_\delta
\subset (\tilde{U}_j^\prime)_P$ for large enough $\delta$.   Then for large enough
$\delta$, $X_\delta$ is a component of $M_j - C_\delta(M_j)$ which embeds in
$N$.  Then by Lemma \ref{GFproj}, $V_j^\prime = p_j(\tilde{V}_j^\prime)$ is the
component of
$N - C_{2 \epsilon}(N)$ with boundary
$p_j(\tilde{A}_j)$.  Because $\tilde{x} \in M_j - C_{2 \epsilon}(M_j)$,
$\tilde{x} \in int \ \tilde{V}_j^\prime$.  Then
$p_j(\tilde{x})
\in int \ V_j^\prime \subset N - C_{2 \epsilon}(N)$.  But by hypothesis,
$x = p_j(\tilde{x}) \in C_{2 \epsilon}(N)$ so this is a contradiction.  Thus, for each
$j$, 
$\tilde{U}_j \subset C_{2 \epsilon}(M_j)$.

%

For each $j$, choose $\tilde{F}_j
\subset M_j$ to be a level surface in the product structure of $\tilde{U}_j$ such that
the compact surface $\tilde{H}_2(\bar{S}_j, 2-\eta)$ lies to one side of
$\tilde{F}_j$. Let
$A_{\tilde{P}_A,j}$ be embedded annuli in
$\tilde{P}_A \cap C_{2 \epsilon}(M_j)$ with boundary components
$\partial \tilde{F}_j \cap \partial \tilde{P}_A$ and
$\tilde{H}_2(\gamma,2-\eta)$.  Also, for each component of $\tilde{P}_S$, let
$A_{\tilde{P}_S,j} = (\phi|_{\tilde{P}_S})^{-1}(\partial
\tilde{F}_j \cap \partial \tilde{P}_S)$ be infinite annuli in
$\tilde{P}_S$.  Note that because $\partial \tilde{F}_j \cap \partial \tilde{P}_S
\subset C_{2 \epsilon}(M_j)$ and $C_{2 \epsilon}(M_j)$ is convex, $A_{\tilde{P}_S,j}
\subset C_{2
\epsilon}(M_j)$.
Let 
$\tilde{\mathcal F}_j = \tilde{F}_j \cup A_{\tilde{P}_A,j} \cup A_{\tilde{P}_S,j}$.  By
construction, the surface $\tilde{\mathcal F}_j$ is embedded in $C_{2 \epsilon}(M_j)$ and
properly homotopic in $M_j$ to $\tilde{H}_2(\bar{S}_j,2-\eta)$.

For each $j$, we can define $\alpha_j: \bar{S}_j \times [2,3] \rightarrow M_j$ to be a
proper homotopy between
$\tilde{H}_2(\bar{S}_j,2-\eta)$ and $\tilde{\mathcal F}_j$ so that
$\alpha_j(\bar{S}_j,2) =
\tilde{H}_2(\bar{S}_j,2-\eta)$ and $\alpha_j(\bar{S}_j,3) = \tilde{\mathcal F}_j$.  Let
$\beta_j: \bar{S}_j \times [2,3] \rightarrow M_j$ be a homotopy constructed from
$\alpha_j$ such that for $x \in \bar{S}_j$, $\beta_j(x,[2,3])$ is a geodesic arc with
the same endpoints and in the same homotopy class as $\alpha_j(x,[2,3])$. Then
$\beta_j$ is a {\it ruled homotopy\/} between
$\tilde{H}_2(\bar{S}_j,2-\eta)$ and
$\tilde{\mathcal F}_j$.  Because $C_{2 \epsilon}(M_j)$ is convex,
$\beta_j(\bar{S}_j,[2,3]) \subset C_{2 \epsilon}(M_j)$.  By the induction hypothesis,
for $x \in \beta_j(\bar{S}_j,[2,3])$, $inj_{M_j}(x) \leq L_{S_j} + 2 \epsilon$.  

For each $j$, let $\bar{F}_j = p_j(\tilde{F}_j)$ so that $\bar{F}_j$ is a level
surface in the product structure of $U_j$, let
$A_{P_A,j} = p_j(A_{\tilde{P}_A,j})$, let $A_{P_S,j} = p_j(A_{\tilde{P}_S,j})$, and let
$\bar{\mathcal F}_j^\prime = \bar{F}_j \cup A_{P_A,j} \cup A_{P_S,j}$.  Then
$\bar{\mathcal F}_j^\prime = p_j(\tilde{\mathcal F}_j)$.

For each $j$, we can let $\delta_j = p_j
\circ \beta_j$.  This map will be a proper homotopy between
$\delta_j(\bar{S}_j,2) = H_2(\bar{S}_j,2-\eta)$ and
$\delta_j(\bar{S}_j,3) = \bar{\mathcal F}_j^\prime$, where
$\delta_j(\bar{S}_j,[2,3]) \subset C_{2 \epsilon}(N)$.  By Lemma
\ref{coverinj}, for $x \in
\delta_j(\bar{S}_j,[2,3])$, $inj_N(x) \leq L_{S_j} + 2 \epsilon$.  

Now let $T_1^\prime = (\bigcup_j \bar{\mathcal F}_j^\prime)$ glued along
$H_2(\gamma,2-\eta)$.  

We can construct the analogous surface $T_2^\prime = (\bigcup_k
\bar{\mathcal F}_k^\prime)$ by repeating this process so that each $\bar{\mathcal
F}_k^\prime$ is the union of embedded annuli $A_{P_B,k}$ in
$P_B$, infinite annuli $A_{P_S,k}$ in $P_S$, and a level surface $\bar{F}_k$ in the
product structure of
$U_k$ such that $\bar{F}_k$ lies to one side of $H_3(\bar{S}_k,-1 + \eta)$.  

For $i= 1,2$, the surfaces $T_i^\prime$ are embedded and properly homotopic to
${\mathcal T}_i$ by construction.  Construct a proper homotopy $\zeta_1: S
\times [2,3]
\rightarrow N$ between $\zeta_1(S,2) = {\mathcal T}_1$ and $\zeta_1(S,3) =
T_1^\prime$ such that $\zeta_1|_{\bar{S}_j} = \delta_j$ for all $j$.  By construction,
the image of $\delta_j$ lies in $p_j(C_{2 \epsilon}(M_j)) \subset C_{2 \epsilon}(N)$. 
Thus for $x \in \zeta_1(S,[2,3])$, $inj_N(x) \leq max \{L_{S_j} + 2 \epsilon \}$.
Construct the analogous proper homotopy $\zeta_2: S \times [-2,-1] \rightarrow N$
between $\zeta_2(S,-1) = {\mathcal T}_2$ and $\zeta_2(S,-2) = T_2^\prime$ such that
$\zeta_2|_{\bar{S}_k} = \delta_k$ for all $k$.

Let $T_1 = T_1^\prime - (\bigcup A_{P_S,j})$  Then by construction, $T_1 = \bar{F}_j \cup A_{P_A,j}$.  Because
$p_j(\tilde{P}_A) = P_A$, we can conclude that
$T_1 \cap U_j = T_1 \cap U_j^\prime = \bar{F}_j$.  In particular,
$\partial \bar{F}_j \subset
\partial {\mathcal P}$.   Also, by construction,
$(p_j|_{\tilde{U}_j})^{-1}(\bar{F}_j) = (p_j|_{\tilde{U}_j^\prime})^{-1}(\bar{F}_j)
\subset
\overline{C_{2 \epsilon}(M_j)}$.  The analogous conclusions hold for  $T_2 = T_2^\prime
- (\bigcup A_{P_S,k})$ and $T_2 \cap U_k = \bar{F}_k$ for each $k$.

By construction $T_1$ and
$T_2$ are properly homotopic  in
$N - P_S$ to the components of $\partial R^* - P^*_S$, where $P^*_S$ denotes the
components of $P^*$ associated to $P_S$.  So
$T_1$, $T_2$ and the boundary components of $R^* - P^*_S$ span a product structure in
$C_{2
\epsilon}(N)$. (Thm 10.5, Hempel
\cite{hempel}) Hence, $T_1$ and $T_2$ are boundary components of a new set $R \subset
C_{2 \epsilon}(N)$ which is a compact core of $N - P_S$.

Now we will show that $$R \subset [W \cup \zeta_1(S,[2,3]) \cup
\zeta_2(S,[-2,-1])]$$  Because $T_1^\prime$ and $T_2^\prime$ are embedded,
incompressible surfaces that are properly homotopic in $N$, we can conclude that $T_1$
and $T_2$ are embedded, incompressible surfaces that are properly homotopic in $N -
P_S$. So $T_1$ and $T_2$ span a product structure in $N - P_S$. (Thm 10.5, Hempel
\cite{hempel}) Then there exists a proper isotopy ${\nu}: \bar{S} \times [0,1]
\rightarrow N$ between $T_1$ and $T_2$ such that ${\nu}$ is an embedding in $N - P_S$.
 Using the inverse map $(\phi|_{P_S})^{-1}$, we can extend
$\bar{\nu}$ to a proper isotopy $\nu: S \times [0,1] \rightarrow N$ between
$T_1^\prime$ and
$T_2^\prime$ such that $\nu$ is also an embedding.  So
$\nu(S,[0,1])$ is homeomorphic to
$S
\times [0,1]$. Using standard degree arguments, 
$\nu(S,[0,1])$ is contained in the image of any proper homotopy between $T_1^\prime$
and $T_2^\prime$. (see Thm 2.14, Lloyd \cite{lloyd})  In particular, $R
\subset
\nu(S,[0,1])
\subset [W
\cup \zeta_1(S,[2,3]) \cup
\zeta_2(S,[-2,-1])]$.

So for $x \in R$, $inj_N(x) \leq max \{ K_S, L_{S_j} + 2\epsilon, L_{S_k} + 2\epsilon
\}$. 

Note that in Case IB, we must use Lemma \ref{para to para} to find $\eta$, and replace
the maps $H_i$ and $\tilde{H}_i$ for
$i=2,3$ with $H_1$ and $\tilde{H}_1$, keeping the domains the same.  This completes the
proof of Lemma
\ref{surface core}. 
\end{proof}

Now we have bounded the injectivity radius for points in a compact core $R$ of
$N - P_S$ which contains $R^*$.

Consider again the  
surface constructed in the previous lemma 
$T_1 \cap U_j^\prime = \bar{F}_j \subset \partial R$ which was a $\pi_1$-injective,
separating, level surface in the product structure of $U_j$.  Let
$W_j^\prime$ be the closure of the component of
$(N - {\mathcal P}) -
\bar{F}_j$ such that $W_j^\prime \subset U_j^\prime$. 
Let $\tilde{W}_j^\prime$ be the subset of $\tilde{U}_j^\prime$ such that
$p_j|_{\tilde{W}_j^\prime}: \tilde{W}_j^\prime \rightarrow W_j^\prime$ is a
homeomorphism.  Then because $\partial
\bar{F}_j \subset
\partial {\mathcal P}$ and 
$(p_j|_{\tilde{W}_j^\prime})^{-1}(\bar{F}_j) =
(p_j|_{\tilde{U}_j^\prime})^{-1}(\bar{F}_j)
\subset \overline{C_{2 \epsilon}(M_j)}$, by Lemma
\ref{not in C(N)}, we know that  $(p|_{\tilde{W}_j^\prime})^{-1}[(N - C_{2
\epsilon}(N))
\cap W_j^\prime] = (M_j - C_{2 \epsilon}(M_j)) \cap \tilde{W}_j^\prime$.  
The analogous conclusions hold for
$T_2
\cap U_k = \bar{F}_k \subset \partial R$ and $W_k^\prime$, the
closure of the component of
$(N - {\mathcal P}) -
\bar{F}_k$ such that $W_k^\prime \subset U_k^\prime$.

If $x \in [C(N) - {\mathcal P}] - R$, then
$x
\in W_j^\prime
\cap C(N)$ for some
$j$, or $x \in W_k^\prime \cap C(N)$ for some $k$. 
Because $(p|_{\tilde{W}_j^\prime})^{-1}[(N - C_{2
\epsilon}(N))
\cap W_j^\prime] = (M_j - C_{2 \epsilon}(M_j)) \cap \tilde{W}_j^\prime$, we know that if
$x \in W_j^\prime \cap C(N)$, then $x \in p_j(C(M_j))$.  So by the
induction hypothesis and Lemma
\ref{coverinj}, $inj_N(x) \leq L_{S_j}$.  Similarly, if $x \in U_k^\prime \cap C(N)$,
then
$inj_N(x) \leq L_{S_k}$.

Finally, if $x \in C(N) \cap {\mathcal P}$, then $inj_N(x) \leq \epsilon$. 

Therefore, for $x \in C(N)$, $inj_N(x) \leq max \{K_S, L_{S_j} + 2
\epsilon, L_{S_k} + 2 \epsilon\}$.   We can do this for all 
$0 <  \epsilon < \epsilon_3$, so that for $x \in C(N)$, $inj_N(x)
\leq max \{K_S, L_{S_j}, L_{S_k}\}$.  This completes the proof of Case I. \\

Case II.  Accidental Parabolics Associated Only to $B$. \\

Now let us assume that there are accidental parabolics associated
only to $B$.  Then  $[P \cap \bar{h}(\bar{S},0)]$ is empty so that the end $E$
of $\NOE - \hat{R}$ associated to $\bar{h}(\bar{S},0)$ is either simply degenerate or
geometrically finite. \\

Case IIA.  End $E$ Is Simply Degenerate. \\

We begin with a sketch of the proof of Case IIA.  Again we will divide the convex
core into a compact core and its complement.  Using a continuous family of simplicial
hyperbolic surfaces and projections of convex cores of covers with fewer accidental
parabolics, points in the compact core will have bounded injectivity radius.  The
remaining portions will lie in an $\epsilon$-thin part, lie in the projection of a
cover which has fewer accidental parabolics, or lie in the image of a simplicial
hyperbolic surface.  This completes the sketch of the proof in this case.

The first step in the proof of Case IIA will be to find a continuous family of
simplicial hyperbolic surfaces which shrinks an accidental parabolic associated to $B$.
We will want points in the image of this continuous family to lie in the convex core and
to have bounded injectivity radius.

In Theorem 4.1 \cite{canary2}, Canary has shown that if $E$ is a simply degenerate end
of $N^\circ_\epsilon$, then there exists a sequence of useful simplicial hyperbolic
surfaces which are properly homotopic to $h: S \rightarrow N$ and which exit every
compact set in $E$.  Recall that $\hat{R}$, a relative compact core of
$N^\circ_\epsilon$ with associated parabolic locus $\hat{P}$, is homeomorphic to
$\bar{S}
\times [0,1]$, where
$\bar{S}$ is a compact core of $S$.  Therefore, the sequence of useful simplicial
hyperbolic surfaces have the form $f: (S,T) \rightarrow N$.  Let
$f: (S,T_A)
\rightarrow N$ be a useful simplicial hyperbolic surface in
$E$ with distinguished edge $\hat{e}_1$.

In the triangulation $T_A$, either there exists an edge which represents an accidental
parabolic associated to $B$ or not.  In the latter case, let $T_2$ be a triangulation
of
$S$ such that there exists an edge in
$T_2$ which is an accidental parabolic associated to $B$.   We can perform the modified
$HH$ sequence of elementary moves altering $T_A$ to $T_2$.  Let $T_B$ be the ``first''
triangulation in the sequence in which there is an edge $\gamma_B$ representing an
accidental parabolic associated to $B$.  There are two cases:  (1) $T_A$ and $T_B$ are
adjacent, or (2) they are separated by triangulations whose edges are not accidental
parabolics associated to $B$.

In either of these cases, we say that the sequence of triangulations has
$Property \ T^B$, that is, there exist triangulations $T_A$ and $T_B$ in our
sequence such that: 
\begin{enumerate}
\item the edges of $T_A$ do not represent accidental parabolics,
\item the edges of $T_B$ may only represent accidental parabolics associated to $B$, and
\item the edges of each intermediary triangulation do not represent accidental
parabolics associated to $B$.
\end{enumerate}
Note that when $T_A$ and $T_B$ are adjacent, there are no intermediary triangulations.

Thus, we can see that either a sequence of triangulations has $Property \ T^B$, or
there exists an edge in the triangulation $T_A$ which represents an accidental
parabolic associated to
$B$.
\\

Case IIA(1).   There Exists a Sequence of Triangulations with $Property \ T^B$. \\

Suppose our sequence of triangulations has $Property \ T^B$.  Now there are two cases: 
either $T_A$ and $T_B$ are separated by at least one intermediary triangulation, or
not. \\

Case IIA(1)a. $T_A$ and $T_B$ Separated by at Least one Intermediary Triangulation.
\\

In this case, there exist triangulations
$T_A$ and $T_B$ such that no edge of $T_A$ is an accidental parabolic, and exactly one
edge $\gamma_B$ of $T_B$ is an accidental parabolic associated to $B$, such that the
edges of each intermediary triangulation are not accidental parabolics associated to
$B$. The fact that there is exactly one edge of $T_B$ associated to an
accidental parabolic follows directly from the fact that an elementary move alters
exactly one edge in the triangulation.

If $S$ is a 4-punctured sphere, then there is exactly one nonperipheral edge in its
suitable, one-vertex triangulation.  Thus, in
$T_A$, this nonperipheral edge is
$\hat{e}_1$; in $T_B$, this nonperipheral edge is $\gamma_B$.  Let
$\hat{e}_1^\prime$ and $\hat{e}_2^\prime$ be any two edges in $T_A$ joining distinct
ideal vertices to the sole internal vertex.  We can use Lemma
\ref{makepractical} to construct a $\pi_1$-injective practical simplicial hyperbolic
surface
$f^\prime:(S,T_A)
\rightarrow N$ with marked edges $\hat{e}_1^\prime$ and $\hat{e}_2^\prime$ such that
$f^\prime$ is properly homotopic to $h:S \rightarrow N$.   Let
$\hat{e}_3^\prime$ and $\hat{e}_4^\prime$ be any two edges in $T_B$ joining distinct
ideal vertices to the sole internal vertex.   Similarly, we can construct a
$\pi_1$-injective practical simplicial hyperbolic surface $g:(S,T_B)
\rightarrow N$ with marked  edges $\hat{e}_3^\prime$ and $\hat{e}_4^\prime$ such that
$g$ is properly homotopic to $h:S \rightarrow N$.  

First we can use Lemma \ref{usefulpractical} to construct a continuous family of
simplicial hyperbolic surfaces which interpolates between the useful simplicial
hyperbolic surface
$f:(S,T_A) \rightarrow N$ with distinguished edge $\hat{e}_1$ and the practical
simplicial hyperbolic surface $f^\prime: (S,T_A) \rightarrow N$ with marked edges
$\hat{e}_1^\prime$ and $\hat{e}_2^\prime$. Then by repeatedly using Lemmas \ref{3 elem
move} and \ref{3 edge} and concatenating, we can
construct the $HH$ interpolation 
$H_1: S
\times [0,1] \rightarrow N$ such that
$H_1|_{(S,0)}:S \rightarrow N$ is the map $f:(S,T_A) \rightarrow N$ and 
$H_1|_{(S,1)}:S \rightarrow N$ is the map $g:(S,T_B) \rightarrow N$.  

Now we can use Lemma \ref{3 length to 0}, setting $H_1|_{(S,1)} = H_2|_{(S,0)}$, to
obtain a continuous family of simplicial hyperbolic surfaces $H_2: S \times [1,2)
\rightarrow N$ that shrinks the accidental parabolic $\gamma_B$.

If $S$ is not a 4-punctured sphere, then  there
exists a nonperipheral edge $\hat{e}_2 \neq \gamma_B$ in $T_B$ which does not represent
an accidental parabolic, and hence can be mapped to a closed geodesic.
Then by Lemma \ref{makeuseful}, we can construct a useful simplicial
hyperbolic surface $g:(S,T_B) \rightarrow N$ with distinguished edge $\hat{e}_2$.

By repeatedly using Lemma \ref{elem move}, we can
construct a continuous family of simplicial hyperbolic surfaces that performs the
elementary moves indicated in the modified $HH$ sequence of
triangulations that interpolates between $T_A$ and $T_B$.  We can also use Lemma
\ref{edge} to construct a continuous family that alters the distinguished edge from
$\hat{e}_1$ to $\hat{e}_2$ as needed in the modified $HH$ sequence.  By concatenating
all of these continuous families, we have the
$HH$ interpolation 
$H_1: S
\times [0,1] \rightarrow N$ such that
$H_1|_{(S,0)}:S \rightarrow N$ is the map $f:(S,T_A) \rightarrow N$ and 
$H_1|_{(S,1)}:S \rightarrow N$ is the map $g:(S,T_B) \rightarrow N$.

Now we can use Lemma \ref{length to 0}, setting
$H_1|_{(S,1)} = H_2|_{(S,0)}$, to obtain a continuous family of simplicial hyperbolic
surfaces $H_2: S
\times [1,2)
\rightarrow N$ that shrinks the accidental parabolic $\gamma_B$.\\

Case IIA(1)b.  $T_A$ and $T_B$ Not Separated by at Least one Intermediary
Triangulation. 
\\

In this case, we can conclude that $T_A$ and $T_B$ are adjacent triangulations, that is,
the elementary move alters an edge of $T_A$ to the edge $\gamma_B$. 

If  $S$ is a 4-punctured sphere, then because its suitable, one-vertex triangulation
contains exactly one nonperipheral edge,  we know that in $T_A$,
$\hat{e}_1$ is that nonperipheral edge, and in $T_B$, $\gamma_B$ is that nonperipheral
edge.  Let
$\hat{e}_1^\prime$ and $\hat{e}_2^\prime$ be any two edges in $T_A$ joining distinct
ideal vertices to the sole internal vertex.  We can use Lemma
\ref{makepractical} to construct a $\pi_1$-injective practical simplicial hyperbolic
surface
$f^\prime:(S,T_A)
\rightarrow N$ with marked edges $\hat{e}_1^\prime$ and $\hat{e}_2^\prime$ such that
$f^\prime$ is properly homotopic to $h:S \rightarrow N$.   Let
$\hat{e}_3^\prime$ and $\hat{e}_4^\prime$ be any two edges in $T_B$ joining distinct
ideal vertices to the sole internal vertex.  Similarly, construct a $\pi_1$-injective
practical simplicial hyperbolic surface $g:(S,T_B)
\rightarrow N$ with marked  edges $\hat{e}_3^\prime$ and $\hat{e}_4^\prime$ such that
$g$ is properly homotopic to $h:S \rightarrow N$.  

We can use Lemma \ref{usefulpractical} to construct a continuous family of simplicial
hyperbolic surfaces which interpolates between the useful simplicial hyperbolic surface
$f:(S,T_A) \rightarrow N$ with distinguished edge $\hat{e}_1$ and the practical
simplicial hyperbolic surface $f^\prime: (S,T_A) \rightarrow N$ with marked edges
$\hat{e}_1^\prime$ and $\hat{e}_2^\prime$. By using Lemma
\ref{3 elem move}, we can construct a continuous family of simplicial hyperbolic
surfaces that performs the elementary move that separates the triangulations
$T_A$ and
$T_B$.   We can also use Lemma \ref{3 edge} to construct a continuous family  that
alters the marked edges from $\hat{e}_1^\prime$ and $\hat{e}_2^\prime$ to
$\hat{e}_3^\prime$ and $\hat{e}_4^\prime$.   After concatenating all of these continuous
families, we construct  the
$HH$ interpolation 
$H_1: S
\times [0,1] \rightarrow N$ such that
$H_1|_{(S,0)}:S \rightarrow N$ is the map $f:(S,T_A) \rightarrow N$ and 
$H_1|_{(S,1)}:S \rightarrow N$ is the map $g:(S,T_B) \rightarrow N$.

Now we can use Lemma \ref{3 length to 0}, setting $H_1|_{(S,1)} = H_2|_{(S,0)}$, to
obtain a continuous family of simplicial hyperbolic surfaces $H_2: S \times [1,2)
\rightarrow N$ that shrinks the accidental parabolic $\gamma_B$.

If $S$ is not a 4-punctured sphere, then there exists a nonperipheral edge
$\hat{e}_2 \neq \gamma_B$ in
$T_B$ which is unaltered by the sole
elementary move and hence can be mapped to a closed geodesic.  Therefore, we can
use Lemma
\ref{makeuseful} to construct
a $\pi_1$-injective useful simplicial hyperbolic surface $g:(S,T_B) \rightarrow N$ with
distinguished edge $\hat{e}_2$ such that $g$ is
properly homotopic to $h:S \rightarrow N$.  

We can use the Lemma \ref{elem move} to perform the elementary move that separates the
triangulations $T_A$ and $T_B$.  We can also use Lemma \ref{edge} to alter the
distinguished edge from $\hat{e}_1$ to $\hat{e}_2$.  After concatenating all of these
continuous families, we can  construct the
$HH$ interpolation 
$H_1: S
\times [0,1] \rightarrow N$ such that
$H_1|_{(S,0)}:S \rightarrow N$ is the map $f:(S,T_A) \rightarrow N$ and 
$H_1|_{(S,1)}:S \rightarrow N$ is the map $g:(S,T_B) \rightarrow N$.

Now we can use Lemma \ref{length to 0}, setting $H_1|_{(S,1)} = H_2|_{(S,0)}$, to obtain
a continuous family of simplicial hyperbolic surfaces $H_2: S \times [1,2) \rightarrow
N$ that shrinks the accidental parabolic $\gamma_B$.\\

Case IIA(2).  There Does Not Exist A Sequence of Triangulations with $Property \ T^B$.
\\

If we cannot construct a sequence of triangulations with $Property \ T^B$, then as we
argued before, there exists an edge $\gamma_B$ in the triangulation $T_A$ that
represents an accidental parabolic associated to $B$.  Using Lemma
\ref{length to 0} we can construct a continuous family of simplicial hyperbolic
surfaces $H_1: S
\times [0,2) \rightarrow N$ that shrinks the accidental parabolic $\gamma_B$. \\

Let $W = H_1(S,[0,1]) \cup H_2(S,[1,2))$ (or $H_1(S,[0,2))$ in Case IIA(2)) be the union
of the images of the continuous families of simplicial hyperbolic surfaces.   
Lemmas \ref{edge}, \ref{elem move}, \ref{length to 0}, \ref{3 edge}, \ref{3 elem move},
\ref{3 length to 0}, and \ref{usefulpractical}
 guarantee that $W \subset C(N)$. 
 Furthermore, we
also know that for any $x \in W$ lies in the image of a
simplicial hyperbolic surface, and hence, by Lemma
\ref{shsinj}, $inj_N(x) \leq K_S$ where the constant $K_S$ depends on the Euler
characteristic of $S$. Now we have bounded the injectivity radius in a portion of the
convex core.  Note that $\phi(W) \subset C(N)$ because $C(N)$ is convex, and that
by a result of Canary (Lem 4.4, \cite{canary1}) and Lemmas \ref{length to 0} and \ref{3
length to 0},
$\overline{\phi(W)}$ is compact. 

The next step in the proof of Case IIA involves finding a compact
core $R$ of $N - P_S$, where $P_S$ represents the collection of rank one cusps of
$\NOE$ associated to $\partial S$, such that points in $R$ have bounded injectivity
radius.

Let $\xi: S \rightarrow N$
denote the map $H_1|_{(S,1)}: S \rightarrow N$.  Let \{$S_{j}$\} be the components of
$S - \gamma_B$.  Let \{$M_j=\HHH / \xi_*(\pi_1(S_{j}))$\} be covers of $N$ with covering
maps \{$p_j$\}.  Let $P_B$ be the rank one cusp of $\NOE$ associated
to $\gamma_B$.  Let ${\mathcal P} = P_B \cup P_S$. 

Given $0 < \epsilon < \epsilon_3$.  By  Lemma \ref{compactincore}, there exists a
relative compact core
$R^*$ of
$\overline{C_{\epsilon}(N)} - {\mathcal P}$ such that $\partial C_\epsilon(N) \cap
C^\circ_\epsilon(N)
\subset
\partial R^*$ and such that $\overline{\phi(W)} \subset R^*$.  Let
\{$U_j^\prime$\} be the components of $[N - {\mathcal P}] - int \ R^*$
such that each $U_j^\prime$ is adjacent to
$P_A$, and let $U_k^\prime$ be the component of $[N - {\mathcal P}] -
int \ R^*$ which is a neighborhood of the simply degenerate end $E$.  Let
\{$U_j$\} be the components of $[C_{2 \epsilon}(N) - {\mathcal P}] - int \ R^*$ such
that each $U_j$ is adjacent to
$P_B$.  Then $U_j \subset U_j^\prime$. Let $U_k$ be the component of $[C_{2 \epsilon}(N)
- {\mathcal P}] - int
\ R^*$ such that $U_k \subset U_k^\prime$.

In the following lemma,  we will choose a compact core
$R$ of $N - P_S$ such that $R$ contains $R^*$,
$R$ lies in $C_{2 \epsilon}(N)$, and points in $R$ have bounded injectivity radius. 

\begin{lemma}
\label{surface core A}  For all $0 < \epsilon < \epsilon_3$, there exists a compact
core $R$ of $N - P_S$ bounded by surfaces $T_1$ and $T_2$ such that the following are
true: 
\begin{enumerate}
\item $R^* \subset R \subset C_{2 \epsilon}(N)$,
\item if $x \in R$, then $inj_N(x) \leq max \{ K_S, L_{S_j} + 2 \epsilon\}$
\item for each $j$, $\partial (T_1 \cap U_j^\prime) \subset \partial {\mathcal P}$, and
\item for each $j$, there exist $\tilde{U}_j^\prime \subset M_j$ such that
$p_j|_{\tilde{U}_j^\prime}: \tilde{U}_j^\prime \rightarrow U_j^\prime$ is a
homeomorphism, and
$(p_j|_{\tilde{U}_j^\prime})^{-1}(T_1 \cap U_j^\prime) \subset \overline{C_{2
\epsilon}(M_j)}$. 
\end{enumerate}
\end{lemma}

\begin{proof} The proof of this lemma mimics that of Lemma \ref{surface core}.  We
give an outline of the steps.  Let ${\mathcal T}_1 = H_2(S,1-\eta)$, where
$H_2(S_j,1-\eta)$ lifts to
$C_\epsilon(M_j)$, and let ${\mathcal T}_2 = H_1(S,0)$.
By Lemma \ref{prodstructure}, we know that $U_k^\prime$, $U_k$,
and each $U_j^\prime$ and $U_j$ possesses a product structure. 

Construct a surface $T_1^\prime = \bigcup_j \bar{\mathcal F}_j^\prime$ with the
following properties: \ (1) each $\bar{\mathcal F}_j^\prime$ is the union of embedded
annuli
$A_{P_B,j}$ in $P_B \cap C_{2 \epsilon}(N)$, embedded infinite annuli $A_{P_S,j}$
in
$P_S
\cap C_{2 \epsilon}(N)$, and a level surface
$\bar{F}_j$ in the product structure on $U_j$ such that
${H}_2(\bar{S}_j,1-\eta)$ lies to one side of
$\bar{F}_j$; and (2) the surface $T_1^\prime$ is properly homotopic to
${\mathcal T}_1$ via a proper homotopy $\zeta_1$ whose image lies in
$\bigcup p_j(C_{2 \epsilon}(M_j))$.  Then for a point $x$ in the image of $\zeta_1$,
$inj_N(x) \leq max \{L_{S_j} + 2 \epsilon\}$. 

Let $T_2^\prime$ be the union of a level surface $\bar{F}_k$ in the product structure
of $U_k$ such that ${\mathcal T}_2$ lies to one side, and infinite annuli $A_{P_S,k}$ in
$P_S$.  Let
$\zeta_2$ be a proper homotopy between $T_2^\prime$ and ${\mathcal T}_2$.  

Consider a
point $x$ that is in the image of $\zeta_2$, but not in the image of
$\zeta_1$ nor $H_1$. 
Note that $U_k$ has a product structure, that
${\mathcal T}_2$ is disjoint from
$U_k$, and that ${\mathcal T}_2$ is properly homotopic to $\partial (U_k)_P$.  Then by a
result of Bonahon (Prop 1.23,
\cite{bonahon}) and standard degree arguments (see Thm 2.14, Lloyd, \cite{lloyd}), there
exists a useful simplicial hyperbolic surface
$g:S \rightarrow N$ with image in $(U_k)_P$ such that $x$ lies in the image of any
proper homotopy between
${\mathcal T}_2$ and $g(S)$.  


Then, in particular, $x$ lies in the image of the $HH$ 
interpolation $J_1: S \times [0,1]
\rightarrow N$ joining $H_1|_{(S,0)}$ and $g$.    Lemmas \ref{edge}, \ref{elem move},
\ref{3 edge}, \ref{3 elem move}, and \ref{usefulpractical}
 guarantee that
$J_1(S,[0,1])
\subset C(N)$.  Furthermore, by Lemma \ref{shsinj}, if $x \in J_1(S,[0,1]))$, 
$inj_N(x) \leq K_S$.  Therefore, for a point $x$ that is contained in the image of
$\zeta_2$ but that is not contained in the image of $\zeta_1$ nor $H_1$,
$inj_N(x)
\leq K_S$.  


Let $T_1 = T_1^\prime - (\bigcup A_{P_S,j})$ and $T_2 = T_2^\prime - (\bigcup
A_{P_S,j})$.  We can show that $T_1 \cap U_j = T_1 \cap U_j^\prime = \bar{F}_j$ and that
$\partial (T_1
\cap U_j) =
\partial
\bar{F}_j
\subset \partial {\mathcal P}$.  Furthermore, we can show that there exists
$\tilde{U}_j^\prime
\subset M_j$ such that
$p_j|_{\tilde{U}_j^\prime}: \tilde{U}_j^\prime \rightarrow U_j^\prime$ is a
homeomorphism, and
$(p_j|_{\tilde{U}_j})^{-1}(\bar{F}_j)
= (p_j|_{\tilde{U}_j^\prime})^{-1}(\bar{F}_j)
\subset \overline{C_{2 \epsilon}(M_j)}$. 


Then we can show that $T_1$ and $T_2$ bound a compact core $R$ of $N - P_S$ such that
$R^* \subset R \subset C_{2 \epsilon}(N)$.  By covering $R$ with $W$ and the images of
$\zeta_1$ and $\zeta_2$, we can show that for
$x
\in R$,
$inj_N(x)
\leq max
\{K_S, L_{S_j} + 2
\epsilon\}$.

Note that in Case IB, we must replace
the maps $H_2$ and $\tilde{H}_2$ with $H_1$ and $\tilde{H}_1$, keeping the domains the
same.  This completes the sketch of the proof of Lemma \ref{surface core A}. 
\end{proof}

Now we have bounded the injectivity radius for points in a compact core $R$ of $N -
P_S$ which contains $R^*$.

Consider again the  
surface constructed in the previous lemma 
$T_1 \cap U_j^\prime = \bar{F}_j \subset \partial R$ which was a $\pi_1$-injective,
separating, level surface in the product structure of $U_j$.  Let
$W_j^\prime$ be the closure of the component of
$(N - {\mathcal P}) -
\bar{F}_j$ such that $W_j^\prime \subset U_j^\prime$.  Let
$W_j^\prime$ be the closure of the component of
$(N - {\mathcal P}) -
\bar{F}_j$ such that $W_j^\prime \subset U_j^\prime$. 
Let $\tilde{W}_j^\prime$ be the subset of $\tilde{U}_j^\prime$ such that
$p_j|_{\tilde{W}_j^\prime}: \tilde{W}_j^\prime \rightarrow W_j^\prime$ is a
homeomorphism.  Then because $\partial
\bar{F}_j \subset
\partial {\mathcal P}$ and 
$(p_j|_{\tilde{W}_j^\prime})^{-1}(\bar{F}_j) =
(p_j|_{\tilde{U}_j^\prime})^{-1}(\bar{F}_j)
\subset \overline{C_{2 \epsilon}(M_j)}$, by Lemma
\ref{not in C(N)}, we know that  $(p|_{\tilde{W}_j^\prime})^{-1}[(N - C_{2
\epsilon}(N))
\cap W_j^\prime] = (M_j - C_{2 \epsilon}(M_j)) \cap \tilde{W}_j^\prime$. 

If $x \in [C(N) - {\mathcal P}] - R$, then
$x
\in W_j^\prime
\cap C(N)$ for some
$j$, or $x \in W_k^\prime \cap C(N)$. 
Because $(p|_{\tilde{W}_j^\prime})^{-1}[(N - C_{2
\epsilon}(N))
\cap W_j^\prime] = (M_j - C_{2 \epsilon}(M_j)) \cap \tilde{W}_j^\prime$, we know that if
$x \in W_j^\prime \cap C(N)$, then $x \in p_j(C(M_j))$.  So by the
induction hypothesis and Lemma
\ref{coverinj}, $inj_N(x) \leq L_{S_j}$.

Suppose $x \in W_k^\prime \cap C(N)$.   Recall that $W_k^\prime$ has a product
structure, that
${\mathcal T}_2$ is disjoint from
$W_k^\prime$, and that ${\mathcal T}_2$ is properly homotopic to $\partial
(W_k^\prime)_P$. Then by a result of Bonahon (Prop 1.23,
\cite{bonahon}) and standard degree arguments (see Thm 2.14, Lloyd,
\cite{lloyd}), there exists a useful simplicial hyperbolic surface
$h:S \rightarrow N$ such that $x$ lies in the image of any proper homotopy between 
${\mathcal T}_2$ and $h(S)$.  Then, in particular, $x$ lies in the image of the $HH$
interpolation $J_2: S
\times [0,1]
\rightarrow N$ joining $H_1|_{(S,0)}$ and $h$.   Lemmas \ref{edge}, \ref{elem move},
\ref{3 edge}, \ref{3 elem move}, and \ref{usefulpractical} guarantee that
$J_2(S,[0,1])
\subset C(N)$.  Furthermore, by Lemma \ref{shsinj}, if $x \in J_2(S,[0,1]))$, 
$inj_N(x) \leq K_S$.  

Finally, for $x \in {\mathcal P}$, $inj_N(x) \leq \epsilon$.  Therefore for $x \in C(N)$,
$inj_N(x) \leq max \{K_S, L_{S_j} + 2 \epsilon\}$.  We can do this for all
$0 < \epsilon < \epsilon_3$, so that for $x \in C(N)$, $inj_N(x) \leq max \{K_S,
L_{S_j}$\}.  This completes the proof of Subcase A. \\

Case IIB.  The End $E$ is Geometrically Finite. \\

The basic outline of the argument for Case IIB is still the same.  We will divide
the convex core into a compact core and its complement.  However, to fill up the
convex core this time, we will use an interpolation between a boundary component of
the convex core and a simplicial hyperbolic surface, a continuous family of simplicial
hyperbolic surfaces, and projections of convex cores of covers with fewer accidental
parabolics.  The remaining portions will, again, lie in an
$\epsilon$-thin part, or lie in the projection of a cover which has fewer accidental
parabolics.

The first step in the proof of Case IIB will be to find a homotopy between the
boundary of the convex core and a useful or practical simplicial hyperbolic surface, and
a continuous family of simplicial hyperbolic surfaces which shrinks an accidental
parabolic associated to $B$. We will want points in the image of the homotopy and 
continuous family to lie in the convex core and to have bounded injectivity radius.

Let $f_1: S \rightarrow N$ be the pleated surface whose image is the boundary
component of $C(N)$ associated to the geometrically finite end $E$.  The following
result of Canary-Minsky guarantees the existence of a ``short'' homotopy
between a pleated surface and a simplicial hyperbolic surface.  

\begin{lemma}
\label{pleatshs} (Lem 6.1, Canary-Minsky \cite{can/min})  Given a pleated
surface $f_1: S \rightarrow N$ and $\epsilon > 0$, there exists a homotopy
$H_1: S \times [0,1] \rightarrow N$ such that $H_1|_{(S,0)} = f_1$, the map
$H_1|_{(S,1)} = f_2$ is a simplicial hyperbolic surface, and the trajectories
$H_1(x,[0,1])$ have lengths at most $\epsilon$.  Furthermore, the triangulation
$T$ of $S$ corresponding to $f_2$ contains a closed curve that is mapped to a closed
geodesic, and there is an upper bound ${\mathcal G}$ on the number of vertices in $T$,
depending only on the Euler characteristic of $S$.  
\end{lemma}

Without loss of generality, by pre-composing with the nearest point retraction map, we
can assume that
$H_1(S,[0,1]) \subset C(N)$.  As the image of a simplicial hyperbolic surface, by Lemma
\ref{shsinj}, for
$x \in H_1(S,1)$,
$inj_N(x)
\leq K_S$.  Then because the trajectories $H_1(x,[0,1])$ have lengths at most
$\epsilon$, for $x \in H_1(S,[0,1])$, $inj_N(x) \leq K_S + \epsilon$.   

We can use the following lemma of Canary-Minsky
to find a bounded length homotopy between the simplicial
hyperbolic surface
$f_2:S
\rightarrow N$ obtained above and a useful simplicial hyperbolic surface $f_3: S
\rightarrow N$.

\begin{lemma}
\label{shsuseful}  (Lem 6.2, Canary-Minsky \cite{can/min})  Let $f_2: S
\rightarrow N$ be a simplicial hyperbolic surface whose associated triangulation
$T_A$ has a closed curve $\hat{e}$ which is mapped to a closed geodesic and has at most
${\mathcal G}$ vertices.  Then there is a homotopy $H_2: S \times [1,2] \rightarrow N$
such that the following are true: \begin{enumerate} 
\item $H_2(S,[1,2]) \subset C(N)$, 
\item $H_2|_{(S,1)} = f_2$, 
\item $H_2|_{(S,2)} = f_3:
(S, T_A)
\rightarrow N$ is a useful simplicial hyperbolic surface with distinguished edge
$\hat{e}_1$, and 
\item $H_2(S,[1,2])$ is contained in the 
$[({\mathcal G} - 1)c]$-neighborhood of the image of $f_2$, for an independent constant
$c$. \end{enumerate} \end{lemma}

Because $H_2$ is a homotopy whose flow lines are contained in the $[({\mathcal G} -
1)c]$-neighborhood of the image of $f_2$, for
$x
\in H_2(S,[1,2])$,
$inj_N(x)
\leq K_S + A_S$, where the constant $A_S$ depends only on $S$.  Now we have bounded the
injectivity radius in a portion of $C(N)$.

The triangulation $T_A$ on $S$ from Lemma \ref{shsuseful} either contains an
edge which represents an accidental parabolic associated to $B$ or not.  In the latter
case, let $T_2$ be a triangulation of $S$ such that there exists an edge in
$T_2$ which is an accidental parabolic associated to $B$.   We can perform the modified
$HH$ sequence of elementary moves altering $T_A$ to $T_2$.  Let $T_B$ be the ``first''
triangulation in the sequence in which there is an edge $\gamma_B$ representing an
accidental parabolic associated to $B$.  

There are two cases:  the sequence of triangulations
has $Property \ T^B$ or not.  If the sequence has $Property \ T^B$, then using the
argument given in Case IIA(1), we can construct a continuous family of
simplicial hyperbolic surfaces $H_3(S,[2,3])$ and $H_4(S,[3,4))$ such that $H_3$
interpolates between the useful simplicial hyperbolic surface $f_3:(S,T_A) \rightarrow
N$ and a useful or practical simplicial hyperbolic surface $g:(S,T_B) \rightarrow N$,
and such that 
$H_4$ shrinks the accidental parabolic $\gamma_B$.  If the sequence does not have
$Property \ T^B$, then using the argument given in  Case IIA(2),
we can construct a continuous family of simplicial hyperbolic surfaces $H_3(S,[2,4))$
such that $H_3$ shrinks the accidental parabolic
$\gamma_B$ which is represented by an edge in $T_A$.

Let $W = H_1(S,[0,1]) \cup H_2(S,[1,2]) \cup H_3(S,[2,3]) \cup H_4(S,[3,4))$ (or
$W=H_1(S,[0,1]) \cup H_2(S,[1,2]) \cup H_3(S,[2,4))$ if the sequence does not have
$Property \ T^B$) be the union of the images of homotopies and continuous families of
simplicial hyperbolic surfaces.    Lemmas
\ref{edge}, \ref{elem move},
\ref{length to 0}, \ref{3 edge}, \ref{3 elem move}, \ref{3 length to 0}, and
\ref{shsuseful}
 guarantee that $W \subset C(N)$. 
 Furthermore, we also know that for any $x \in W$ that also lies in the image of a
simplicial hyperbolic surface, by Lemma
\ref{shsinj}, $inj_N(x) \leq K_S$ where the constant $K_S$ depends on the Euler
characteristic of $S$. To summarize, for $x \in W$, $inj_N(x) \leq max \{K_S + A_S, K_S
+ \epsilon\}$. Now we have bounded the injectivity radius in a portion of the
convex core.  Note that $\phi(W) \subset C(N)$ because $C(N)$ is convex, and that
by a result of Canary (Lem 4.4, \cite{canary1})  and Lemmas \ref{length to 0} and \ref{3
length to 0}, $\overline{\phi(W)}$ is compact. 

The next step in the proof of Subcase B involves finding a compact core
$R$ of $N - P_S$, where $P_S$ represents the collection of rank one cusps of $\NOE$
associated to $\partial S$, such that points in $R$ have bounded injectivity radius.

Let $\xi: S \rightarrow N$
denote the map $H_1|_{(S,1)}: S \rightarrow N$.  Let \{$S_{j}$\} be the components of
$S - \gamma_A$.  Let \{$M_j=\HHH / \xi_*(\pi_1(S_{j}))$\} be covers of $N$ with covering
maps \{$p_j$\}.  Let $P_A$ be the rank one cusp of $\NOE$ associated
to $\gamma_A$.  Let ${\mathcal P} = P_A \cup P_S$. 

Given $0 < \epsilon < \epsilon_3$.  By  Lemma \ref{compactincore}, there exists a
relative compact core
$R^*$ of
$\overline{C_{\epsilon}(N)} - {\mathcal P}$ such that $\partial C_\epsilon(N) \cap
C^\circ_\epsilon(N)
\subset
\partial R^*$ and such that $\overline{\phi(W)} \subset R^*$.  Let
\{$U_j^\prime$\} be the components of $[N - {\mathcal P}] - int \ R^*$
such that each $U_j^\prime$ is adjacent to
$P_A$, and let $U_k^\prime$ be the component of $[N - {\mathcal P}] -
int \ R^*$ which is a neighborhood of the geometrically finite end $E$.
Let
\{$U_j$\} be the components of $[C_{2 \epsilon}(N) - {\mathcal P}] - int \ R^*$ such
that each $U_j$ is adjacent to
$P_B$.  Then $U_j \subset U_j^\prime$. Let $U_k$ be the component of $[C_{2 \epsilon}(N)
- {\mathcal P}] - int
\ R^*$ such that $U_k \subset U_k^\prime$.

In the following lemma,  we will choose a compact core
$R$ of $N - P_S$ such that $R$ contains $R^*$,
$R$ lies in $C_{2 \epsilon}(N)$, and points in $R$ have bounded injectivity radius. 

\begin{lemma}
\label{surface core B}  For all $0 < \epsilon < \epsilon_3$, there exists a compact
core $R$ of $N - P_S$ bounded by surfaces $T_1$ and $T_2$ such that the following are
true: 
\begin{enumerate}
\item $R^* \subset R \subset C_{2 \epsilon}(N)$,
\item if $x \in R$, then $inj_N(x) \leq max \{K_S + A_S, K_S + 2 \epsilon, L_{S_j} + 2
\epsilon\}$,
\item for each $j$, $\partial (T_1 \cap U_j^\prime) \subset \partial {\mathcal P}$, and
\item for each $j$, there exist $\tilde{U}_j^\prime \subset M_j$ such that
$p_j|_{\tilde{U}_j^\prime}: \tilde{U}_j^\prime \rightarrow U_j^\prime$ is a
homeomorphism, and
$(p_j|_{\tilde{U}_j^\prime})^{-1}(T_1 \cap U_j^\prime) \subset \overline{C_{2
\epsilon}(M_j)}$. 
\end{enumerate}
\end{lemma}

\begin{proof}  The proof of this lemma mimics that of Lemma \ref{surface core}. 
We give an outline of the steps.  Let
${\mathcal T}_1 = H_4(S,4-\eta)$, where $H_4(S_j,4-\eta)$ lifts to $C_\epsilon(M_j)$, and
let ${\mathcal T}_2 = H_1(S,0)$.  By Lemma \ref{prodstructure}, we know that $U_k^\prime$, $U_k$,
and each $U_j^\prime$ and $U_j$ possesses a product structure.

Construct a surface $T_1^\prime = \bigcup \bar{\mathcal F}_j^\prime$ with the following
properties: \ (1) each $\bar{\mathcal F}_j^\prime$ is the union of embedded annuli
$A_{P_A,j}$ in $P_A \cap C_{2 \epsilon}(N)$, embedded infinite annuli $A_{P_S,j}$ in
$P_S 
\cap C_{2 \epsilon}(N)$, and a level surface
$\bar{F}_j$ in the product structure on $U_j$ such that ${H}_4(\bar{S}_j,1-\eta)$
lies to one side of $\bar{F}_j$; and (2) the surface
$T_1^\prime$ is properly homotopic to ${\mathcal T}_1$ via a proper homotopy $\zeta_1$
whose image lies in $\bigcup p_j(C_{2 \epsilon}(M_j))$.  Then for a point $x$ in the
image of
$\zeta_1$, $inj_N(x) \leq max \{L_{S_j} + 2 \epsilon\}$. 


Let $T_1 = T_1^\prime - (\bigcup A_{P_S,j})$ and $T_2 = {\mathcal T}_2 \cap (N - P_S)$. 
We can show that $T_1 \cap U_j = T_1 \cap U_j^\prime = \bar{F}_j$ and that $\partial
(T_1 \cap U_j) =
\partial
\bar{F}_j
\subset
\partial {\mathcal P}$.  Furthermore, we can show that there exists
$\tilde{U}_j^\prime
\subset M_j$ such that
$p_j|_{\tilde{U}_j^\prime}: \tilde{U}_j^\prime \rightarrow U_j^\prime$ is a
homeomorphism, and
$(p_j|_{\tilde{U}_j^\prime})^{-1}(\bar{F}_j) =
(p_j|_{\tilde{U}_j^\prime})^{-1}(\bar{F}_j)
\subset
\overline{C_{2 \epsilon}(M_j)}$,   


Then we can show that $T_1$ and $T_2$ bound a compact core $R$ of $N - P_S$ such that
$R^* \subset R \subset C_{2
\epsilon}(N)$.  By covering 
$R$ with $W$ and the image of $\zeta_1$, we can show that for
$x
\in R$,
$inj_N(x) \leq max \{K_S + A_S, K_S + 2 \epsilon, L_{S_j} + 2
\epsilon\}$.

Note that if the sequence does not have $Property \ T^B$, then we must replace
the maps $H_4$ and $\tilde{H}_4$ with $H_3$ and $\tilde{H}_3$, keeping the domains the
same. This completes the sketch of the proof of Lemma \ref{surface core B}. 
\end{proof}

Now we have bounded the injectivity radius for points in a compact core $R$ of $N -
P_S$ which contains $R^*$.

Consider again the  
surface constructed in the previous lemma 
$T_1 \cap U_j^\prime = \bar{F}_j \subset \partial R$ which was a $\pi_1$-injective,
separating, level surface in the product structure of $U_j$.  Let
$W_j^\prime$ be the closure of the component of
$(N - {\mathcal P}) -
\bar{F}_j$ such that $W_j^\prime \subset U_j^\prime$.  Let
$W_j^\prime$ be the closure of the component of
$(N - {\mathcal P}) -
\bar{F}_j$ such that $W_j^\prime \subset U_j^\prime$. 
Let $\tilde{W}_j^\prime$ be the subset of $\tilde{U}_j^\prime$ such that
$p_j|_{\tilde{W}_j^\prime}: \tilde{W}_j^\prime \rightarrow W_j^\prime$ is a
homeomorphism.  Then because $\partial
\bar{F}_j \subset
\partial {\mathcal P}$ and 
$(p_j|_{\tilde{W}_j^\prime})^{-1}(\bar{F}_j) =
(p_j|_{\tilde{U}_j^\prime})^{-1}(\bar{F}_j)
\subset \overline{C_{2 \epsilon}(M_j)}$, by Lemma
\ref{not in C(N)}, we know that  $(p|_{\tilde{W}_j^\prime})^{-1}[(N - C_{2
\epsilon}(N))
\cap W_j^\prime] = (M_j - C_{2 \epsilon}(M_j)) \cap \tilde{W}_j^\prime$.  

Recall that $H_1(S,0)$ is the component of $\partial C(N)$ associated to the
geometrically finite end $E$.  Note that we chose $H_1(S,0) \subset W \subset R^*$,
and that $U_k$ is a component of $[C_{2 \epsilon}(N) - {\mathcal P}] - int \ R^*$ not
adjacent to $P_A$.  So 
$U_k \cap C(N) = \emptyset$.

Suppose $x
\in [C(N) - {\mathcal P}] - R$.  Then $x
\in U_j
\cap C(N)$ for some
$j$.  
Because $(p|_{\tilde{W}_j^\prime})^{-1}[(N - C_{2
\epsilon}(N))
\cap W_j^\prime] = (M_j - C_{2 \epsilon}(M_j)) \cap \tilde{W}_j^\prime$, we know that if
$x \in W_j^\prime \cap C(N)$, then $x \in p_j(C(M_j))$.  So by the
induction hypothesis and Lemma
\ref{coverinj}, $inj_N(x) \leq L_{S_j}$.  

If $x \in {\mathcal P}$, then $inj_N(x) \leq \epsilon$.  Therefore for $x \in C(N)$,
$inj_N(x) \leq max \{K_S + A_S, K_S + 2 \epsilon, L_{S_j} + 2 \epsilon\}$.  We can do
this for all $0 < \epsilon < \epsilon_3$, so that for $x \in C(N)$, $inj_N(x) \leq max
\{K_S + A_S, L_{S_j}$\}.  

This completes the proof of Subcase B.  Thus, the proof of Theorem \ref{surface
bound} is complete.
\end{proof}

Now we will show that the injectivity radius in convex cores is bounded for all
hyperbolic 3-manifolds (weakly type-preserving) homotopy equivalent to a fixed compact
surface.

\begin{theorem}
\label{twisted surface bound} Let $S$ be a (possibly non-orientable) compact
hyperbolic surface. Then there exists a constant $K_S$ such that if $N$ is a hyperbolic
3-manifold such that
$h: S \rightarrow N$ induces a weakly type-preserving homotopy equivalence, 
and
$x
\in C(N)$, then
$inj_N(x) \leq K_S$.
\end{theorem}

\begin{proof}
If $S$ is orientable, then Theorem \ref{surface bound} proves the existence of $K_S$.

If $S$ is non-orientable, then let $\tilde{S}$ be an orientable double cover of $S$,
and let
$\tilde{N}$ be a double cover of $N$ associated to $\pi_1(\tilde{S})$. 
Using Theorem
\ref{surface bound}, there exists a constant $K$ such that for $x \in C(\tilde{N})$,
$inj_{\tilde{N}}(x) \leq K$.  

Let $\tilde{N} = \HHH / \tilde{\Gamma}$, and let
$N = \HHH / \Gamma$.  Then since $[\tilde{\Gamma}: \Gamma]$ is finite, we know that
the two limit sets $\Lambda_\Gamma$ and $\Lambda_{\tilde{\Gamma}}$ are equal. (II.K.15
\cite{maskit})  Thus,
$C(\tilde{N})$ is a double cover of $C(N)$.  Then, by Lemma
\ref{coverinj}, for $x \in C(N)$, $inj_N(x) \leq K$.
\end{proof}

The main result Theorem \ref{main result} follows from Theorem \ref{twisted surface
bound} since an $I$-bundle over a (possibly non-orientable) closed surface is homotopy
equivalent to a (possibly non-orientable) closed surface. (Thm 10.2, Hempel
\cite{hempel})

\bibliographystyle{amsplain}

\begin{thebibliography}{10}





%

\bibitem[Be]{beardon} A.\ Beardon, {\it The Geometry of Discrete Groups\/},
Graduate Texts in Mathematics 91, Springer-Verlag, 1983.


\bibitem[BP]{benedetti} R.\ Benedetti and C.\ Petronio, {\it Lectures on
Hyperbolic Geometry\/}, Universitext, Springer-Verlag, 1992.

\bibitem[Bi]{bielefeld} B.\ Bielefeld, ``Conformal Dynamics Problems Set,''
Institute of Mathematical Sciences---Stony Brook preprint.

\bibitem[Bo]{bonahon} F.\ Bonahon, ``Bouts des vari\'{e}t\'{e}s hyperboliques de
dimension 3,'' {\it Ann.\ of Math.\/} {\bf 124}(1986), 71--158.

\bibitem[C1]{canary1} R.D.\ Canary, ``A Covering Theorem for Hyperbolic
3-Manifolds,'' {\it Topology\/} {\bf 35}(1996), 751--778.

\bibitem[C2]{canary2} R.D.\ Canary, ``Ends of Hyperbolic 3-Manifolds,'' {\it
J.\ A.\ M.\ S.\/} {\bf 6}(1993), 1--35.
 
\bibitem[CEG]{can/ep/green} R.D.\ Canary, D.B.A.\ Epstein, and P.\ Green,
``Notes on notes of Thurston,'' in {\it Analytical and Geometrical Aspects of
Hyperbolic Spaces\/}, ed. by D.B.A.\ Epstein, Cambridge University
Press, 1987, 3--92.

\bibitem[CM]{can/min} R.D.\ Canary and Y.N.\ Minsky, ``On limits of tame
hyperbolic 3-manifolds,'' {\it J.\ Diff.\ Geom.\/} {\bf
43}(1996), 1--41.


\bibitem[EM]{epstein/marden} D.B.A.\ Epstein and A.\ Marden, ``Convex Hulls in
Hyperbolic Space, a Theorem of Sullivan, and Measured Pleated Surfaces,'' in
{\it Analytical and Geometrical Aspects of Hyperbolic Spaces\/}, ed. by D.B.A.\
Epstein, Cambridge University Press, 1987, 113--253.

\bibitem[F1]{fan1} C.\ Fan, ``Injectivity Radius Bounds in Hyperbolic Convex
Cores I,'' preprint.

\bibitem[F2]{fan2} C.\ Fan, ``Injectivity Radius Bounds in Hyperbolic Convex
Cores II,'' in preparation.

\bibitem[Har1]{harer1} J.L.\ Harer, ``Stability of the mapping class groups of
orientable surfaces,'' {\it Ann.\ of Math.\/} {\bf 121}(1985), 214--249.

\bibitem[Har2]{harer2} J.L.\ Harer, ``The virtual cohomological dimension of the
mapping class group of an orientable surface,'' {\it Invent.\ Math.\/}
{\bf 84}(1986), 157--176.

\bibitem[Hat]{hatcher} A.\ Hatcher, ``On triangulations of surfaces,'' {\it
Topology Appl.\/} {\bf 40}(1991), 189--194.

\bibitem[He]{hempel} J.\ Hempel, {\it 3-manifolds\/}, Annals of
Mathematics Studies 86, Princeton University Press, 1976.

%
%
%

\bibitem[KT]{kerck/thur} S.P.\ Kerckhoff and W.P.\ Thurston, ``Non-continuity of
the action of the modular group at Bers' boundary of Teichmuller space,''
{\it Invent.\ Math.\/} {\bf 100}(1990), 25--47.

\bibitem[KS]{kul/shalen} R.S.\ Kulkarni and P.B.\ Shalen, ``On Ahlfors'
Finiteness Theorem,'' {\it Adv.\ Math.\/} {\bf 76}(1989), 155--169.

\bibitem[L]{lloyd} N.G.\ Lloyd, {\it Degree Theory\/}, Cambridge Tracts in
Mathematics 73, Cambridge University Press, 1978.

\bibitem[Ma]{maskit} B.\ Maskit, {\it Kleinian Groups\/}, Grundlehren der
mathematischen Wissenschaften 287, Springer-Verlag, 1988.

\bibitem[McC]{mccullough} D.\ McCullough, ``Compact submanifolds of 3-manifolds
with boundary,'' {\it Quart.\ J.\ Math.\ Oxford Ser.\ (2)\/}
{\bf 37}(1986), 299--307.

\bibitem[MMS]{mmswarup} D.\ McCullough, A.\ Miller, and G.A.\ Swarup,
``Uniqueness of cores of non-compact 3-manifolds,'' {\it J.\ London Math.\
Soc.\/} {\bf 61}(1985), 548--556.

\bibitem[McM]{mcmullen} C.T.\ McMullen, {\it Renormalization and 3-Manifolds
which Fiber over a Circle\/}, Annals of
Mathematics Studies 142, Princeton University Press, 1996.

\bibitem[Mor]{morgan} J.W.\ Morgan, ``On Thurston's uniformization theorem for
three-dimensional manifolds,'' in {\it The Smith Conjecture\/}, ed. by J.\
Morgan and H.\ Bass, Academic Press, 1984, 37--125.
%
%
%
\bibitem[P]{perry} P.\ Perry, ``The Laplace operator on a hyperbolic manifold
II:  Eisenstein series and the scattering matrix,'' {\it J.\ Reine Angew.\
Math.\/} {\bf 398}(1989), 67--91.

\bibitem[Sc]{scott} G.P.\ Scott, ``Compact submanifolds of 3-manifolds,'' {\it
J.\ London Math.\
Soc.\/} {\bf 7}(1973), 246--250.

\bibitem[Su1]{sullivan} D.\ Sullivan, ``A finiteness theorem for cusps,'' {\it
Acta Math.\/} {\bf 147}(1981), 289--299.


\bibitem[Th1]{thurston} W.P.\ Thurston, {\it The Geometry and Topology of
3-manifolds\/}, lecture notes.

\bibitem[W]{wald} F.\ Waldhausen, ``On irreducible 3-manifolds which are
sufficiently large,'' {\it Ann. of Math.\/} {\bf 87}(1968), 56--88.

\end{thebibliography}

\end{document}